\documentclass[opre,nonblindrev]{informs3} % current default for manuscript submission

\usepackage{caption}
\usepackage{anyfontsize}
\usepackage{subfiles}
\usepackage{xr}
\usepackage{float}
\def\biblio{\bibliographystyle{plainnat}\bibliography{../DataDrivenRobustBibliography}}

\usepackage{accents}

\usepackage{natbib}
 \bibpunct[, ]{(}{)}{,}{a}{}{,}%
 \def\bibfont{\small}%
\TheoremsNumberedThrough     % Preferred (Theorem 1, Lemma 1, Theorem 2)
\ECRepeatTheorems

\EquationsNumberedThrough    % Default: (1), (2), ...
\RequirePackage[shortlabels]{enumitem}
\usepackage{mathtools} % for triangleq
\usepackage{enumitem} % for enumeration
\usepackage{dsfont} % for indicator function
\usepackage{mathrsfs}  

\newtheorem{innercustomgeneric}{\customgenericname}
\providecommand{\customgenericname}{}
\newcommand{\newcustomtheorem}[2]{%
	\newenvironment{#1}[1]
	{%
		\renewcommand\customgenericname{#2}%
		\renewcommand\theinnercustomgeneric{##1}%
		\innercustomgeneric
	}
	{\endinnercustomgeneric}
}

\newcustomtheorem{customthm}{Theorem}
\newcustomtheorem{customlemma}{Lemma}

\allowdisplaybreaks

\usepackage{epsfig}

\usepackage[ruled, vlined, linesnumbered]{algorithm2e}
\let\oldnl\nl% Store \nl in \oldnl
\newcommand{\nonl}{\renewcommand{\nl}{\let\nl\oldnl}}

\usepackage{IEEEtrantools}
\usepackage{multicol}
\usepackage{caption}
\usepackage{subcaption}
\usepackage[margin=2.54cm]{geometry}
\usepackage{lipsum}
\usepackage{bigints}
\usepackage{framed}
\usepackage{colortbl}
\usepackage{calligra}
\usepackage[export]{adjustbox}
\usepackage{booktabs}
\usepackage{multirow}
\usepackage{makecell}

 \usepackage{wrapfig}
\usepackage{algpseudocode}

\usepackage{bbm}
\usepackage[normalem]{ulem}
\usepackage[colorlinks=true,linkcolor=black,citecolor=black]{hyperref}
\usepackage{cleveref}

\usepackage{hyperref}
\hypersetup{colorlinks=false, pdfborder={0 0 0}}

\newcommand{\R}{{\mathcal R}}
\newcommand{\K}{{\mathcal K}}

\newcommand{\N}{{\mathbb N}}
\newcommand{\I}{{\mathbb I}}

\newcommand{\PP}{\mathbb P}
\renewcommand{\SS}{\mathbb S}

\newcommand{\Nhat}{{\hat{\mathbb N}}}

\newcommand{\Mhat}{{\hat{\mathbb M}}}

\newcommand{\mP}{\mathcal{P}}

\newcommand{\Q}{\mathbb Q}

\newcommand{\eps}{\epsilon}

\newcommand{\bm}{\boldsymbol}

\def\RR{ {\mathbb{R}}}

\newcommand{\X}{{\mathcal X}}

\newcommand{\DD}{{\mathscr D}}

\newcommand{\M}{{{\mathbb M}} }

\newcommand{\xii}{{\bm{\xi}}}

\newcommand{\bmu}{{\bm{\mu}}}
\newcommand{\bSigma}{{\bm{\Sigma}}}
\newcommand{\Sigmat}{{\hat{\bSigma}}}
\newcommand{\btheta}{{\bm{\theta}}}

\newcommand{\txi}{\tilde{\bm{\xi}}}

\newcommand{\x}{{\bm{x}}}
\newcommand{\z}{{\bm{z}}}

\newcommand{\s}{{\bm{s}}}
\newcommand{\phat}{{\hat{p}}}
\newcommand{\pl}{{\underline{p}}}
\newcommand{\lbar}{{\overline{\ell}}}
\newcommand{\fhat}{{\hat{f}}}
\newcommand{\fl}{{\underline{f}}}
\newcommand{\fu}{{\overline{f}}}
\newcommand{\fbar}{{\overline{f}}}
\newcommand{\xis}{{\xii|\s}}

\newcommand{\ts}{{\tilde{\bm{s}}}}

\newcommand{\bxi}{{\bm\xi}}

\newcommand{\EE}{{\mathbb{E}}}

\newcommand{\Wass}{\mathds W}

\let\Q\Rational
\let\R\Real
\let\Z\Integer

\def\st{\textup{ s.t. }}
\def\M{\mathbb{M}}
\def\Ms{\M_{\bxi|\s}}

\def\Mhats{\Mhat_{\bxi|\s}}

\def\tr{\textup{tr}}

\renewcommand{\qed}{{\hfill\halmos}}
\usepackage{graphicx}
\usepackage{tikz,ifthen}

   \usetikzlibrary{math}
   \usetikzlibrary{shapes.misc}
   \usetikzlibrary{shapes.multipart, positioning, patterns,backgrounds}

\usepackage{multirow}
\usepackage{array}
\usepackage{cleveref} % used for Cref
\usepackage{booktabs}
\usepackage{multirow}
\usepackage{makecell}

\usepackage{array}
\newcolumntype{x}[1]{>{\centering\arraybackslash}p{#1}}
\newcommand{\revise}[1]{{\color{blue}{#1}}}

\begin{document}

\def\biblio{}

%%%%%%%%%%%%%%%%

% Outcomment only when entries are known. Otherwise, leave as is and
%\setcounter{page}{1}
%\DOI{10.1287/xxxx.0000.0000}%

% \RUNAUTHOR{Yoon, Hanasusanto, Wang}

\RUNAUTHOR{Yoon et al.}
\RUNTITLE{Data-Driven Contextual Optimization with Gaussian Mixtures}

% \ARTICLEAUTHORS{YoungChul Yoon\thanks{Department of Industrial and Enterprise Systems Engineering, University of Illinois Urbana-Champaign, Urbana, IL 61801, USA. Email: {\tt ycyoon2, gah@illinois.edu}.},  
% Grani A. Hanasusanto\footnotemark[1], 
% and
% Yijie Wang\thanks{School of Economics and Management, 
% Tongji University, Shanghai, China. Email: {\tt yijiewang@tongji.edu.cn}.}
% }

\ARTICLEAUTHORS{YoungChul Yoon\thanks{Department of Industrial and Enterprise Systems Engineering, University of Illinois Urbana-Champaign, Urbana, IL 61801, USA. Email: {\tt ycyoon2, gah@illinois.edu}.},  Ling Dai\thanks{Department of Data Science, City University of Hong Kong, Hong Kong, China. Email: {\tt ling.dai@my.cityu.edu.hk, clint.ho@cityu.edu.hk}.}, Grani A. Hanasusanto\footnotemark[1], Chin Pang Ho\footnotemark[2], and Yijie Wang\thanks{School of Economics and Management, 
Tongji University, Shanghai, China. Email: {\tt yijiewang@tongji.edu.cn}.} 
}

% \TITLE{Bundling Information Goods of Decreasing Value}
%\TITLE{Nonparametric confidence intervals for the sample mean: complexity, exact algorithms, and deterministic approximations}
%\TITLE{Computation of exact bootstrap confidence intervals: complexity and deterministic algorithms}
%\TITLE{Sample Robust Optimization: Theory, Algorithms, and Applications}

\TITLE{Data-Driven Contextual Optimization with Gaussian Mixtures: Flow-Based Generalization, Robust Models, and Multistage Extensions}

\ABSTRACT{Contextual optimization enhances decision quality by leveraging side information to improve predictions of uncertain parameters. However, existing approaches face significant challenges when dealing with multimodal or mixtures of distributions. The inherent complexity of such structures often precludes an explicit functional relationship between the contextual information and the uncertain parameters, limiting the direct applicability of parametric models. Conversely, while non-parametric models offer greater representational flexibility, they are plagued by the ``curse of dimensionality," leading to unsatisfactory performance in high-dimensional problems. To address these challenges, this paper proposes a novel contextual optimization framework based on Gaussian Mixture Models (GMMs). This model naturally bridges the gap between parametric and non-parametric approaches, inheriting the favorable sample complexity of parametric models while retaining the expressiveness of non-parametric schemes. By employing normalizing flows, we further relax the GM assumption and extend our framework to more general distributions. Finally, inspired by the structural properties of GMMs, we design a novel GMM-based solution scheme for multistage stochastic optimization problems with Markovian uncertainty. This method exhibits significantly better sample complexity compared to traditional approaches, offering a powerful methodology for solving long-horizon, high-dimensional multistage problems. We demonstrate the effectiveness of our framework through extensive numerical experiments on a series of operations management problems. The results show that our proposed approach consistently outperforms state-of-the-art methods, underscoring its practical value for complex decision-making problems under uncertainty.
% By employing normalizing flows, we relax the GM assumption and extend our framework to arbitrary distributions. 
}%

%\HISTORY{This paper was first submitted in February 2018.}
% \HISTORY{This paper was first submitted on February 11, 2018 and underwent three revisions. It was accepted for publication on November 1, 2020. }

\maketitle

% \noindent\textbf{Code and Data Disclosure.}
% The AsCollected disclosure for this manuscript is available at:
% \url{https://ascollected.org/EP7_I76}.

\section{Introduction} \label{sec:intro}
Business decision-making and operations management are fundamentally intertwined with the presence of uncertainty. Across different industries, managers are constantly tasked with making critical choices whose outcomes depend on future events that cannot be perfectly predicted. Therefore, developing solution schemes that can effectively model and respond to such uncertainty is a critical research and practical endeavor. Classical stochastic optimization approaches solely focus on the probability distributions of the underlying random variables. A standard formulation of a stochastic optimization problem can be expressed as 
\[
\min_{\x \in \mathcal X} \EE_{\PP_\bxi}\left[\ell(\x,\txi)\right],
\]
where $\x$ represents the decision variables, $\mathcal X$ denotes its feasible region, and $\PP_\bxi$  characterizes the distribution of the random variable $\txi \in \RR^R$. While stochastic optimization provides a powerful framework for solving decision-making problems under uncertainty, this approach often overlooks valuable information that could provide a more accurate description of the uncertain problem. In many practical settings, decision-makers have access to observable exogenous factors, often referred to as contextual information or side information, that carry significant predictive power over the uncertain outcomes.  For instance, in inventory management, future demand for certain products is often correlated with observable covariates like past sales history, promotional activities,  or even macroeconomic trends~\citep{chen2008dynamic,zhang2020product}. Similarly, in energy systems, the electricity generation from renewable sources like wind or solar, as well as its demand, are heavily dependent on weather conditions and seasonal patterns~\citep{conejo2005locational,ward2013effect,bhatti2018making}. In revenue management, customer characteristics, browsing behavior on the e-commerce platform, and competitors' prices can provide valuable context for predicting willingness-to-pay and improving pricing strategies~\citep{chen2022statistical}. Furthermore, within financial applications like asset management, stock returns are known to be influenced by a wide array of macroeconomic indicators such as interest rates, inflation rates, and GDP growth, as well as firm-specific data like market capitalization and book-to-market ratios~\citep{fama2015five,gu2020empirical,leippold2022machine}. 

Leveraging this contextual information allows for a more refined understanding of the underlying uncertainties,  reflecting its conditional distribution rather than just its overall average behavior. This class of problem is known as contextual optimization. In this setting, prior to making a decision, the decision-maker observes a realization $\s$ of a vector of random exogenous covariates $\ts \in \RR^Q$. Subsequently, the objective is to minimize the expected loss conditioned on the observed contextual information $\ts = \s$. This gives rise to the following contextual stochastic optimization problem:
\begin{equation}\label{eq:contextual}
    \min_{\x \in \mathcal X} \EE_{\PP}\left[\ell(\x,\txi)|\ts=\s\right],
\end{equation}
where $\PP$ is the joint distribution of $\txi$ and $\ts$.

The core idea of contextual optimization is to leverage the observed side information $\s$ to inform the decision $\x$, recognizing the fact that the distribution of the uncertain parameter $\tilde \bxi$ may depend on it.  By minimizing the conditional expectation, the decision-maker obtains a policy that is optimally tailored to the specific context provided by $\s$.  Unfortunately, the true conditional distribution $\PP(\txi|\ts=\s)$ is rarely accessible to decision-makers in practice. Instead, they must rely on historical data $\{(\s_n,\bxi_n)\}_{n\in[N]}$ to learn this conditional distribution and then optimize for the conditional expectation. Existing approaches in the literature employ various techniques for this learning task.  Non-parametric methods, such as Nadaraya-Watson kernel regression~\citep{nadaraya:64,watson:64} or k-nearest neighbors~\citep{fix1985discriminatory}, offer flexibility without assuming a specific functional form for the conditional distribution. However, these methods can suffer from the ``curse of dimensionality," performing poorly when the dimension of $\s$ is high~\citep{srivastava2019robust,wang2024generalization}. On the other hand, parametric regression methods usually presume a functional relationship between the uncertain parameter and the contextual information, e.g., $\txi = g(\s) + \tilde{\bm\epsilon}$. 
While these methods can yield better convergence results, such strong structural assumptions may not hold in many real-world problems, and may struggle to capture complex multimodal data distributions.

Driven by the limitations of existing learning techniques and inspired by the inherent characteristics of real-world data in various operational contexts, we propose a novel approach that models the joint distribution of the uncertain parameter $\txi$ and the contextual information $\ts$ using a Gaussian Mixture Model (GMM). This modeling choice offers a compelling balance between the flexibility of non-parametric methods and the structural advantages of parametric models. In many operations management and statistical analysis problems, the Gaussian distribution serves as a cornerstone assumption due to its tractability as well as its empirical adequacy in capturing the characteristics of a wide range of real-world data. The GMM, composed of multiple Gaussian components, can be seen as a natural extension of the Gaussian distribution. At one extreme, the GMM with a single cluster reduces to the classical Gaussian model. At the other extreme, when the number of clusters reaches the number of data points, GMM reduces to the non-parametric kernel density estimation method~\citep{silverman:86}, effectively bridging the gap between parametric and non-parametric approaches. 

More importantly, GMMs are particularly adept at capturing complex, multi-modal data distributions and identifying underlying ``hidden states" or ``regimes" within the data, a feature highly relevant to contextual decision-making. This ability enables the model to adapt its understanding of uncertainty to different contexts based on the observation of the side information. For example, in the financial market, empirical studies have long observed that asset returns exhibit characteristics like skewness and multimodality, suggesting the presence of multiple underlying states or regimes~\citep{beedles1986asymmetry,fabozzi2005fat,schaller1997regime}. To address the deviation from the Gaussian distribution,~\cite{kon1984models} shows that a handful of Gaussian mixture components can accurately approximate the distribution of stock market returns.  More recently, ~\cite{twosigma} from Two Sigma study the integration of side information, such as interest rates and exchange rates, with financial market returns into a GMM framework. Their empirical analysis demonstrates that when contextual information is included, GMMs can effectively identify and model significant financial market regimes throughout history. 

A key advantage of adopting the GMM in contextual optimization is that it yields an analytical expression for the conditional distribution. Specifically, if $(\txi,\ts)$ follows a GM distribution $\M$, then the conditional distribution is also a GM distribution with parameters that can be explicitly derived from the parameters of $\M$. With this closed-form expression of the conditional distribution, the contextual optimization problem can be reformulated into a structured problem that is amenable to efficient solution.  Beyond this analytical tractability, the GMM framework fits nicely into contextual optimization due to its learning property. Contrary to the `curse of dimensionality' issue often encountered in non-parametric methods, learning GM distributions can sometimes become easier in higher dimensions~\citep{anderson2014more}. Intuitively, this is because distributions that are inseparable or overlapping in a low-dimensional space may become separable in a higher-dimensional space. This phenomenon aligns closely with the core idea of contextual optimization: when considering only the uncertain parameter $\txi$, the decision-maker might overlook important information to fully understand its behavior or depict its distribution. However, by introducing the contextual information $\ts$, the increased dimensionality provides a richer and more accurate characterization of the uncertainty, thereby leading to improved decision quality. 

Furthermore, the analytical properties of the conditional GM distributions enable our framework to tackle multistage stochastic optimization problems. In operations management, many critical problems, such as multistage energy planning and dynamic portfolio optimization, are sequential. In these settings, the underlying uncertainties often exhibit a Markovian structure, which naturally fits into the broader contextual optimization literature. Traditional solution schemes typically rely on sequential conditional Monte Carlo sampling, which suffers from severe computational bottlenecks because the size of the resulting scenario tree grows exponentially with the planning horizon $T$. Our GMM framework gracefully circumvents this issue by providing closed-form analytical expressions for the likelihood ratio weights, completely bypassing the need for intractable numerical integration. Thus, this structural property allows us to design a novel solution scheme that operates directly on observable historical data.

The major contribution of this paper can be summarized as follows:
\begin{enumerate}
    \item \textbf{A GMM framework for contextual optimization}: We propose a GMM framework for contextual optimization where the uncertain parameter and contextual information jointly follow a GM distribution. This framework offers a structured approach to leverage side information by modeling the conditional distribution in an analytical expression. Additionally, we provide a theoretical analysis of the approximation quality for the empirical GMM estimation. Our result demonstrates that the approximation error of the empirical GMM is only linearly dependent on its parameter estimation errors. Under reasonable assumptions, learning GMM parameters to a certain accuracy only requires polynomially many samples. This favorable sample complexity contrasts with the ``curse of dimensionality" often inherent in non-parametric methods, where achieving similar approximation quality typically requires a number of samples that grows exponentially with the dimension of the contextual information.
    \item \textbf{Handling general distributions via normalizing flows}: While our GMM framework provides strong modeling power, real-world data may originate from more complex, or even intractable distributions that cannot be perfectly captured by a pure GMM. To extend the applicability of our framework to such general data distributions, we propose a novel approach employing normalizing flows. Inspired by its representational ability for complex distributions, we train a normalizing flow to transform the original random vector into a latent space, where the transformed random vector follows a Gaussian or a GM distribution. This allows us to implicitly work with our GMM framework in the latent space, preserving analytical advantages for optimization even when the original data distribution is complex. To the best of our knowledge, this is the first work combining normalizing flows with the analytical tractability of GMMs for approximating general conditional distributions in stochastic optimization.
    \item \textbf{Distributionally robust optimization}: To mitigate the overfitting issue from empirical GMM estimations and improve decision reliability, we integrate our framework with Distributionally Robust Optimization (DRO). Leveraging the analytical availability of our GMM framework, we construct a Wasserstein ambiguity set centered around the empirical conditional distribution. This formulation avoids the degeneracy issues encountered by many conditional DRO models when they condition on a measure-zero singleton $\{\s\}$. We further derive distribution coverage results and establish a performance guarantee for the DRO model. Finally, we extend our DRO framework to handle situations where the true number of components in the underlying GMM is unknown.  By considering an expanded ambiguity set that encompasses candidate conditional GM distributions with different plausible cluster numbers, our framework offers performance guarantees against potential misspecification of the model complexity.
     \item \textbf{Multistage setting}: Building upon the structural insights provided by our GMM framework for conditional distributions, we further develop a novel approximation scheme tailored for multistage stochastic optimization problems with Markovian uncertainty. Existing solution schemes mainly rely on Sample Average Approximation (SAA) or kernel regression to handle the evolving conditional distributions. However, the SAA method involves sequential conditional Monte Carlo sampling whose problem size grows exponentially with the time horizon $T$, while kernel regression methods also incur a sample complexity that grows exponentially with the dimension of the uncertain parameter.  In contrast, our proposed GMM-based approximation scheme offers an attractive alternative whose sample complexity grows only polynomially with the time horizon $T$ and does not suffer from the curse of dimensionality issue, rendering the effectiveness of our approach for large-scale high-dimensional multistage problems.
\end{enumerate}

\subsection{Literature Review}

\textbf{Contextual optimization}:
The increasing availability of data correlated with uncertain outcomes has spurred significant research into contextual optimization. As formalized in~\eqref{eq:contextual}, contextual optimization problems aim to optimize a decision based on the conditional distribution of uncertain parameters given observed side information.  As highlighted in the survey by~\cite{sadana2025survey}, the literature on solving data-driven contextual optimization problems can be broadly categorized into three paradigms: decision rule approach, sequential learning and optimization, and integrated learning and optimization. We review the relevant literature within these frameworks.

The decision rule approach directly seeks a policy from the covariate $\s$ to the decision $\x$ by minimizing an empirical risk or a related objective over the historical data. Various functional forms have been explored for these decision rules. The linear decision rule approach, favored for its interpretability and tractability, was studied by~\cite{ban2019big} in newsvendor problems. Unfortunately, linear decision rules may not achieve asymptotic optimality for general problems. To obtain greater flexibility,~\cite{bertsimas2022data} propose to approximate the optimal policy with a linear policy within the reproducing kernel Hilbert space (RKHS). More complex non-linear decision rules have also been implemented, including decision tree~\citep{bertsimas2019optimal} and neural networks~\citep{zhang2017assessing, huber2019data, oroojlooyjadid2020applying}. Finally, recognizing that minimizing empirical risk can be prone to overfitting, particularly with limited data or complex functional classes, distributionally robust variants have been designed to offer robustness against data uncertainty. While the decision rule approach is very computationally efficient at the decision stage, it still faces several challenges. First, the decision rule method typically restricts the potential solution space to the chosen class of decision rules. Additionally, ensuring that the learned decision rule outputs decisions can be difficult. Standard empirical risk minimization does not inherently enforce feasibility, and various problem-specific techniques are needed to project the rule's output onto the feasible region~\citep{zhang2021universal,chen2023end}. 

The sequential learning and optimization paradigm, also known as predict-then-optimize or prescriptive analytics, follows a two-stage process~\citep{bertsimas2020predictive}.  In the first stage, a model is trained to predict the conditional distribution or sufficient statistics (e.g., conditional mean and covariance) given the context, while in the second stage, the decision-maker simply solves an optimization problem based on the predicted distribution or statistics. There exist various solution schemes for the prediction stage. Non-parametric methods such as Nadaraya-Watson kernel regression~\citep{nadaraya:64,watson:64} and k-nearest neighbors~\citep{fix1985discriminatory} offer flexibility without strong distributional assumptions. However, they can suffer from the curse of dimensionality, performing poorly when the dimension of the contextual information is high~\citep{srivastava2019robust, wang2024generalization}. Parametric regression models, by assuming a specific functional form for the relationship between the contextual information $\s$ and $\bm \xi$, achieve a faster convergence rate~\citep{sen2018learning, kannan2025data}. However, such a functional form assumption may not hold in many practical problems and may struggle to learn complex multimodal data distributions.  A key challenge in the sequential learning and optimization framework is the ``optimizer's curse," where small errors in the prediction model can lead to significantly suboptimal decisions in the downstream optimization problem. To mitigate these issues, robust optimization, distributionally robust optimization, and variance regularization techniques have been applied to enhance models' performance in the out-of-sample circumstances~\citep{ bertsimas2017bootstrap,chenreddy2022data, kannan2024residuals,wang2026data}. Our work belongs to this category and contributes to the literature by its ability to handle complicated high-dimensional multimodal distributions.

The integrated learning and optimization framework follows an end-to-end approach that trains the prediction model based on the downstream task loss~\citep{donti2017task, elmachtoub2022smart,qi2023practical}. Compared with using a decision-blind loss such as mean squared error (MSE), the end-to-end model can yield a prediction that aligns better with the decision stage.  For linear programs, the gradient of the decision vector with respect to the predicted cost vector is either nonexistent or zero. To address this issue, surrogate loss functions and other smoothing techniques are employed to create differentiable proxies for the optimization objective~\citep{wilder2019melding,blondel2020learning,elmachtoub2022smart, huang2024decision}. While conceptually appealing for its potential to yield high-quality decisions by directly minimizing over the task loss, end-to-end models can be significantly more computationally demanding and complex to implement compared with the previous two paradigms, particularly for large-scale or combinatorial optimization problems~\citep{tang2024pyepo}.

We remark that contextual optimization is a rapidly growing area, and it is impossible to review all related studies within this section. Therefore, we refer readers who are interested in this topic to a comprehensive review by~\cite{sadana2025survey}.

\noindent \textbf{Normalizing flows}: Normalizing flow is a powerful class of generative models that allow for both efficient sampling and exact likelihood evaluation~\citep{dinh2014nice,rezende2015variational,kobyzev2020normalizing,papamakarios2021normalizing}. A normalizing flow defines a complex probability distribution by applying a sequence of invertible and differentiable transformations (diffeomorphisms) to a simple base distribution (e.g., Gaussian). The probability density of a sample then can be computed using the change of variables formula, which involves the density of the transformed sample under the base distribution and the determinant of the Jacobian of the inverse transformation~\citep{villani2008optimal}. Various normalizing flow architectures have been developed with different properties and objectives. Early examples include simple elementwise flows and linear flows~\citep{dinh2014nice}. More expressive non-linear transformations were introduced in planar and radial flows~\citep{rezende2015variational}. Borrowing ideas from residual networks, residual flows offer another way to construct invertible mappings, where the residual connection can be viewed as a discretization of a first-order ordinary differential equation~\citep{chen2019residual}. Continuous flows move a step forward by directly learning the continuous dynamical system, which is also known as infinitesimal flows~\citep{chen2018neural, grathwohl2018ffjord}. 

Among the most successful and widely adopted normalizing flow architectures are coupling flows and autoregressive flows, which maintain a balance between representation power and computational efficiency~\citep{kingma2016improved,papamakarios2017masked}.  Coupling flows partition the input and apply transformations to one part conditioned on the other~\citep{dinh2014nice,dinh2016density}. Autoregressive flows, which are particularly relevant to our work, structure the transformation such that each output dimension depends only on its previous dimensions in a fixed ordering~\citep{kingma2016improved,papamakarios2017masked,huang2018neural}. This enables the projection of the contextual information $\s$ into the latent space without the realization of uncertain parameters $\bm\xi$. Additionally, this ordering structure results in a triangular Jacobian matrix, making its determinant easily computable. Finally, the universality property of autoregressive flows has been proven, indicating their ability to learn any target probability density~\citep{huang2018neural,jaini2019sum}. Prominent examples include masked autoregressive flows~\citep{papamakarios2017masked}, neural autoregressive flows~\citep{huang2018neural}, and inverse autoregressive flow~\citep{kingma2016improved}.

\noindent \textbf{Multistage stochastic optimization problem with Markovian uncertainty}: Multistage stochastic programming (MSP) problems involve a sequence of decisions made over time in the face of evolving uncertainty~\citep{shapiro2021lectures}. These problems are significantly more complex than the single or two-stage problems due to the need to determine policies that are functions of the revealed uncertain parameters at each stage~\citep{birge2011introduction}.  We are particularly interested in a class of MSP problems where the uncertain process follows a Markovian structure, i.e., the distribution of future uncertainty depends only on the current state~\citep{kallenberg1997foundations}. This Markovian property directly fits into the setting of contextual optimization and is prevalent in many real-world applications, such as financial modeling and inventory control. 

A standard approach for solving MSP problems is the SAA method~\citep{kleywegt2002sample}. It assumes the knowledge of the true distribution, and employs a discretization scheme to approximate the MSP problem by a scenario tree generated by conditional Monte Carlo Sampling. However, the size of the scenario tree grows exponentially as the planning horizon increases, making the problem computationally intensive~\citep{shapiro:03,reaiche2016note,jiang2021complexity,shapiro2021lectures}. To alleviate this issue, kernel regression methods have also been designed to solve MSP problems with Markovian uncertainty~\citep{park2022data}. Although the problem size of kernel regression methods does not grow exponentially with $T$, it suffers from the curse of dimensionality with respect to the dimension of the uncertain parameter. Hence, there is merit in designing a new solution scheme that offers a better sample complexity. 

\subsection{Paper Structure and Notation} \label{subsec:notation}
The remainder of the paper is organized as follows. In Section~\ref{sec:contextual_GMM}, we formally introduce our Gaussian Mixture Model framework for contextual optimization, providing a theoretical analysis of the approximation quality. We further present a novel approach leveraging normalizing flows to extend the applicability of our framework to general data distributions that may not strictly follow a GMM. Section~\ref{sec:DRO} proposes to address parameter estimation uncertainty and potential misspecification of the cluster number using distributionally robust optimization. Section~\ref{sec:observational} extends our approach to MSP problems under Markovian uncertainty, establishing a novel GMM-based approximation scheme. We present extensive numerical experiments on synthetic and real-world datasets to demonstrate the practical effectiveness of our proposed approach in Section~\ref{sec:experiment}.

\noindent \paragraph{\textbf{Notation}} We use boldface letters to denote vectors and matrices. Random variables are indicated with a tilde (e.g., $\txi$), while their realizations appear without the tilde (e.g., $\bxi$). The multivariate Gaussian distribution with mean vector $\bmu$ and covariance matrix $\bSigma$  is denoted by $\N(\bmu,\bSigma)$, and its density function is written as $\mathcal N(\bm z|\bmu,\bSigma)$. We denote the vector of all ones by $\mathbf{e}$. The probability simplex in $\RR^K_+$ is denoted as $\Delta^K$. The little-o notation $o(\cdot)$ is used for asymptotic analysis: $h(\epsilon)=o(\epsilon)$ if $h(\epsilon)/\epsilon\to 0$, that is, $h(\epsilon)$ becomes negligible relative to $\epsilon$ as $\epsilon\to 0$. 

\section{Contextual Stochastic Optimization with Mixture Models}\label{sec:contextual_GMM}

\begin{figure}
    \centering
    \includegraphics[width=0.95\linewidth]{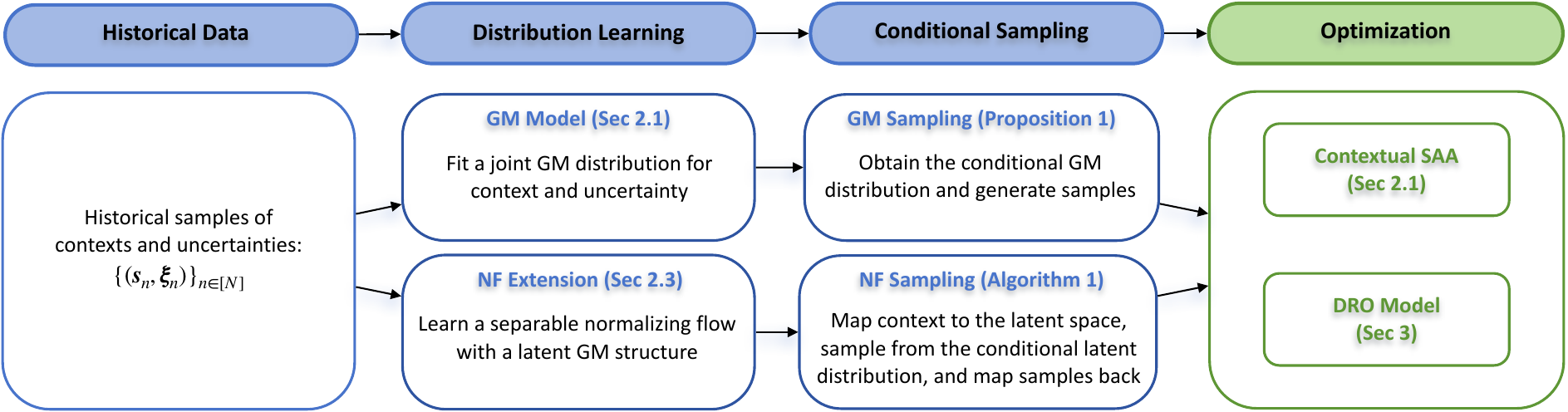}
    \caption{An overview of the GMM contextual optimization framework for static optimization problems.}
    \label{fig:structure_static}
\end{figure}

\subsection{Contextual Stochastic Optimization with GMMs}
We consider the contextual stochastic optimization problem given by
\begin{equation}
\label{eq:cso}
\min_{\x\in\X} \EE_{\M}\left[\ell(\x,\txi)\big|\ts=\s\right].
\end{equation}
The random vector $(\ts,\txi)\in\RR^{Q+R}$ follows a joint distribution 
$$\M\coloneqq\sum_{k\in[K]}p^{k}\N(\bmu^{k},\bSigma^{k}),$$ which is a Gaussian mixture (GM) distribution  with component means and covariances given by
\begin{equation*}
\bmu^k\coloneqq\begin{bmatrix}\bmu^k_\s\\\bmu^k_\bxi\end{bmatrix}\in\RR^{Q+R}\quad\textup{and}\quad
\bSigma^k\coloneqq\begin{bmatrix}\bSigma^k_{\s\s} & \bSigma^k_{\s \bxi} \\\bSigma^k_{\bxi\s} & \bSigma^k_{\bxi\bxi}\end{bmatrix}\in\SS_+^{Q+R}\qquad\forall k\in[K].
\end{equation*}
The mixture weights are given by $\bm p\in\Delta^K$, where $p^k\in[0,1]$ represents the weight of the $k$-th component. 
At first glance, problem \eqref{eq:cso} may appear challenging because it involves computing the conditional expectation of $\ell(\x,\txi)$ given that $\tilde\s=\s$ is observed. However, a recent result establishes that the conditional distribution of a GM distribution remains a GM distribution, thereby enabling the development of a principled solution scheme. 
\begin{proposition}
\label{prop:conditional_param}
The conditional distribution  of $\txi$ given $\ts=\s$ is a GM 
\begin{equation*}
\Ms\coloneqq\sum_{k\in[K]}p^{k}_{\xii|\s} \N(\bmu_{\xii|\s}^{k},\bSigma_{\xii|\s}^{k}),
\end{equation*}
 where
\begin{equation}\label{eq:conditional_param}
\begin{array}{rl}
\displaystyle\bm{\mu}_{\xii|\s}^k &\displaystyle = \;\; \bm{\mu}^k_{\xii} + \bm{\Sigma}_{\xii\s}^k( \bm{\Sigma}_{\s\s}^{k})^{-1}(\bm{s}-\bm{\mu}_\s^k),\\
\displaystyle\bm{\Sigma}^k_{\xii|\s} &\displaystyle = \;\; \bm{\Sigma}_{\bxi\bxi}^k- \bm{\Sigma}_{\xii\s}^k( \bm{\Sigma}_{\s\s}^{k})^{-1} \bm{\Sigma}_{\s\xii}^k, \ \textup{and}\\
\displaystyle p^k_{\xii|\s} &\displaystyle = \;\; \frac{ p^k \mathcal{N} \left(\bm{\s}| \bm{\mu}_{\s}^k, \bm{\Sigma}_{\s\s}^k\right)}{\sum_{j=1}^K p^j \mathcal{N} \left(\bm{s}| \bm{\mu}_{s}^{j}, \bm{\Sigma}_{\s\s}^{j}\right)}.
\end{array}
\end{equation}
\end{proposition}
This result simplifies problem \eqref{eq:cso} into an \emph{un}conditional stochastic optimization problem
\begin{equation}\label{eq:unconditional}
\min_{\x\in\X} \EE_{\Ms}\left[\ell(\x,\txi)\right].
\end{equation}
For specific loss functions, such as the mean-variance utility used in portfolio optimization, problem~\eqref{eq:unconditional} can be solved directly by plugging the conditional parameters. For more general loss functions, problem~\eqref{eq:unconditional} is amenable to a high-quality solution via the sample average approximation (SAA). Since the exact conditional distribution $\Ms$ is still a GM distribution, the decision maker can efficiently draw arbitrarily many conditional samples from it. Then, the expectation is replaced with the average of the loss function over generated samples.

In practice, however, the joint distribution $\M$ is rarely available to the decision maker. Instead, one must rely on $N$ historical observations $\{(\s_n,\xii_n)\}_{n\in[N]}$ to estimate these parameters. Under this data-driven setting, a standard approach is to fit a GM distribution $\Mhat$ to these observations, and solve the approximate contextual stochastic optimization problem 
\begin{equation*}
\min_{\x\in\X} \EE_{\Mhats}\left[\ell(\x,\txi)\right], 
\end{equation*}
where the approximate conditional distribution $\Mhats$ is obtained by applying Proposition \ref{prop:conditional_param} to $\Mhat$. While this empirical approach is straightforward to implement, the significance of its approximation error is still not fully understood. In the following section,  we provide a theoretical analysis of this empirical approach to quantify its approximation quality.

\subsection{Approximation Quality}
We analyze the approximation quality of this empirical approximation. Quantifying this approximation error is crucial for understanding the reliability of the empirical approach. To facilitate our analysis, we make the following standard assumptions:
\begin{enumerate}[(A)]
\item \label{as:bounded_loss} The loss function is bounded, i.e.,  $-\lbar\leq\ell(\x,\xii)\leq\lbar$ for some $\lbar\in\RR_+$.
\item \label{as:covariance} There exists a constant $\gamma\in\RR_+$ such that $\|\bmu^k\|,\|\hat\bmu^k\|\leq \gamma$. In addition, the covariance matrices $\bSigma^k$ and $\hat\bSigma^k$ are positive definite, and there exist constants $\alpha,\beta\in\RR_{++}$ such that $\alpha\I\preceq\bm\Sigma^k,\hat\bSigma^k\preceq\beta\I$ for all $k\in[K]$.
\item \label{as:weights} There exists a constant $\pl\in\RR_{++}$ such that $p^k\geq \pl$ for all $k\in[K]$. 
\end{enumerate}
The requirement that the loss function is bounded in Assumption \ref{as:bounded_loss} is made to simplify the exposition. 
The analysis in the paper can be extended to losses with suitable growth conditions. In particular, finite GM distributions have Gaussian tails, so for losses with at most polynomial growth, one can localize the analysis to a high-probability ball and control the complement using Gaussian tail bounds. This yields analogous high-probability bounds with additional tail terms, but the resulting constants are less transparent. We maintain the bounded-loss assumption to keep the presentation concise.
Assumption~\ref{as:covariance} imposes standard boundedness and nondegeneracy conditions on the GM component means and covariance matrices. The lower eigenvalue bound on the covariance matrices ensures that the Gaussian components admit well-defined densities on the ambient space and that the conditional distributions remain well conditioned. The upper eigenvalue and mean bounds prevent the components from becoming excessively diffuse or drifting arbitrarily far. These conditions ensure that the constants in the perturbation analysis below remain finite.
The lower bound \(\pl\) in Assumption~\ref{as:weights} prevents mixture components from having arbitrarily small weights, which is important for stable learning and for controlling the conditional mixture weights. This is consistent with the GMM learning literature, where small mixture weights are known to worsen the conditioning and sample complexity of the learning problem~\citep{kalai2012disentangling}.

\begin{theorem}
\label{thm:approx_error}
Let $\M=\sum_{k\in[K]}p^{k}\N(\bmu^{k},\bSigma^{k})$ be the true underlying Gaussian mixture distribution and $\Mhat=\sum_{k\in[K]}\hat p^k\Nhat(\hat\bmu^k,\hat\bSigma^k)$ be the estimated one. If $|p^{k}-\hat p^k|\leq \epsilon_p$, $\|\bmu^{k}-\hat\bmu^k\|\leq\epsilon_\bmu$, and $\|\bSigma^{k}-\hat\bSigma^k\|\leq\epsilon_\bSigma$ for all $k\in[K]$, then 
\begin{equation}
\label{eq:approx_error}
\begin{array}{rl}
&\left|\EE_{\Ms}\left[\ell(\x,\txi)\right]- \EE_{\Mhats}\left[\ell(\x,\txi)\right]\right|
\leq\lbar \left(2\left(C_{Q}\| \s\|^2+C'_{Q}\| \s\| + C''_{Q}\right)+      c_R\|\s\| + c'_R\right)\quad{\forall \x\in\X}, 
\end{array}
\end{equation}
where
  $C_Q=\left(1 + \frac{\epsilon_p}{\underline{p}} \right)  \frac{Q \epsilon_\bSigma}{2 \alpha^2}$, $C'_Q=\left(1 + \frac{\epsilon_p}{\underline{p}} \right)  \left( \frac{\epsilon_\bmu}{\alpha} + \frac{Q \epsilon_\bSigma  \gamma}{\alpha^2} \right)$, and $C''_Q=\frac{\epsilon_p}{\underline{p}} 
+ \left(1 + \frac{\epsilon_p}{\underline{p}} \right)
\left( \frac{\epsilon_\bmu}{\alpha} \, \gamma 
+ \frac{Q \epsilon_\bSigma}{2 \alpha^2} (\gamma^2 + \alpha) 
+ o(\epsilon_\mu + \epsilon_\bSigma) \right)$,  $c_R=\frac{\sqrt{R}(\alpha + \beta)}{\alpha^3}\epsilon_\bSigma$, and $c'_R=\frac{\sqrt{R}}{\alpha} \left( \left( \frac{\beta}{\alpha} + 1 \right) \epsilon_\bmu 
+ \frac{\alpha + \beta}{\alpha^2} \gamma \epsilon_\bSigma \right) 
+ \frac{1}{2}  \frac{R^2}{\alpha} \left( \frac{\beta}{\alpha} \right)^2 \epsilon_\bSigma$. 
\end{theorem}

\begin{remark}
\label{rem:decision_performance}
Theorem~\ref{thm:approx_error} directly implies a bound on the out-of-sample performance of the decision obtained from the estimated conditional GM distribution. Let
\[
\widehat \x \in \arg\min_{\x\in\X}
\EE_{\Mhats}\left[\ell(\x,\txi)\right]
\qquad\text{and}\qquad
\x^\star \in \arg\min_{\x\in\X}
\EE_{\Ms}\left[\ell(\x,\txi)\right].
\]
Denote the right-hand side of \eqref{eq:approx_error} by
$\varepsilon
\coloneqq
\lbar \left(2\left(C_{Q}\| \s\|^2+C'_{Q}\| \s\| + C''_{Q}\right)
+ c_R\|\s\| + c'_R\right)$.
Then
\[
\EE_{\Ms}\left[\ell(\widehat\x,\txi)\right]
-
\EE_{\Ms}\left[\ell(\x^\star,\txi)\right]
\leq 2\varepsilon.
\]
Indeed, by the uniform bound in Theorem~\ref{thm:approx_error} and the optimality of \(\widehat \x\) under \(\Mhats\), we have
\[
\begin{aligned}
\EE_{\Ms}\left[\ell(\widehat\x,\txi)\right]
&\leq
\EE_{\Mhats}\left[\ell(\widehat\x,\txi)\right]+\varepsilon(\s) \\
&\leq
\EE_{\Mhats}\left[\ell(\x^\star,\txi)\right]+\varepsilon(\s) \\
&\leq
\EE_{\Ms}\left[\ell(\x^\star,\txi)\right]+2\varepsilon(\s).
\end{aligned}
\]
Thus, the conditional-expectation approximation bound in Theorem~\ref{thm:approx_error} directly translates estimation errors in the fitted GMM parameters into a decision-performance guarantee. 
\end{remark}

The bound accounts for the errors arising from the estimation of mixture probabilities, given by $2\lbar\left(C_{Q}\| \s\|^2+C'_{Q}\| \s\| + C''_{Q}\right)$, and from the estimation of the mean and covariance of each Gaussian component, captured by $\lbar(c_R\|\s\|+c_R')$. %{\color{red} The error from the latter is dominant with cubic dependence on $\alpha$}  
We observe that the approximation quality deteriorates as ${\alpha}$ decreases and $\gamma$ increases. These constants are introduced to conservatively upper bound $ \|(\bSigma_{\s\s}^k)^{-1}\|\leq \frac{1}{\alpha}$ and $\|\s-\bmu_{s}^k\|\leq \|\s\|+\gamma$.   Thus, when there is no Gaussian component in which the covariates $\tilde{\s}$ exhibit high variability or whose mean is close to $\s$, the approximation error becomes significant. This corresponds to cases where samples $\txi$ in the vicinity of $\tilde{\s}=\s$ are rarely observed, which slows the learning rate of the conditional distribution.  The dependence in $\|\s\|$ stems from the same reasoning: since the GM distribution is sub-Gaussian, the density $f(\s)$ decreases as the covariates $\s$ move further from the origin. 
The bound also depends on $\beta$, with the proof showing that the error increases as $\|\bm{\Sigma}_{\xii\s}^k\|\leq \beta$ grows. This reflects strong dependence between $\tilde\s$ and $\txi$, which complicates the estimation of the conditional distribution. In contrast, in the limiting case of zero correlation, the marginal distribution of $\txi$ alone suffices, making learning easier.

Finally, this result confirms that the error introduced by the data-driven approximation is linear in $\epsilon_p$, $\epsilon_\bmu$, and $\epsilon_\bSigma$. This linear dependency is particularly significant when considering the sample complexity of GMM learning. Under standard regularity conditions, GMM algorithms typically yield parameter estimation errors that decay polynomially in the sample size $N$ and the dimension $Q+R$ (See Remark~\ref{remark:GMM} for details). This stands in contrast to existing data-driven methods for general contextual stochastic optimization, which often suffer from approximation errors that decay slowly with the dimension $Q$ of the contextual covariates $\ts$ \citep{wang2024generalization,wang2026data}.  Moreover, in the special case of a single Gaussian component ($K=1$), standard concentration results imply that $\epsilon_\bmu$ and $\epsilon_\bSigma$ are of the order $O(\frac{1}{\sqrt{N}})$. While our framework requires learning a GMM, the primary focus of this work is contextual optimization rather than the development of new GMM estimation algorithms. We therefore build on established GMM learning methods and their statistical guarantees, which are summarized in the following remark.

\begin{remark}\label{remark:GMM}
Learning mixtures of Gaussians is a highly active research area, where numerous algorithms have been developed, accompanied by theoretical results that establish their sample complexity. The work  \citep{moitra2010settling,kalai2012disentangling} develops an algorithm that, for any GM distribution $\M=\sum_{k\in[K]}p^{k}\N(\bmu^{k},\bSigma^{k})$, outputs an $\epsilon$-estimate $\Mhat=\sum_{k\in[K]}p^{k}\Nhat(\hat\bmu^{k},\hat\bSigma^{k})$ satisfying $|p_k-\hat p_k|\leq\epsilon$ and $\mathds{TV}(\N(\bmu^{k},\bSigma^{k}),\Nhat(\hat\bmu^{k},\hat\bSigma^{k}))\leq \epsilon$ with high probability.  Here, the total variation distance is used to determine the closeness of the component Gaussians, which is defined as $\mathds{TV}(f,\fhat)=\frac{1}{2}\int |f(\bxi)-\hat f(\bxi)|{\mathrm d}\bxi$ for any two distributions $f$ and $\hat f$. In the case of two Gaussians $\N(\bmu^{k},\bSigma^{k})$ and $\Nhat(\hat\bmu^{k},\hat\bSigma^{k}))$, a small total variation between them implies that the parameters $\hat\bmu$ and $\hat\bSigma$ have small estimation errors as well (see Proposition \ref{prop:tv_imply_parameters}). The running time and data requirement of this algorithm are polynomial in the dimension $D$, accuracy ${\epsilon}$, and condition number $\kappa=\min\{\pl,\min_{i,j\in[K]:i\neq j}\mathds{TV}(\N(\bmu^{i},\bSigma^{i}),\N(\bmu^{j},\bSigma^{j}))\}^{-1}$. 
% \yj{More recently, \cite{ashtiani2018nearly} established nearly tight sample complexity bounds, showing that $\tilde{O}(K(Q+R)^2/\epsilon^2)$ samples are sufficient and necessary to learn a K-component GM distribution up to error $\epsilon$ in total variation distance. }
Hence, the number of samples $N$ required to sustain an approximation error of at most $\tau$ in \eqref{eq:approx_error} will be polynomial in the desired accuracy $\tau$ and problem dimensions $Q$ and $R$. 
More recent literature provides tighter bounds for various problem settings and algorithms. For instance, when the number of components $K=2$, ~\cite{hardt2015tight} proposes an algorithm that can learn an arbitrary GM distribution to a certain accuracy with polynomially many samples. For a spherical GM distribution that is $\Omega(\sqrt{\log K})$ separable, ~\cite{kwon2020algorithm} shows that $\tilde{O}(K(Q+R)/\epsilon^2)$ samples suffice for the Expectation Maximization (EM) algorithm to achieve $\epsilon$ parameter estimation accuracy. Using the PCA learning method, ~\cite{ashtiani2018sample} derive a sample-efficient learning algorithm that can learn an arbitrary GM distribution to certain TV accuracy with $\tilde{O}(K(Q+R)^2/\epsilon^4)$ samples, and an axis-aligned GM distribution with $\tilde{O}(K(Q+R)/\epsilon^4)$ samples. This result is later improved by their subsequent work~\citep{ashtiani2018nearly} with a sample complexity of $\tilde{O}(K(Q+R)^2/\epsilon^2)$ for general GM distributions, and $\tilde{O}(K(Q+R)/\epsilon^2)$ for axis-aligned GMs. For high dimensional GM distributions,~\cite{kwon2020algorithm} show that with a suitable initialization and  $N\geq\tilde{O}((\min_{k\in[K]}{p}^{k \star})^{-1}(Q+R)/\epsilon^2)$, the Expectation Maximization (EM) algorithm converges in $T=O(\log(1/\epsilon))$ iterations to estimates $ \hat {p}^{k}$, $\hat{\mu}^{k}$, $({\hat\sigma}^{k})^2$ with accuracies $\epsilon_p=\max_{k\in[K]}p^{k\star}\epsilon$, $\epsilon_\mu=\max_{k\in[K]}\sigma^{k\star}\epsilon$, $\epsilon_\Sigma=(\max_{k\in[K]}\sigma^{k\star})^2\epsilon/\sqrt{Q+R}$, respectively.

\end{remark}

So far, we have introduced our GMM framework for contextual optimization and provided a theoretical analysis of the approximation quality.  However, real-world data may originate from more complex distributions that are not strictly GMMs. In the next section, we extend the applicability of our framework to general distributions.

\subsection{Handling General Probability Distributions}\label{subsec:normalizing_flow}
While GMMs possess expressive power that can model many real-world distributions fairly effectively, it may not accurately capture the precise characteristics of these distributions. Our solution scheme can, in principle, be extended to a mixture of parametric distributions since their conditional distributions are also closed-form, which can be tractable~\citep{cambanis1981theory}. However, in some situations, these assumptions remain restrictive, and the data may be more faithfully represented as (a mixture of) non-elliptical distributions. In this work, we strive to handle even more general distributions while still exploiting the unique property of GMM in computing the conditional distributions in closed form.
We propose employing normalizing flows \citep{dinh2014nice,rezende2015variational}, a powerful mechanism for representing complex distributions through a sequence of mappings from a simple base distribution, such as a Gaussian or a mixture of Gaussians. To our knowledge, the use of normalizing flows and Gaussian mixture models to compute conditional distributions of an intractable distribution has not previously been considered in the literature. 

Let $(\ts',\txi')\sim\PP$  be the data-generating random vector and $(\ts,\txi)\sim\M$ be the latent random vector governed by a tractable base distribution  $\M$, which in our case is set to a mixture of Gaussians. The framework seeks a differentiable and invertible function $T\in\RR^{Q+R}\rightarrow\RR^{Q+R}$ such that $T(\ts,\txi)$ follows the distribution $\PP$, which enables density estimation and sampling via the base distribution $\M$.  Under suitable regularity assumptions, normalizing flows are highly expressive: broad classes of target distributions can be represented or approximated by an invertible transformation of a simple base distribution. In practice, the function is typically restricted to a neural network representation $T_{\btheta}$ with parameters $\btheta$, enabling large-scale training with state-of-the-art algorithms. 

The normalizing flow framework has been successfully implemented in various large-scale tasks, including image generation \citep{ho2019flow++}, video generation \citep{kumar2019videoflow}, reinforcement learning \citep{mazoure2020leveraging,touati2020randomized}, and computer graphics \citep{muller2019neural}. We refer the reader to the survey papers \citep{kobyzev2020normalizing,papamakarios2021normalizing} for comprehensive reviews of normalizing flows. 
% \yc{We should add the $\bmu, \bSigma$ from the latent space can be trained by the algorithm(showing the structure of loss function)}

Let $\Q^\btheta$ be the pushforward distribution of $\M$ by $T_\btheta$.
Since $T_\btheta$ is a \emph{diffeomorphism}, i.e., both $T_{\bm \theta}$ and its inverse mapping $T_{\bm \theta}^{-1}$ are continuously differentiable, the pushforward density function can be evaluated in closed form via the change-of-variables formula \cite[Theorem 7.24]{rudin1987real}:
\begin{equation*}
f_{\Q^\btheta}(\s',\xii')
=
f_{\M}\left(T_{\btheta}^{-1}(\s',\bxi')\right)  
|\det \mathbf J_{T_{\btheta}^{-1}}(\s',\bxi') |,
\end{equation*}
where $\det \mathbf J_{T_{\btheta}^{-1}}(\s',\bxi')$ denotes the Jacobian determinant of $T_{\btheta}^{-1}$ at $(\s',\bxi')$. 
To find the best neural network parameters $\btheta$, we minimize the KL-divergence between the target density $f_{\PP}$ and the pushforward density induced by $T_\btheta$:
\begin{align*}
&\min_{\btheta} 
\DD\left(
f_{\PP}(\s',\xii')\; \big|\big|\;
f_{\M}\left(T_\btheta^{-1}(\s',\bxi')\right) 
|\det \mathbf J_{T_\btheta^{-1}}(\s',\bxi') |
\right).
\end{align*}
By the definition of the KL divergence, this is equivalent, up to an additive constant independent of $\btheta$, to the following optimization problem:
\begin{align*}
\min_{\btheta} 
-\EE_{\PP}\left[ 
\log f_{\M}\left(T_\btheta^{-1}(\ts',\txi')\right) 
+ \log |\det \mathbf J_{T_\btheta^{-1}}(\ts',\txi') | 
\right]. 
\end{align*}
Using the observed samples $\{(\s'_n,\bxi'_n)\}_{n\in[N]}$, the expectation in the objective is then approximated using the sample-average approximation
\begin{align}
\label{eq:KL_SAA_Transform}
\min_{\btheta} 
-\frac{1}{N}\sum_{n\in[N]}
\left(
\log f_{\M}\left(T_\btheta^{-1}(\s_n',\bxi_n')\right) 
+ \log |\det \mathbf J_{T_\btheta^{-1}}(\s_n',\bxi_n') |
\right), 
\end{align}
which is solved at scale using stochastic gradient descent. 

In this work, we further assume a separable structure in the transformation 
\[
T_\btheta(\s,\bxi)
=
\bigl(H_\btheta(\s),F_\btheta(\bxi;\s)\bigr),
\]
where \(H_\btheta\) is an invertible transformation of the contextual covariate and, for each fixed \(\s\),
\(F_\btheta(\cdot;\s)\) is invertible with respect to the uncertainty variable.
%which enforces the transformation for contextual covariates $\s'$ to depend solely on $\s$. 
This structure is relatively general and includes, as a special case, the popular \emph{autoregressive flows}~\citep{kingma2016improved,papamakarios2017masked,huang2018neural}. 
Note that, by construction, the inverse transformation is given by
\[
T_\btheta^{-1}(\s',\bxi')
=
\left(
H_\btheta^{-1}(\s'),
F_\btheta^{-1}\left(\bxi';H_\btheta^{-1}(\s')\right)
\right).
\]
% where $T_{\btheta,\s'}^{-1}(\s')$ is the inverse mapping with respect to $\s$, acting as a retransformation for the contextual information, while $T_{\btheta,\bxi'}^{-1}\left(T_{\btheta,\s'}^{-1}(\s'),\bxi'\right)$ denotes the inverse transformation with respect to $\bm \xi$. 

This normalizing flow architecture, including pushforward mapping and inverse mapping, allows for two key benefits. First, the transformations of $\s'$ or $\s$ are independent of the uncertain parameter, meaning we can transform the contextual information between the original and latent spaces without knowing the realization of $\txi'$. This aligns perfectly with the setup of contextual optimization, where decision-makers can only observe the realization of contextual information when making the decision. Second, the inverse transformation of $\bm \xi$ utilizes the projected contextual information $\s$, facilitating the computation of the conditional distribution of $\txi'$ given $\ts'=\s'$ in closed form, as described below.
\begin{proposition}\label{prop:normalizing_flow}
Let  $\s=H_\btheta^{-1}(\s')$. Then the conditional distribution \(\Q^\btheta_{\bxi'|\s'}\) has density
\[
f_{\Q^\btheta_{\bxi'|\s'}}(\bxi'|\s')
=
f_{\M_{\bxi|\s}}\left(F_\btheta^{-1}(\bxi';\s)\mid \s\right)
\left|
\det \mathbf J_{F_\btheta^{-1}}(\bxi';\s)
\right|,
\]
where \(\mathbf J_{F_\btheta^{-1}}(\bxi';\s)\) denotes the Jacobian of
\(F_\btheta^{-1}(\bxi';\s)\) with respect to \(\bxi'\), holding \(\s\) fixed.\end{proposition}
Proposition~\ref{prop:normalizing_flow} shows how the conditional density induced by the learned normalizing flow can be obtained from the latent conditional distribution. This result is useful for multistage stochastic problems, as illustrated in Section~\ref{sec:observational}. For the single-stage problem, however, the decision maker does not need to compute this density explicitly. Instead, one can solve the stochastic optimization problem by conditional sampling. Specifically, the normalizing flow is trained on the joint observations \((\s',\bxi')\), rather than by directly minimizing a conditional negative log-likelihood. Once the joint flow has been trained, the separable structure allows us to recover the latent contextual covariate \(\s=H_\btheta^{-1}(\s')\), sample from the latent conditional GM distribution \(\Ms\), and map the resulting samples back to the original space. The procedure is summarized in Algorithm~\ref{alg:conditional_sampling_flow}.
\begin{algorithm}[h!]
\caption{Conditional Sampling via a Separable Normalizing Flow}
\label{alg:conditional_sampling_flow}
\KwIn{Joint observations \(\{(\s_n',\bxi_n')\}_{n\in[N]}\), latent GM distribution \(\M\), observed contextual covariate \(\s'\), and sample size \(M\).}
\KwOut{Conditional samples \(\{\bxi_m'\}_{m\in[M]}\) drawn from \(\Q^\btheta_{\bxi'|\s'}\).}

Train the separable normalizing flow
\[
T_\btheta(\s,\bxi)=\left(H_\btheta(\s),F_\btheta(\bxi;\s)\right)
\]
by solving the joint maximum-likelihood problem \eqref{eq:KL_SAA_Transform}.

Compute the latent contextual covariate
\[
\s \leftarrow H_\btheta^{-1}(\s').
\]

Use Proposition~\ref{prop:conditional_param} to obtain the latent conditional GM distribution $\Ms$.

\For{\(m\in[M]\)}{
    Draw a latent conditional sample
    \[
    \bxi_m\sim\Ms.
    \]
    
    Map the latent sample back to the original space:
    \[
    \bxi_m'\leftarrow F_\btheta(\bxi_m;\s).
    \]
}

\Return \(\{\bxi_m'\}_{m\in[M]}\).
\end{algorithm}

We now provide a concrete example below.

\begin{example}
\begin{figure}[h]
\centering
\begin{subfigure}[t]{0.49\textwidth}
\centering
\includegraphics[width=1.0\textwidth]{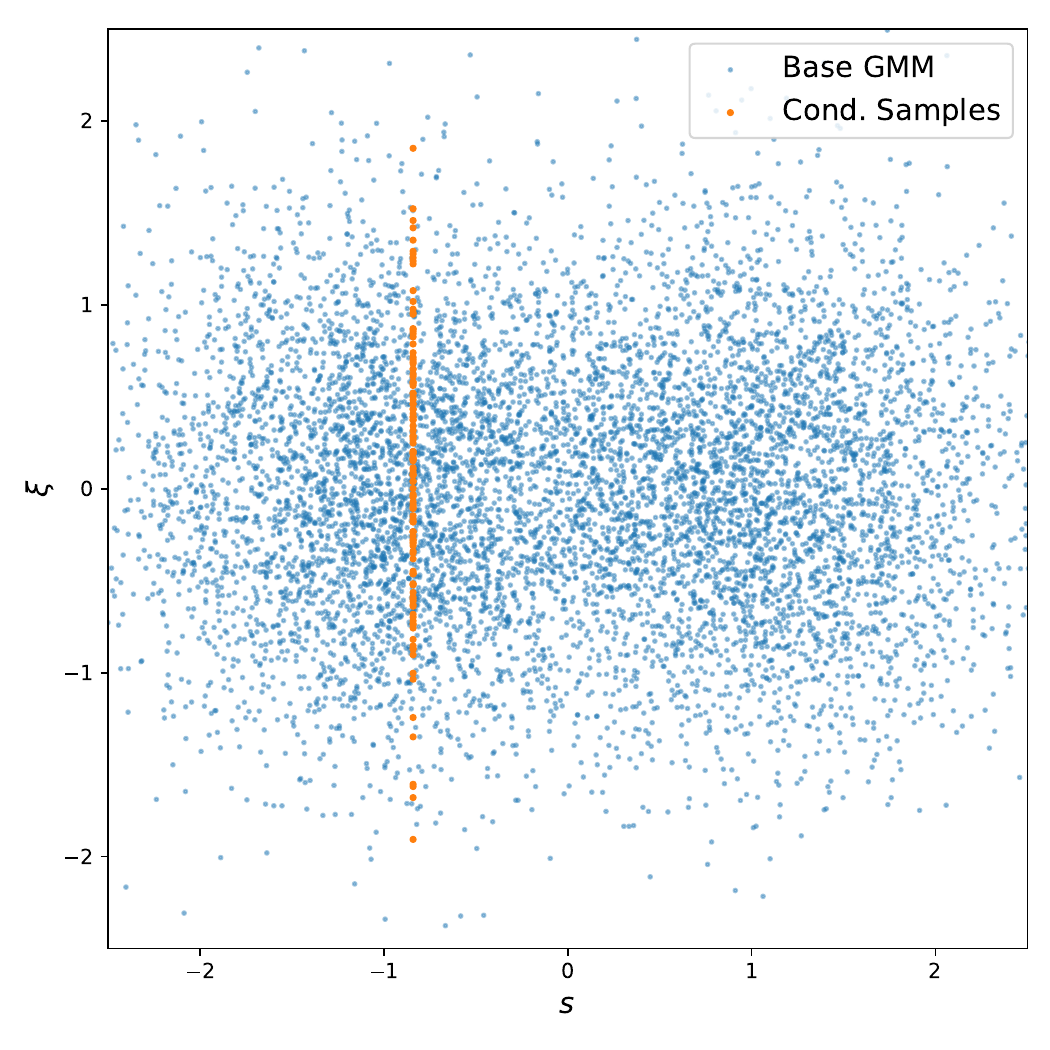}
\caption{Latent Space}
\end{subfigure}
\begin{subfigure}[t]{0.49\textwidth}
\centering
\includegraphics[width=1.0\textwidth]{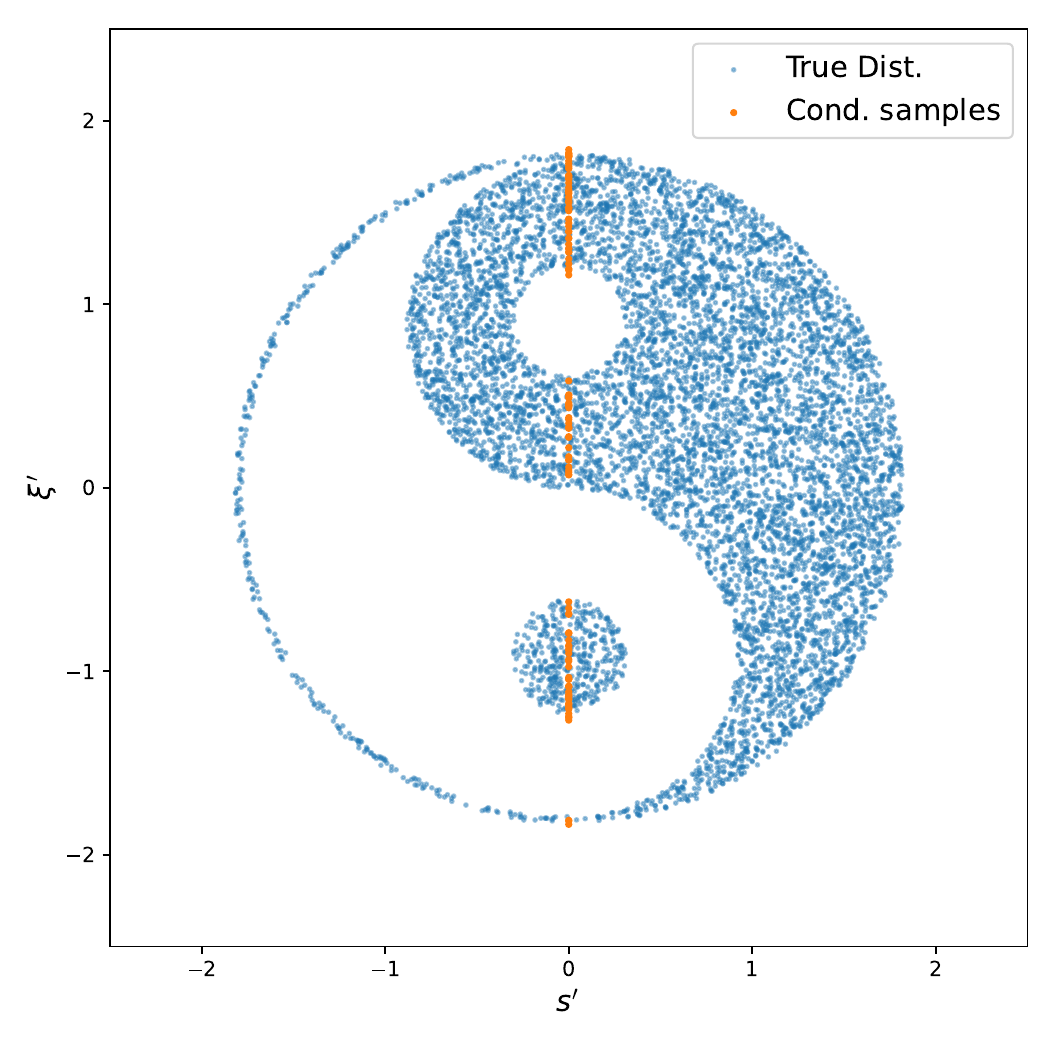}
\caption{Original Space}
\end{subfigure}
\caption{A visualization of the conditional sample generating scheme using normalizing flow. }
\label{fig:NF_example}
\end{figure}

This example illustrates the conditional sampling procedure enabled by our normalizing flow framework. As visualized by Figure~\ref{fig:NF_example}, the original samples, denoted by blue points in Figure~\ref{fig:NF_example}(b), come from a complex `Yin-Yang' shaped distribution\footnote{Yin and Yang is a Chinese philosophy that views the universe as governed by opposing yet complementary forces observable in nature, like the moon and sun, darkness and light, cold and heat.}. A normalizing flow is then trained to learn an invertible mapping between this true distribution and a base GM distribution in the latent space. Figure~\ref{fig:NF_example}(a) shows the mapped samples (blue dots) in the latent space. 

Suppose now we observe the side information $s' =0$. Then, using the trained normalizing flow, we can recover the latent covariate representation via $s=H_{\btheta}^{-1}(s')=-0.84$. Within the latent space, we leverage the GMM's analytical properties. Using the result from Proposition~\ref{prop:conditional_param}, we compute the conditional base distribution in closed form. Since this conditional distribution is also a GMM, an arbitrary number of samples can be efficiently generated, which is denoted by the orange points in Figure~\ref{fig:NF_example}(a). In the final step, these latent samples are propagated through the forward transformation, yielding the desired conditional samples in the original space.
\end{example}

\begin{remark}
A single Gaussian latent distribution is commonly used in the normalizing flow literature. In principle, one could rely entirely on the expressiveness of the normalizing flow to transform such a simple base distribution into the target distribution. We use a GM latent distribution because it provides a more structured and expressive base model. This reduces the burden on the flow when part of the target distribution can already be captured by a GMM, while preserving the closed-form conditional distribution \(\M_{\bxi|\s}\) in the latent space.
\end{remark}

Finally, we establish theoretical guarantees for our normalizing flow framework.  It is well-documented in the existing literature that general normalizing flows possess universal approximation properties. That is, given a sufficiently wide or deep enough transformation network $T_{\btheta}$, normalizing flows can approximate any target probability density to an arbitrary precision. Such universal approximation guarantees have been established across various flow architectures, such as Neural Autoregressive Flows~\citep{huang2018neural}, Block Neural Autoregressive Flows~\citep{de2020block}, and Coupling-based Flows~\citep{teshima2020coupling}. However, while these existing results ensure that the joint distribution can be approximated perfectly, our contextual optimization framework requires accurately sampling from the \emph{conditional} distribution. Therefore, it is important to verify that the universal approximation property is preserved under conditioning.

The following theorem formally bridges this gap, extending the universal approximation guarantees to the conditional distribution setting. It establishes that as the KL-divergence between the joint distributions tends to zero, the total variation distance between the conditional distributions also converges to zero, thereby offering the desired theoretical guarantee.

\begin{theorem}\label{thm:tv_conditional}
    Let $\PP$ and $\mathbb Q$ be two probability distributions on $\R^{Q+R}$ with strictly positive joint Lebesgue densities $f_{\PP},f_{\mathbb Q}$, respectively. Assume that both $f_{\PP}$ and $f_{\mathbb Q}$ are with H\"older continuous slice functions; that is, for all $\s_{1},\s_{2}\in\R^{Q}$, it holds that
\begin{equation*}
\int_{\RR^{R}} |f_{\PP}(\s_{1},\bxi) - f_{\PP}(\s_{2},\bxi)| {\rm d}\bxi \;\le\;K_{\PP}\,\|\s_{1}-\s_{2}\|^\alpha, \quad \int_{\RR^{R}} |f_{\mathbb Q}(\s_{1},\bxi) - f_{\mathbb Q}(\s_{2},\bxi)| {\rm d}\bxi \;\le\;K_{\mathbb Q}\,\|\s_{1}-\s_{2}\|^\alpha,
\end{equation*}
for some $\alpha\in (0, 1]$ and $K_{\PP} > 0$ and $K_{\mathbb Q}>0$. If the KL-divergence between $\PP$ and $\Q$ satisfies $\DD\left(f_{\PP} \big|\big|f_{\mathbb Q}\right) < \eps$, then 
\begin{equation*}
    \mathds{TV}(\PP_{\bxi | \s},\mathbb Q_{\bxi | \s}) \leq \frac{1}{f_{\PP}(\s)}\left(\left(\displaystyle\frac{\alpha}{Q}\right)^{\frac{Q}{Q+\alpha}} + \left(\displaystyle\frac{Q}{\alpha}\right)^{\frac{\alpha}{Q+\alpha}}\right) \left( K_{\PP} + K_{\mathbb Q}\right)^{\frac{Q}{Q+\alpha}}\left(\displaystyle\frac{\sqrt{2\eps}}{V_Q}\right)^{\frac{\alpha}{Q+\alpha}} \quad \forall \s\in \R^Q,
\end{equation*}
where $V_Q = \pi^{Q/2}/\Gamma(Q/2+1)$ is the volume of the unit Euclidean ball in $\R^Q$.
\end{theorem}

\section{Distributionally Robust Optimization Framework}\label{sec:DRO}

The GMM framework presented in Section~\ref{sec:contextual_GMM}, while powerful, relies on empirical estimation from limited historical observations. This can lead to solutions that overfit the training sample, yielding poor performance in out-of-sample scenarios. Furthermore, when learning the GM distribution, the model's structure (the number of mixture components K) may be misspecified. To systematically address these issues and enhance decision reliability, we integrate our approach with Distributionally Robust Optimization (DRO). In what follows, we will formulate this robust model, establish its reformulation, and derive performance guarantees that also account for potential misspecification of the model complexity.

\subsection{Reformulation and Performance Guarantee}
The DRO model seeks a solution that hedges against the worst-case distribution within a prescribed ambiguity set, i.e.,
\begin{equation}
\label{eq:DRO}
\min_{\x\in\X} \sup_{\Q\in\mP_\varepsilon}\EE_{\Q}\left[\ell(\x,\txi)\big|\ts=\s\right], 
\end{equation}
where the ambiguity set $\mathcal P_\varepsilon$ is a neighborhood of radius $\varepsilon\in\RR_+$, centered at the estimated $\Mhat$, defined with respect to a suitable probability metric. However, such a formulation suffers from degeneracy due to the conditioning on the singleton $\{\s\}$, which has zero measure. 
A common workaround is to condition on a set $\mathcal S$ of positive measure that contains $\s$ \citep{van2023generalized}. Such a scheme, however, may lead to difficulty in defining the structure and the size of $\mathcal S$ in practice. This limitation motivates the development of an alternative formulation that provides both finite-sample performance guarantees and asymptotic consistency.

A key advantage of our GMM framework is that the empirical conditional distribution $\Mhats$ is analytically available given the empirical estimation of $\M$. Leveraging this, we can directly construct a data-driven ambiguity set centered at this empirical conditional distribution, i.e., 
\begin{equation}
\label{eq:Wasserstein}
\mathcal P_\varepsilon\coloneqq \left\{\Q\in\mathcal P_0:\Wass_2\left(\Q,\Mhats\right)\leq\varepsilon\right\}. 
\end{equation}
Here, we adopt the type-2 Wasserstein distance $\Wass_2$ \citep{kantorovich1958space,gao2023distributionally} defined as
\begin{equation*}
	\label{eq:wasserstein}
\Wass_{2}(\Q_1, \Q_2) \coloneqq\inf_{\pi \in \Pi(\Q_1, \Q_2)} \left( \int_{\mathcal Z \times \mathcal Z} \|\bm z_1- \bm z_2\|^2 \, \pi ({\rm d}{\bm z_1}, {\rm d}{\bm z_2}) \right)^{\frac{1}{2}},
\end{equation*}
where $\Pi(\Q_1 \times \Q_2)$ is the set of all joint probability distributions of random vectors $\tilde\z_1$ and $\tilde\z_2$ with marginals $\Q_1$ and $\Q_2$, respectively.  The use of the type-2 Wasserstein distance facilitates theoretical performance guarantees, as we will discuss in the following section. Moreover, compared to the type-1 Wasserstein distance, it avoids pathological worst-case distributions and often yields more robust and higher-quality decisions \citep{byeon2025comparative}. 

Using this ambiguity set, we obtain a distributionally robust optimization problem that minimizes the worst-case \emph{unconditional} expectation with respect to all distributions within a Wasserstein ball centered at the estimated conditional distribution: 
\begin{equation}
\label{eq:cond_DRO}
\min_{\x\in\X} \sup_{\Q\in\mP_\varepsilon}\EE_{\Q}\left[\ell(\x,\txi)\right].
\end{equation}
The practical utility and theoretical performance guarantee of model~\eqref{eq:cond_DRO} require coverage of the true distribution.  To ensure that $\Ms\in\mathcal P_\varepsilon$, the radius $\varepsilon$ must be chosen appropriately. The following Theorem demonstrates that the radius required to ensure coverage scales polynomially with the estimation errors $\epsilon_p$, $\epsilon_\bmu$, and $\epsilon_\bSigma$.
% as established in the following theorem.
\begin{theorem}
\label{thm:coverage_radius}
Assume the conditions in Theorem \ref{thm:approx_error} hold. Then it is sufficient to set
\begin{equation}
\label{eq:coverage_radius}
\begin{array}{rl}
\varepsilon^2=&\left(\left(\frac{\beta}{\alpha}+1\right)\epsilon_\bmu+\frac{\alpha+\beta}{\alpha^2}\left(\|\bm s\|+\gamma \right)\epsilon_\bSigma\right)^2 + R\left(\frac{\beta}{\alpha} \right)^2 \epsilon_\bSigma+ 2\left(4\gamma^2 +2R\beta\right)\left(C_Q\|\s\|^2 + C'_Q\|\s\|+C''_Q \right). 
\end{array}
\end{equation}
to guarantee that $\Ms\in\mathcal P_\varepsilon$.
\end{theorem} 
This result yields the following out-of-sample performance guarantee for the solution to problem \eqref{eq:cond_DRO}. 
\begin{corollary}\label{coro:performance_bound}
Let $\hat\x$ and $\hat J$ be an optimal solution and the optimal value of the distributionally robust optimization problem \eqref{eq:cond_DRO}. If the radius $\varepsilon$ of the ambiguity set $\mP_\varepsilon$ is chosen according to \eqref{eq:coverage_radius}, then the following guarantee holds:
\begin{equation*}
\EE_{\Ms}[\ell(\hat\x,\txi)] \leq \hat J. 
\end{equation*}
\end{corollary}

Having established the performance guarantees of our DRO framework, we now turn to its practical implementation. By \cite[Remark 1]{blanchet2019quantifying}, if the loss function $\ell(\x,\bxi)$ is upper semicontinuous in $\bxi$, this problem is equivalent to the following minimization problem involving the estimated conditional distribution $\Mhats$:
\begin{equation*}
\label{eq:cond_DRO_min}
\min_{\x\in\X,\lambda\in\RR_+} \varepsilon^2\lambda+\EE_{\Mhats}\left[\sup_{\bm\omega\in\RR^R}\ell(\x,\bm\omega)-\lambda \|\bm\omega-\txi\|^2\right].
\end{equation*}
This formulation is amenable to the sample average approximation:
\begin{equation}
\label{eq:cond_DRO_min_SAA}
\min_{\x\in\X,\lambda\in\RR_+} \varepsilon^2\lambda+\frac{1}{M}\sum_{m\in[M]}\sup_{\bm\omega\in\RR^R}\ell(\x,\bm\omega)-\lambda \|\bm\omega-\xii_m\|^2,
\end{equation}
where $\{\bxi_m\}_{m\in[M]}$ are samples drawn from the conditional distribution $\Mhat_{\bxi|\s}$. When the loss function is piecewise affine, the resulting DRO problem admits a tractable second-order conic programming (SOCP) reformulation that can be efficiently solved using standard optimization solvers. Although the derivation is straightforward, to our knowledge, it does not appear explicitly in the existing literature. We therefore provide the following result.
\begin{proposition}\label{prop:dro_reformulation}If the loss function takes the form $\ell(\x,\xii)=\max_{j\in[J]}\bm a_j(\x)^\top\xii+b_j(\x)$, then the distributionally robust optimization problem \eqref{eq:cond_DRO_min_SAA} is equivalent to the second-order conic program
\begin{equation*}
\begin{array}{rl}
\min_{}\;\;\; &\displaystyle\varepsilon^2\lambda+\frac{1}{M}\sum_{m\in[M]}\gamma_m\\
\st\;\;& \displaystyle\x\in\X,\;\lambda\in\RR_+,\;\bm \gamma\in\RR^M\\
&\displaystyle \left\|\begin{bmatrix} 2\lambda\bxi_m + \bm a_j(\x)  \\
\lambda\xii_m^\top\xii_m + \gamma_m-b_j(\x)-\lambda \end{bmatrix}\right\| \leq \lambda\xii_m^\top\xii_m + \gamma_m-b_j(\x) + \lambda\quad\forall j\in[J] \quad\forall m\in[M]\\
&\displaystyle \lambda\xii_m^\top\xii_m + \gamma_m\geq b_j(\x)\quad\forall j\in[J] \quad\forall m\in[M].
\end{array}
\end{equation*}
\end{proposition}

\begin{remark}
The DRO model and its tractable reformulation presented in Proposition~\ref{prop:dro_reformulation} can be extended to cases where the random vector \((\ts',\txi')\) follows a general, non-GM distribution \(\PP\). Specifically, the sample average approximation in~\eqref{eq:cond_DRO_min_SAA} can be constructed using conditional samples \(\{\bxi'_m\}_{m\in[M]}\) generated by the normalizing-flow procedure described in Section~\ref{subsec:normalizing_flow}. Given the observed context \(\s'\) and the trained flow \(T_\btheta\), this procedure first computes the latent context
$\s=H_{\btheta}^{-1}(\s')$,
then draws latent samples 
$\bxi_m\sim \M_{\bxi|\s}$,
and finally maps them back to the original space through
$\bxi'_m=F_{\btheta}(\bxi_m;\s)$, $m\in[M]$. These generated samples are then used in~\eqref{eq:cond_DRO_min_SAA}.
\end{remark}

Together, Proposition~\ref{prop:dro_reformulation} and Remark~\ref{remark:GMM} demonstrate that our DRO framework is both computationally tractable and broadly applicable to general distributions. In the next part of this section, we will address the challenge of model misspecification.

 \subsection{Robustifying Against Unknown Mixture Size $K$} 

A practical challenge in applying the GMM framework is that the true number of mixture components $K$ is typically unknown and must be chosen as a hyperparameter. Selecting an incorrect number $K'$ can lead to model misspecification, which invalidates the performance guarantees derived in the previous section. Our DRO framework can be naturally extended to hedge against this structural uncertainty. Specifically, although the exact number of $K$ is unknown, the decision maker could generally determine a finite set $\mathcal K\subseteq \Z_+$ that constitutes all plausible candidate numbers.  
In this case, one could train GMM algorithms under different mixture sizes, resulting in $|\K|$ candidate conditional distributions: $\Mhats^L$ for $L\in\K$. Then, the decision maker can solve the distributionally robust optimization problem \eqref{eq:cond_DRO} using an ambiguity set of an appropriate radius centered at one of these distributions. The following proposition provides a recommendation for designing such an ambiguity set.

\begin{proposition}\label{prop:dro_k}
Define $\overline{\Wass}^2_2(\Mhats^K,\Mhats^{K'})$ for $K,K'\in\K$  to be the optimal value of the following linear program:  
\begin{equation*}
    \begin{array}{rl}
\overline{\Wass}^2_2(\Mhats^K,\Mhats^{K'}) \coloneqq\min & \displaystyle \sum_{i\in[K]}\sum_{j\in[K']}\pi_{ij} \left(\|\hat\bmu^i_{\bxi|\s}-\hat\bmu^j_{\bxi|\s}\|^2+\tr\left(\hat\bSigma_{\bxi|\s}^i+\hat\bSigma_{\bxi|\s}^j-2(((\hat\bSigma_{\bxi|\s}^i)^{\frac{1}{2}}\hat\bSigma_{\bxi|\s}^j(\hat\bSigma_{\bxi|\s}^i)^{\frac{1}{2}})^{\frac{1}{2}}\right)\right)\\
 \displaystyle \st& \displaystyle\bm\pi\in\Delta^K\times\Delta^{K'}\\
  & \displaystyle\sum_{i\in[K]}\pi_{ij}=\hat p^j\quad \forall j\in[K']\\
  &  \displaystyle \sum_{j\in[K']}\pi_{ij}= \hat p^i\quad \forall i\in[K].
  \end{array}
\end{equation*}
Assume the conditions in Theorem \ref{thm:approx_error} hold and the true GMM mixture size $K$ belongs to $\mathcal K$. Let $\varepsilon$ be defined as in Theorem \ref{thm:coverage_radius}. Then, centering the ambiguity set \eqref{eq:Wasserstein} at the conditional distribution $\Mhats^{K'}$, where
\begin{equation*}
K'\in\argmin_{K'\in\K} \max_{L\in\K} \overline{\Wass}^2_2(\Mhats^L,\Mhats^{K'}),
\end{equation*}
and setting the radius to 
\begin{equation*}
\varepsilon'\coloneqq \varepsilon+ \max_{L\in\K}  \overline{\Wass}_2(\Mhats^L,\Mhats^{K'}),
\end{equation*}
ensure that $\Ms^K\in\mathcal P_{\varepsilon'}$.
\end{proposition}

Proposition~\ref{prop:dro_k} offers a theoretically sound method for making decisions that are robust not only to data-driven estimation errors but also to the misspecification of the GMM's complexity. To the best of our knowledge, this is the first work to establish a formal coverage guarantee for GMM under misspecification. 

The models and reformulations developed so far leverage the analytical tractability of GMMs, enabling the approximation of the conditional distribution through efficient Monte Carlo sampling. This approach is highly effective for single-stage optimization problems. However, when extended to multistage settings, standard sampling-based dynamic programming methods encounter a significant challenge: the size of the resulting problem often grows exponentially with the number of planning horizons $T$. To overcome this limitation, the next section introduces an alternative approximation scheme that avoids direct sampling, offering a solution method whose complexity scales much more favorably with the planning horizon.

\section{Optimization Models Using Observational Data}\label{sec:observational}
% In this section, we develop optimization models that directly utilize the observational data, bypassing the need for Monte Carlo sampling to approximate the conditional expectations. These models offer significant advantages in multistage settings, particularly when solving stochastic optimization problems using data-driven dynamic programming \citep{park2022data}, where the problem size scales only linearly in the time horizon and data size. 
While our conditional sampling framework is highly effective for single-stage problems, extending it to multistage settings introduces severe computational barriers. In this section, we develop optimization models that directly utilize the observational data, bypassing the need for Monte Carlo sampling to approximate the conditional expectations. Operating within an observational paradigm allows our proposed models to scale only linearly with the planning horizon and data size, ensuring tractability in data-driven dynamic programming.

\begin{figure}
    \centering
    \includegraphics[width=0.95\linewidth]{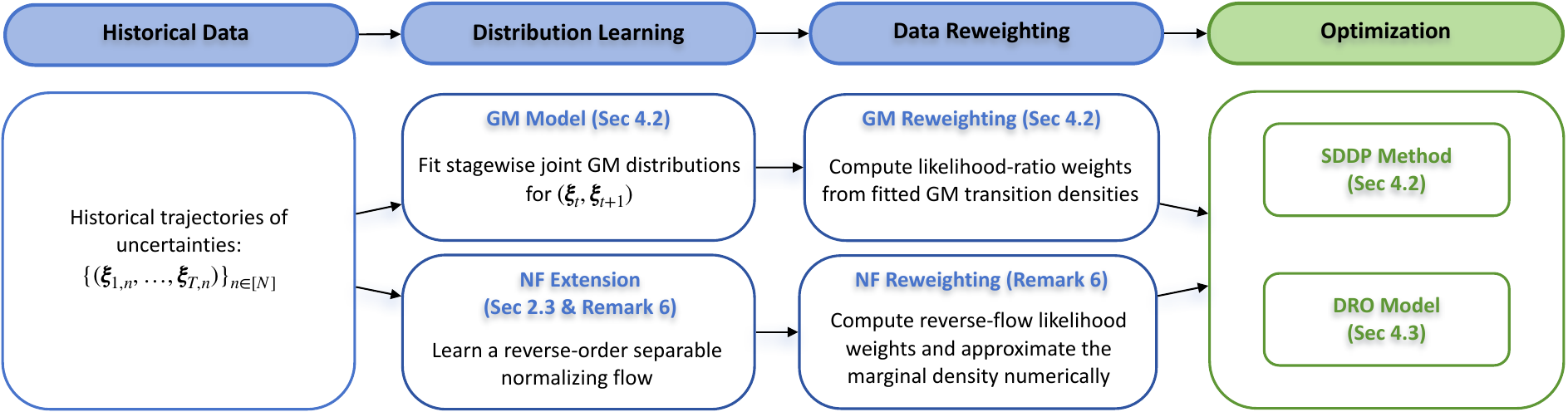}
    \caption{An overview of the GMM contextual optimization framework for multistage optimization problems.}
    \label{fig:structure_multistage}
\end{figure}

The traditional scheme for multistage stochastic programming is the sample-average approximation (SAA) scheme \citep{shapiro2011analysis}. SAA presumes full access to the underlying distribution and relies on sequential conditional Monte Carlo sampling to capture the temporal evolution of the stochastic process. If one were to apply the SAA scheme by sequentially generating new samples from an estimated distribution (such as $\Mhat$) at each stage, the resulting scenario tree would grow exponentially with the time horizon, rendering the dynamic programming computationally intractable. To circumvent this exponential explosion, \cite{park2022data} propose a data-driven scheme that operates directly on historical trajectories by reweighting the transition probabilities between consecutive data points using kernel regression. However, local non-parametric kernel regression is notoriously vulnerable to the ``curse of dimensionality" in high-dimensional state spaces. To address both bottlenecks simultaneously, we propose to replace these kernel-based weights with Gaussian-mixture-based weights, leveraging a GMM to approximate conditional transitions. This substitution preserves the direct-from-data workflow while mitigating the curse of dimensionality inherent in kernel regression, and enables scalable multistage optimization without resorting to large-scale Monte Carlo sampling.

We begin by describing the solution procedure in a single-stage setting. Let $f$ and $\fhat$ be the density functions of $\M$ and $\Mhat$, respectively. Given $N$ historical observations $\{(\s_n,\bxi_n)\}_{n\in[N]}$ drawn from $\M$, we approximate the contextual stochastic optimization problem \eqref{eq:cso} by 
\begin{equation}
\label{eq:observational_cse}
\min_{\x\in\X} \frac{1}{\fhat(\s)N}{\sum_{n\in[N]}\ell(\x,\bxi_n)\hat f(\s|\bxi_n) }.
\end{equation}
We now justify this approximation.
The conditional expectation in the objective can be rewritten as 
\begin{align}
\EE_{\M}\left[\ell(\x,\txi)\big|\ts=\s\right]&={\int \ell(\x,\xii) f(\s,\xii) /f(\s){\rm d}\xii}\\
&=\int \ell(\x,\xii) f(\s|\xii) f(\xii)/f(\s){\rm d}\xii\\
&=\frac{1}{f(\s)}\EE_{\M}\left[\ell(\x,\txi)f(\s|\txi)\right].
\end{align}\label{eq:observable}
Using the samples $\{\bxi_n\}_{n\in[N]}$ and applying sample-average approximation yields:
\begin{align*}
\EE_{\M}\left[\ell(\x,\txi)\big|\ts=\s\right]\approx\frac{1}{N}{\sum_{n\in[N]}\ell(\x,\bxi_n)\frac{f(\s|\bxi_n)}{f(\s)} }.
\end{align*}
This approximation can be interpreted as an importance-sampling estimator for conditional expectations. Specifically, samples from the joint distribution \(\M\) are reweighted by the likelihood ratio $\frac{f(\s|\bxi_n)}{f(\s)}$
to estimate the conditional expectation under \(\M_{\bxi|\s}\).
The factor $\frac{f(\s|\bxi_n)}{f(\s)}=\frac{f(\s,\bxi_n)}{f(\s)f(\bxi_n)}$ quantifies the local dependence between $\s$ and $\bxi_n$. Values greater than one indicate that the co-occurrence $(\s,\bxi_n)$ is more likely than under independence, while values smaller than one indicate the opposite \citep{church1990word}. Lastly, the formulation \eqref{eq:observational_cse} is obtained by replacing the unknown true density $f$ with the approximation $\hat f$. 

\begin{remark}
The identity
\[
\EE\left[\ell(\x,\txi)\mid\ts=\s\right]
=
\frac{1}{f(\s)}
\EE\left[\ell(\x,\txi) f(\s|\txi)\right]
\]
holds for any joint distribution for which the relevant densities are well defined. However, for a general joint distribution, computing the marginal density
\[
f(\s)=\int f(\s,\bxi)\,\rm d\bxi
\]
may require numerical integration that is computationally prohibitive. The GM model avoids this difficulty: when the joint distribution is modeled as a GM, both the marginal density \(f(\s)\) and the conditional density \(f(\s|\bxi)\) are available in closed form. This enables the efficient likelihood-ratio reweighting scheme used in the preceding approximation.
\end{remark}
\begin{remark}
The approximation scheme can be generalized to a random vector \((\ts',\txi')\) governed by a general distribution through a normalizing flow. Specifically, consider a reverse-order separable flow
\[
T_\btheta(\s,\bxi)=\left(G_\btheta(\bxi),R_\btheta(\s;\bxi)\right).
\]
By the same change-of-variables argument as in Proposition~\ref{prop:normalizing_flow}, the conditional density
\(f_{\Q^\btheta_{\s'|\bxi'}}(\s'|\bxi')\) can be evaluated under the learned pushforward distribution \(\Q^\btheta\). Indeed, for
\[
\bxi=G_\btheta^{-1}(\bxi'),\qquad
\s=R_\btheta^{-1}(\s';\bxi),
\]
we have
\[
f_{\Q^\btheta_{\s'|\bxi'}}(\s'|\bxi')
=
f_{\M_{\s|\bxi}}(\s|\bxi)
\left|\det \mathbf J_{R_\btheta^{-1}}(\s';\bxi)\right|,
\]
where the Jacobian is taken with respect to \(\s'\), holding \(\bxi\) fixed. The corresponding likelihood-ratio weight is
\[
\frac{f_{\Q^\btheta_{\s'|\bxi'}}(\s'|\bxi')}
{f_{\Q^\btheta_{\s'}}(\s')}.
\]
Unlike in the GM case, the marginal density \(f_{\Q^\btheta_{\s'}}(\s')\) is generally not available in closed form under this reverse-order flow. It can, however, be approximated by Monte Carlo integration:
\[
f_{\Q^\btheta_{\s'}}(\s')
=
\EE_{\bxi'\sim \Q^\btheta_{\bxi'}}
\left[
f_{\Q^\btheta_{\s'|\bxi'}}(\s'|\bxi')
\right]
\approx
\frac{1}{L}\sum_{m=1}^L
f_{\Q^\btheta_{\s'|\bxi'}}(\s'|\bxi'_m),
\]
where \(\bxi'_m=G_\btheta(\bxi_m)\) and \(\bxi_m\) is sampled independently from the marginal latent distribution \(\M_{\bxi}\). Since samples can be generated independently from the trained flow model, the approximation error due to Monte Carlo integration can be reduced by increasing \(L\).
\end{remark}

\subsection{Approximation Quality}
We now proceed to establish the approximation quality of our proposed scheme, demonstrating its favorable sample complexity properties. To this end, we introduce the following regularity condition:
\begin{enumerate}[(D)]
\item \label{as:density_2} There exist constants $\fbar,\fl\in\RR_{++}$ such that $\fbar\geq \fhat(\s|\bxi)$ for all $\bxi\in\RR^R$, and $\fl\leq \fhat(\s)$ with probability $1-\delta$. 
\end{enumerate} 
This condition is standard in the nonparametric regression literature \citep{gyorfi2006distribution,kohler2009optimal,belkin2019does} and serves to simplify the presentation. Specifically, since  $\hat f(\s|\bxi)$ is a conditional density of a mixture of Gaussians---and hence itself a mixture of Gaussians---it admits a uniform upper bound. Likewise, since $\hat f(\s)$ is a marginal density of a Gaussian mixture, one can always determine a lower bound $\underline{f}$ given a confidence level $\delta$. Based on these assumptions, the following theorem bounds the difference between the true and empirical conditional expectations.

\begin{theorem}\label{thm:error_observe_fix_x}
Let $\M=\sum_{k\in[K]}p^{k}\N(\bmu^{k},\bSigma^{k})$ be the true Gaussian mixture distribution and $\Mhat=\sum_{k\in[K]}\hat p^k\Nhat(\hat\bmu^k,\hat\bSigma^k)$ be the estimated one. Suppose that $|p^{k}-\hat p^k|\leq \epsilon_p$, $\|\bmu^{k}-\hat\bmu^k\|\leq\epsilon_\bmu$, and $\|\bSigma^{k}-\hat\bSigma^k\|\leq\epsilon_\bSigma$ for all $k\in[K]$. Then, for any $\delta\in(0,1)$, we have
\begin{equation}
\label{eq:fix_x_bound}
\begin{array}{rl}
&\displaystyle\;\;\;\;\left|\EE_{\M}\left[\ell(\x,\txi)\big|\ts=\s\right]-\frac{1}{\fhat(\s)N}\sum_{n\in[N]}\ell(\x,\bxi_n)\hat f(\s|\bxi_n)\right|\\
\leq & \lbar\left[\;C_{Q+R}\left( \|\s\|^2  + \beta + \left(\gamma+\frac{\beta}{\alpha}(\|\s\|+\gamma)\right)^2\right)+C'_{Q+R}\sqrt{ \|\s\|^2  + \beta + \left(\gamma+\frac{\beta}{\alpha}(\|\s\|+\gamma)\right)^2} + C''_{Q+R}\right] \\
&\;\; + \frac{\lbar\fbar}{\fl} \left(C_R \left(\beta+\gamma^2\right)+C'_R\sqrt{\beta+\gamma^2}+C''_R+C_{Q}\| \s\|^2+C'_{Q}\| \s\| + C''_{Q}\right)\\
&\quad+\frac{\lbar\fbar}{\fl}\sqrt{\frac{1}{2N}\log\left(\frac{8}{\delta}\right)}.
\end{array}
\end{equation}
with probability at least $1-\delta$. Here, the constants are defined as   $C_D=\left(1 + \frac{\epsilon_p}{\underline{p}} \right)  \frac{D \epsilon_\bSigma}{2 \alpha^2}$, $C'_D=\left(1 + \frac{\epsilon_p}{\underline{p}} \right)  \left( \frac{\epsilon_\bmu}{\alpha} + \frac{D \epsilon_\bSigma  \gamma}{\alpha^2} \right)$, and $C''_D=\frac{\epsilon_p}{\underline{p}} 
+ \left(1 + \frac{\epsilon_p}{\underline{p}} \right)
\left( \frac{\epsilon_\bmu}{\alpha} \, \gamma 
+ \frac{D \epsilon_\bSigma}{2 \alpha^2} (\gamma^2 + \alpha) 
+ o(\epsilon_\mu + \epsilon_\bSigma) \right)$. 
\end{theorem}
In Theorem~\ref{thm:error_observe_fix_x}, the constants $C_{Q+R}$, $C'_{Q+R}$,  $C''_{Q+R}$,  $C_{R}$, $C'_{R}$,  $C''_{R}$ depend linearly on $\epsilon_p$, $\epsilon_{\bmu}$, $\epsilon_{\bSigma}$, and the final term decays as $O(\frac{1}{\sqrt N})$; therefore, the proposed approximation enjoys polynomial sample complexity and avoids the curse of dimensionality. In the next section, we apply this approximation framework to address multistage stochastic optimization problems under Markovian uncertainty.

\subsection{Applications to Multistage Stochastic Programming}
\label{sec:multistage}
We consider a multistage stochastic optimization problem with Markovian uncertainty over a planning horizon of $T$ stages, given by:
\begin{equation}\label{true_multi_form}
    \min_{\substack{\bm x_1 \in \X_1(\x_0,\bm \xi_1) }} \ell(\x_1,\xii_1) + \EE \left[ 
        \min _{ \substack{\bm x_2 \in \X_2(\x_1, \tilde{\bm \xi}_2) }}\ell(\x_2,\txi_2)+ \EE \left[
                    \cdots + \EE \left[ \min_{ \substack{ \bm x_T \in \X_T(\x_{T-1},\tilde{\bm \xi}_T)} } \ell(\x_T,\txi_T) \;\middle\vert\; \txi_{T-1} \right] \cdots  \;\middle\vert\; \txi_2 
                    \right]  \;\middle\vert\; \bm \xi_1 
    \right].
\end{equation}
We follow the standard multistage stochastic programming convention in which the stagewise decision vector \(\x_t\) may include both physical state variables and control/action variables; see, e.g., \cite[Section~3]{shapiro2021lectures}. Under this convention, the state-transition equations are incorporated into the stagewise feasible constraints. Equivalently, one could adopt a stochastic-control notation that separates the state and control variables explicitly; see \cite[Remark~2.4]{fullner2025stochastic}. We use the compact stochastic-programming notation to simplify the exposition.

This problem is traditionally solved via dynamic programming by recursively evaluating the cost-to-go functions from the terminal stage $t=T$ backward to the initial stage $t=1$:
\begin{equation}
\label{eq:non-dro_true_formulation_T}
V_t(\bm x_{t-1}, \bm \xi_t) = \min_{\bm x_t \in \X_t(\x_{t-1},\tilde{\bm \xi}_t)}  \bm c_t^\top \bm x_t + \mathcal V_{t+1} \left( \bm x_t, \bm \xi_t \right)  
\quad\forall \x_{t-1}\in\X_{t-1}\;\;\forall \bxi_t\in\RR^Q,
\end{equation}
where 
\begin{equation}
\label{eq:conditional_expected_cost_to_go}
\mathcal V_{t+1} ( \bm x_t, \bm \xi_t )
= \EE \left[ V_{t+1}(\bm x_{t}, \tilde{\bm \xi}_{t+1} ) \;\middle \vert\; \bm \xi 
_t \right]
\end{equation} 
denotes the conditional expectation of the future cost-to-go function at stage $t+1$, given the most recent realization $\bm \xi_t$. We assume $V_{T+1}(\cdot) = 0$, which implies that no additional costs beyond the terminal stage $T$.

In practice, the underlying stochastic process  $\txi_{[T]} = (\txi_1, \txi_2, \dots, \txi_T)\in\RR^{Q\times T}$ is typically unknown and one only has access to $N$ i.i.d.~sample trajectories $\{\bm \xi_{[T],n} \coloneqq (\bm \xi_{1,n}, \dots, \bm \xi_{T,n})\}_{n\in[N]}$. In \citep{park2022data}, the authors propose to approximate the conditional expectation in a data-driven fashion
\begin{equation} \label{eq:apprx_exptected_ctg}
    \hat{\mathcal V}_{t+1} ( \bm x_{t}, \bm \xi_{t}) 
    \approx
    \sum_{i \in[N]} {w}_{t+1}(\bm \xi_{t},\bm \xi_{t,i}) {V}_{t+1}( \bm x_{t}, \bm \xi_{t+1,i} ),
\end{equation}
where the weights ${w}_{t+1}(\bm \xi_{t},\bm \xi_{t,i})$ are defined through the kernel regression. This scheme avoids costly evaluations over the continuous space by computing the value functions only at the sample points $\{\bxi_{t,n}\}_{n\in[N]}$.   However, this method suffers from the curse of dimensionality: its suboptimality scales as $\tilde{\mathcal O}(T^{\frac{3}{2}}/N^{\frac{2}{Q+4}})$ and thus deteriorates rapidly as the dimension $Q$ increases.

To overcome this limitation, we propose to replace the kernel regression weights in \eqref{eq:apprx_exptected_ctg} with the Gaussian-mixture-based observational weights in \eqref{eq:observational_cse}, yielding the approximation
\begin{equation} \label{eq:apprx_exptected_ctg_new}
    \hat{\mathcal V}_{t+1} ( \bm x_{t}, \bm \xi_{t}) 
    \approx \hat\EE\left[ {V}_{t+1} ( \bm x_{t}, \tilde{\bm \xi}_{t+1}) \vert \tilde{\bm \xi}_{t}= \bm \xi_{t}  \right]= {\frac{1}{N}}\sum_{i \in[N]} \frac{\fhat (\bxi_{t}|\bm \xi_{t+1,i})}{\fhat (\bxi_{t}) } {V}_{t+1}( \bm x_{t}, \bm \xi_{t+1,i} ).
    %\sum_{i \in[N]} \frac{\fhat (\bxi_{t}|\bm \xi_{t+1,i})}{\sum_{j\in[N]}\fhat (\bxi_{t}|\bm \xi_{t+1,j}) } {V}_{t+1}( \bm x_{t}, \bm \xi_{t+1,i} ). 
\end{equation}
This leads to the following suboptimality guarantee for the solution $\widehat{\bm x}_1^N$ obtained via our data-driven dynamic programming approach.
\begin{theorem}\label{thm:out_of_sample}
Assume that every consecutive pair \((\txi_t,\txi_{t+1})\), \(t\in[T-1]\), follows a Gaussian mixture distribution.
Assume further that, for each $t$, the cost-to-go functions $V_t(\x_{t-1},\bxi_t)$ are $L$-Lipschitz continuous in $\x_{t-1}$ for any fixed $\bxi_t$. Furthermore, assume that the feasible regions are nonempty and compact: there exists a constant $\overline D\in\RR_{++}$ such that   $\sup _{\x_{t},\x_{t}' \in \mathcal X_{t}(\x_{t-1}, \bm \xi_t)  } \Vert \x_{t} - \x_{t}' \Vert \leq \overline D$ for all $t\in[T]$, $\x_{t-1}\in\mathcal X_{t-1}(\cdot)$, and $\bm \xi_t\in\Xi_{t}$.
Then for any $\eta\in\RR_{++}$, the following bound holds with probability at least $1-\delta$.
\begin{align*}
\label{eq:out_of_sample} 
    &\bm \ell(\widehat{\bm x}_1^N,\xii_1)+ {\mathcal V}_{2}(\widehat{\bm x}_1^N, \bm \xi_1)
    - \left(\min_{\bm x_1 \in \X_1(\x_0,\bm \xi_1) }
    \ell(\x_1,\xii_1) + {\mathcal V}_{2}({\bm x}_1, \bm \xi_1)\right)\\
    &\quad\quad      \leq 2(T-1)\tau + 4(T-1)L\eta + 2\sum_{t=2}^{T}
        \frac{\lbar\fbar^2}{\fl^2}
         \sqrt{\frac{2}{N} \log \left(\frac{ (T+1)\mathcal O(1) N^{t-2} (\overline D / \eta)^{R(t-1)} }{ \delta }\right) }
\end{align*}
Here,  
\begin{equation}
\label{eq:tau_bound}
\begin{array}{rl}
\tau\coloneqq& \lbar\left[\;C_{Q+R}\left( \bar\xi^2  + \beta + \left(\gamma+\frac{\beta}{\alpha}(\bar\xi+\gamma)\right)^2\right)+C'_{Q+R}\sqrt{ \bar\xi^2  + \beta + \left(\gamma+\frac{\beta}{\alpha}(\bar\xi+\gamma)\right)^2} + C''_{Q+R}\right] \\
&\;\; + \frac{\lbar\fbar}{\fl} \left(C_R \left(\beta+\gamma^2\right)+C'_R\sqrt{\beta+\gamma^2}+C''_R+C_{Q}\bar\xi^2+C'_{Q}\bar\xi + C''_{Q}\right). 
\end{array}
\end{equation}
and $\bar\xi\coloneqq 4\sigma\sqrt{R}+2\sigma\sqrt{2\log\tfrac{NT}{\delta}} + \gamma$. 
\end{theorem}
This result establishes that the suboptimality of the proposed approximation decays at the rate 
\begin{align*}
 &\;\;\mathcal O\left(\sum_{t=2}^{T}
        \sqrt{\frac{2}{N} \log \left(\frac{\mathcal O(1) N^{t-2} (\overline D / \eta)^{R(t-1)} }{ \delta }\right) }\right)\\
         =&\;\;\mathcal O\left(\sum_{t=2}^{T}
         \sqrt{\frac{2}{N} \left(\log \mathcal O(1) + (t-2)\log N +  R(t-1) \log \frac{\overline D}{ \eta} -\log\delta\right)}\right)\\
        =&\;\;\tilde{\mathcal O}\left(\sum_{t=2}^{T}
         \sqrt{\frac{t}{N}}\right)=\tilde{\mathcal O}\left(\frac{T^{3/2}}{N^{1/2}}\right),
\end{align*}
and crucially, it is free from the curse of dimensionality. In contrast, the kernel-based approach yields a slower rate of $\tilde{\mathcal O}(T^{3/2}/N^{2/(Q+4)})$ \citep[Theorem 1]{park2022data}. To the best of our knowledge, this is the first result of its kind, even for the special case where the underlying stochastic process is a single multivariate Gaussian ($K=1$), a setting that frequently arises in multistage portfolio optimization. Notably, the standard sample average approximation approach in this setting only achieves a convergence rate of $\tilde{\mathcal O}(T^{1/2}/N^{1/(2T)})$, which deteriorates significantly as the planning horizon $T$ increases \citep{shapiro2005complexity}.

\subsection{Distributionally Robust Optimization Framework}
To enhance the reliability of our data-driven approach, we now embed our model within a DRO framework. To this end, we define the \emph{empirical measure} as
\[\hat\Q\coloneqq\sum_{n\in[N]}\frac{\hat f(\s|\bxi_n)}{\sum_{u\in[N]} \hat f(\s|\bxi_{u})}\delta_{\xii_n},\]
and use it as the center of the Wasserstein ambiguity set defined in \eqref{eq:Wasserstein}.
Reformulating the resulting DRO problem into a tractable conic program is straightforward by employing Proposition~\ref{prop:dro_reformulation}. For the sake of brevity, we omit this derivation and focus on the illustration of how to determine an appropriate radius for the Wasserstein ambiguity set. 

The practical application and performance guarantees of the DRO model hinge on choosing a radius that ensures the true conditional distribution $\Ms$ is contained within the ambiguity set with high probability.
To achieve this, we need a valid upper bound on $\Wass_2(\hat\Q,\Ms)$. Let $\Mhats^M\coloneqq \frac{1}{M}\sum_{m\in[M]} \delta_{\bxi'_m}$ be the empirical measure constructed from $M$ samples of $\Mhats$.  By the triangle inequality:
\begin{equation*}
\Wass_2(\hat\Q,\Ms)\leq \Wass_2(\hat\Q,\Mhats^M) + \Wass_2(\Mhats^M,\Mhats) +  \Wass_2(\Mhats,\Ms). 
\end{equation*}
The first term, $\Wass_2(\hat\Q, \Mhats^M)$, corresponds to the square root of the optimal value of the following linear program:
\begin{equation*}
    \begin{array}{rl}
 \displaystyle \Wass^2_2(\hat\Q,\Mhats^M)= \min& \displaystyle \sum_{i\in[N]}\sum_{j\in[M]}\pi_{ij} \|\bxi_i-\bxi'_j\|^2\\
 \displaystyle \st& \displaystyle\bm\pi\in\Delta^N\times\Delta^M\\
  & \displaystyle\sum_{i\in[N]}\pi_{ij}=\frac{1}{M}\quad \forall j\in[M]\\
  &  \displaystyle \sum_{j\in[M]}\pi_{ij}= \frac{\hat f(\s|\bxi_i)}{\sum_{u\in[N]} \hat f(\s|\bxi_{u})}\quad \forall i\in[N].
  \end{array}
\end{equation*}
The second term, $\Wass_2(\Mhats^M, \Mhats)$, can be bounded using concentration results for empirical Wasserstein distances \cite[Theorem 2]{fournier2015rate}. 
The final term, $\Wass_2(\Mhats, \Ms)$, is upper bounded by the result established in Theorem~\ref{thm:coverage_radius}.

\begin{remark}
The proposed DRO framework can be seamlessly embedded into the multistage stochastic programming formulation described in Section~\ref{sec:multistage}. Specifically, at each stage $t$, instead of approximating the conditional expectation $\mathcal V_{t+1} ( \bm x_t, \bm \xi_t )$ in \eqref{eq:conditional_expected_cost_to_go}  via a fixed weighted average as in \eqref{eq:apprx_exptected_ctg_new}, we solve a distributionally robust optimization problem over a Wasserstein ball centered at the empirical measure $\Q$. 
\end{remark}

\section{Experimental Results}\label{sec:experiment}
In this section, we conduct numerical experiments to evaluate the out-of-sample performance of different decision-making models. We study three important operations management problems: inventory management, portfolio optimization, and the multistage wind energy planning problem. All experiments are conducted on a server equipped with an AMD Ryzen Threadripper 9970X 32-Core processor, 251 GiB of RAM, and an NVIDIA GeForce RTX 5090 GPU.
% In each experiment, we begin by describing the problem setting, followed by the model training and hyperparameter tuning procedures, and conclude with out-of-sample evaluation results and their interpretation.

\subsection{Inventory Management}

% \textbf{Problem setting:}
We consider an inventory management problem, where the objective of the decision maker is to determine an optimal order quantity under uncertain demand. Specifically, when the order quantity exceeds actual demand, the firm incurs a holding cost. Conversely, when the order quantity is too low, a stock-out cost is penalized due to unmet demand. The loss function under order quantity $q\in\RR_+$ and realized demand $\xi\in\RR_+$ is given by:
\[
\ell(q, \xi) = h(q-\xi)_+ + b(\xi-q)_+,
\]
where $h$ and $b$ denote the per-unit holding and stock-out costs, respectively. In this experiment, we set these costs to $h=10$ and $b=2$. We assume the demand is influenced by a vector of observable side information $\ts \in \RR^Q$, which is available prior to making the ordering decision.

To assess the capabilities of our proposed method, we design a synthetic data-generating process that simulates a complex, state-dependent demand environment. We assume the demand is driven by a $Q$-dimensional vector of observable contextual information, $\ts = [\tilde s_1, \dots, \tilde s_Q]^\top \in \mathbb{R}^Q$. In our experiments, we draw each contextual feature independently from a uniform distribution, such that $\tilde{s}_i \sim \mathbb{U}(-2, 2)^Q$. We construct the conditional distribution of the demand $\txi$ given $\ts=\s$ as a multimodal distribution combining a linear and a quadratic regime:
\[
\xi = \begin{cases}
 0.3\mathbf e^\top \s  + \delta_1 + \mu_1 & \text{with probability } p(\s), \\
 5\sum s_i^2 + \delta_2 + \mu_2 & \text{with probability } 1-p(\s),
\end{cases}
\]
where $\mathbf{e} \in \mathbb{R}^Q$ is a vector of ones, and $\tilde{\delta}_1, \tilde{\delta}_2 \in \mathbb{R}$ are independent, one-dimensional uniform noise terms drawn from $\mathbb{U}(-2, 2)$. The constants $\mu_1 = 50$ and $\mu_2 = 40$ represent the base demands for the two regimes, respectively. To capture the context-driven regime shift, we define the state-dependent mixture probability as $p(\s) = \frac{1}{2}(1 + \tanh(s_1))$, implying that the probability of the demand following the linear regime increases with the first contextual feature $s_1$. This setup presents a highly challenging scenario for contextual optimization. It features an $Q$-dimensional context space, a non-Gaussian multimodal target distribution, and a complex non-linear dependency in which the covariates $\s$ dictate both the value of the demand and the underlying data-generating regime.

We evaluate the performance of our proposed models—both the plain-vanilla GMM framework (GMM) and its extension with normalizing flows (GMM-NF)—against four prominent benchmarks from the literature. The first benchmark is the regularized Nadaraya-Watson regression (RNW) method~\citep{wang2026data}. The RNW method is known for enhancing the robustness of kernel-based approaches by incorporating a conditional variance regularization term. It serves as a strong baseline for non-parametric models that do not assume a specific data structure but can be susceptible to the curse of dimensionality. The second benchmark is the residual-based distributionally robust contextual (ResDRO) optimization approach~\citep{kannan2024residuals}. ResDRO's key strength lies in its flexible residual architecture, which allows it to be effectively paired with various parametric prediction models. This method also employs a DRO framework to improve out-of-sample stability, and it is regarded as a state-of-the-art parametric framework for contextual optimization. We then include two prominent baselines from the decision rule paradigm. The third benchmark is the Linear Decision Rule (LDR) approach, which restricts the decision policy to an affine function of the covariates~\citep{ban2019big}. The fourth benchmark is the Neural Network Decision Rule (NN-DR) approach, which leverages the representational power of deep learning to form a more complex policy. Together, these four benchmarks provide a comprehensive comparison, representing the leading non-parametric and parametric approaches in the field.

Table~\ref{tab:performance_comparison_3dec} summarizes the out-of-sample performance of all competing methods across three levels of contextual dimensionality: $Q \in \{1, 5, 20\}$. To provide a comprehensive evaluation of the models' sample complexity and data efficiency, we evaluate their performance under varying training sample sizes, specifically $N \in \{50, 100, 200, 400\}$. To ensure the statistical significance of our results, the reported metrics are based on 50 independent trials. Within each trial, we train every model on a randomly drawn training set of size $N$. The performance is subsequently assessed on a large, separately generated test set of 10,000 samples, allowing for a precise estimation of the true out-of-sample conditional expected loss. For reproducibility, the detailed model formulations, hyperparameter tuning procedures, and specifics of the normalizing flow training are provided in Appendix~\ref{subsec:inventory_appendix}.

The experimental results in Table~\ref{tab:performance_comparison_3dec} reveal several key findings regarding the models' scalability across varying sample sizes and contextual dimensions. First, our GMM-based frameworks demonstrate exceptional robustness, achieving the lowest average out-of-sample loss in all experimental settings. Furthermore, we observe that as the sample size increases, the GMM-NF model exhibits an even stronger performance. This trend clearly demonstrates that the normalizing flow, as a highly expressive deep generative architecture, is particularly powerful in learning complex conditional distributions under data-rich settings. Nevertheless, the performance of the original GMM model remains respectable, highlighting its inherent flexibility even under non-GM distributions. Second, while the RNW method is competitive in low-dimensional cases ($Q=1$), its performance degrades precipitously as the contextual dimension increases. This phenomenon aligns with the ``curse of dimensionality" issue discussed earlier. Finally, parametric prediction methods (ResDRO) and policy optimization approaches (LDR and NN-DR) struggle to navigate the complex mixture distributions. Notably, even though the NN-DR baseline utilizes deep neural networks to enhance its representational capacity, it consistently falls short of the GMM-NF framework. This observation underscores that simply deploying deep neural networks in contextual optimization is not a panacea. Instead, meticulous structural design, such as our separate NF structure and robustification scheme, is essential to fully unlock the power of deep learning.

In summary, our proposed GMM framework, particularly when enhanced with normalizing flows, establishes its superiority and robustness across all tested scenarios. It effectively overcomes the curse of dimensionality that plagues the non-parametric benchmark and shows remarkable flexibility compared with parametric and policy optimization models.

\begin{table}[htbp]
\centering
\begin{adjustbox}{width=\textwidth}
\begin{tabular}{ll rrr rrr rrr}
\toprule
\multirow{2}{*}{\textbf{N}} & \multirow{2}{*}{\textbf{Method}}
& \multicolumn{3}{c}{\textbf{$Q=1$}}
& \multicolumn{3}{c}{\textbf{$Q=5$}}
& \multicolumn{3}{c}{\textbf{$Q=20$}} \\
\cmidrule(lr){3-5} \cmidrule(lr){6-8} \cmidrule(lr){9-11}
& & mean & $P_{10}$-th & $P_{90}$-th
  & mean & $P_{10}$-th & $P_{90}$-th
  & mean & $P_{10}$-th & $P_{90}$-th \\
\midrule
\multirow{6}{*}{50}
& GMM-NF     & \textbf{12.260} & 4.547 & \textbf{21.834} & 26.017 & 7.678 & \textbf{57.216} & 160.406 & 53.282 & 314.498 \\
& GMM & 12.817 & \textbf{4.459} & 22.342 & \textbf{25.326} & 7.847 & 57.675 & \textbf{159.448} & \textbf{38.698} & 306.120 \\
& ResDRO     & 16.566 & 6.943 & 23.052 & 34.383 & 14.741 & 59.920 & 192.064 & 86.651 & 337.304 \\
& RNW        & 15.273 & 6.763 & 23.780 & 26.869 & \textbf{6.936} & 63.551 & 234.931 & 178.753 & 293.519 \\
& LDR        & 16.539 & 6.364 & 28.286 & 30.031 & 8.968 & 65.579 & 281.934 & 102.599 & 653.420 \\
& NN-DR      & 54.867 & 43.340 & 64.568 & 68.187 & 36.195 & 109.626 & 168.362 & 79.780 & \textbf{267.252} \\
\midrule
\multirow{6}{*}{100}
& GMM-NF     & \textbf{9.736} & 4.127 & \textbf{15.471} & \textbf{22.686} & 6.238 & 58.142 & 133.568 & \textbf{28.535} & 261.308 \\
& GMM & 10.745 & \textbf{4.025} & 16.780 & 23.597 & 6.576 & 57.240 & \textbf{122.810} & 29.414 & \textbf{244.716} \\
& ResDRO     & 16.752 & 6.666 & 23.717 & 29.884 & 10.942 & \textbf{54.084} & 184.617 & 95.553 & 301.145 \\
& RNW        & 15.159 & 6.725 & 21.661 & 26.043 & \textbf{6.223} & 64.354 & 238.851 & 175.894 & 296.181 \\
& LDR        & 16.436 & 5.828 & 29.221 & 26.306 & 7.067 & 61.749 & 172.150 & 76.119 & 298.702 \\
& NN-DR      & 12.544 & 4.075 & 22.048 & 42.140 & 13.802 & 78.871 & 156.608 & 59.492 & 262.902 \\
\midrule
\multirow{6}{*}{200}
& GMM-NF     & \textbf{8.471} & \textbf{3.791} & \textbf{14.901} & \textbf{22.949} & 6.419 & \textbf{57.932} & \textbf{126.193} & \textbf{19.884} & 246.358 \\
& GMM & 10.725 & 3.976 & 19.462 & 23.338 & 6.961 & 58.505 & 133.023 & 30.495 & \textbf{235.075} \\
& ResDRO     & 16.263 & 6.973 & 22.691 & 30.954 & 13.405 & 58.400 & 166.981 & 92.871 & 292.989 \\
& RNW        & 14.555 & 7.715 & 19.806 & 26.349 & \textbf{6.198} & 64.124 & 244.610 & 176.879 & 309.386 \\
& LDR        & 16.433 & 6.347 & 28.784 & 25.582 & 7.043 & 62.222 & 153.120 & 60.439 & 271.033 \\
& NN-DR      & 10.316 & 3.860 & 16.112 & 29.154 & 10.567 & 60.424 & 136.577 & 41.762 & 241.497 \\
\midrule
\multirow{6}{*}{400}
& GMM-NF     & \textbf{7.891} & \textbf{3.741} & \textbf{14.443} & \textbf{22.499} & \textbf{6.152} & 57.618 & \textbf{101.671} & \textbf{18.644} & \textbf{202.233} \\
& GMM & 10.018 & 3.938 & 16.046 & 22.533 & 7.009 & \textbf{50.606} & 110.498 & 26.618 & 221.192 \\
& ResDRO     & 15.793 & 7.400 & 20.863 & 29.349 & 13.356 & 56.493 & 151.186 & 69.332 & 251.526 \\
& RNW        & 13.963 & 7.441 & 18.714 & 26.078 & 6.175 & 66.502 & 122.534 & 20.512 & 262.059 \\
& LDR        & 15.642 & 6.461 & 24.963 & 25.785 & 6.320 & 62.011 & 136.657 & 25.310 & 257.867 \\
& NN-DR      & 9.185 & 4.047 & 14.927 & 27.872 & 10.266 & 57.966 & 132.456 & 30.841 & 235.752 \\
\bottomrule
\end{tabular}
\end{adjustbox}
\caption{Performance comparison of different approaches across varying covariate dimensions ($Q$) and training sample sizes ($N$).}
\label{tab:performance_comparison_3dec}
\end{table}

\subsection{Portfolio Optimization}

In the second experiment, we address the problem of contextual portfolio optimization, a critical task in financial decision-making. To ground our analysis in a state-of-the-art context, we follow the experimental design as in~\citet{nguyen2024robustifying}. This alignment ensures that all methods are tested under identical market conditions, facilitating a fair and rigorous comparison.  For the sake of completeness, we herein detail the full experimental setup.

Following the work of~\cite{NEURIPS2022_3df87436,nguyen2024robustifying}, we adopt the following five indices as contextual information: (i) Volatility Index (VIX), (ii) 10-year Treasury Yield Index (TNX), (iii) Crude Oil Index (CL=F), (iv) S$\&$P 500 (GSPC), (v) Dow Jones Index (DJI).\footnote{The five contextual indices are downloaded from Yahoo Finance. The daily stock-price files used to construct the return panel are downloaded from the Kaggle dataset \url{https://www.kaggle.com/datasets/joebeachcapital/s-and-p500-index-stocks-daily-updated?resource=download}. The time-varying S\&P 500 membership mask is constructed from WRDS/CRSP using \texttt{crsp.msp500list} and \texttt{crsp.msenames}.} To ensure implementability and avoid look-ahead bias, each contextual variable is lagged by one trading day. In particular, the portfolio decision for day $t$ is based only on closing market information observed at day $t-1$.

However, deviating from \cite{nguyen2024robustifying}, our investable universe is updated dynamically at the start of each quarter and held fixed for the quarter's duration. On the first trading day of each quarter, we restrict the universe to current S\&P 500 constituents with complete return observations, thereby avoiding look-ahead bias. We train our GMM-NF window using a 2-year rolling training set. To calibrate the models, the final month of this 2-year set serves as the validation period for hyperparameter tuning. Over our out-of-sample testing period from January 2022 to December 2023, we compute daily portfolio decisions by feeding the daily observed contextual side information into our models to update the conditional distribution estimations.

We consider a mean--CVaR portfolio problem with tail-risk level $\tau=0.1$:
\[
\text{CVaR}_\PP^{1-\tau}[- \bm x^\top \txi \mid \tilde \s = \s] - \eta \EE_\PP[\bm x^\top \txi \mid \tilde \s = \s].
\]
The parameter $\eta$ controls the trade-off between downside protection and return-seeking behavior, and we evaluate model performance for $\eta \in \{1,3,5,7,9\}$. Our proposed methods are benchmarked against a comprehensive suite of established models in portfolio optimization and contextual portfolio optimization literature: (i) the equal-weighted portfolio (EW), (ii) the unconditional Mean-CVaR model (MC), (iii) the Distributionally Robust unconditional Mean-CVaR model (DRMC)~\citep{blanchet2019quantifying,mohajerin2018data}, (iv) the Conditional Mean-CVaR model (CMC)~\citep{nguyen2024robustifying}, (v) the Distributionally Robust Conditional Mean-CVaR model (DRCMC)~\citep{nguyen2020distributionally}, (vi) the Optimal Transport distributionally robust Conditional Mean-CVaR model (OTCMC)~\citep{nguyen2024robustifying}, and (vii) our data-driven contextual optimization method based on Gaussian mixtures with normalizing flows (GMM-NF). The benchmark implementations are based on publicly available code,\footnote{\url{https://github.com/shanshanwang2019/Robustifying-Conditional-Portfolio-Decisions-via-Optimal-Transport}} and other implementation details are reported in Appendix~\ref{subsec:portfolio_appendix}.

The out-of-sample performance of all methods is summarized in Table~\ref{tab:model_comparison}. The results compellingly demonstrate the superiority of our proposed framework. Specifically, for $\eta=3,5$ and $9$, our method achieves the best performance, evidenced by a significantly lower out-of-sample loss and a higher Sharpe ratio compared to the benchmark approaches. Meanwhile, for $\eta=1$, our method also delivers close-to-best performance. This robust outperformance suggests that our framework is adept at modeling the conditional distribution of asset returns, leading to more effective portfolio allocation decisions across all risk preference levels.

\begin{table}[h]
\centering
\begin{tabular}{c l c c c c}
\toprule
$\eta$ & model & risk $\downarrow$ & mean $\uparrow$ & CVaR $\downarrow$ & Sharpe $\uparrow$ \\
\midrule
\multirow{7}{*}{1} 
& EW & 2.084 & \textbf{0.017} & 2.101 & \textbf{0.225} \\
& MC & \textbf{1.664} & -0.015 & \textbf{1.649} & -0.271 \\
& DRMC & 1.992 & 0.012 & 2.004 & 0.166 \\
& CMC & 1.780 & -0.019 & 1.760 & -0.325 \\
& DRCMC & 2.054 & 0.010 & 2.064 & 0.138 \\
& OTCMC & 2.003 & 0.004 & 2.007 & 0.064 \\
& GMM-NF & 1.692 & -0.009 & 1.683 & -0.159 \\
\midrule
\multirow{7}{*}{3}
& EW & 2.050 & 0.017 & 2.101 & 0.225 \\
& MC & 1.948 & -0.031 & 1.856 & -0.500 \\
& DRMC & 2.000 & 0.012 & 2.037 & 0.171 \\
& CMC & 2.021 & -0.018 & 1.967 & -0.264 \\
& DRCMC & 2.033 & 0.010 & 2.064 & 0.141 \\
& OTCMC & 1.972 & 0.012 & 2.007 & 0.163 \\
& GMM-NF & \textbf{1.729} & \textbf{0.018} & \textbf{1.784} & \textbf{0.289} \\
\midrule
\multirow{7}{*}{5}
& EW & 2.016 & 0.017 & 2.101 & 0.225 \\
& MC & 2.041 & -0.010 & 1.992 & -0.149 \\
& DRMC & 1.993 & 0.013 & 2.058 & 0.176 \\
& CMC & 2.470 & -0.033 & 2.307 & -0.410 \\
& DRCMC & 2.009 & 0.011 & 2.064 & 0.152 \\
& OTCMC & 1.954 & 0.012 & 2.013 & 0.166 \\
& GMM-NF & \textbf{1.709} & \textbf{0.018} & \textbf{1.800} & \textbf{0.283} \\
\midrule
\multirow{7}{*}{7}
& EW & 1.982 & \textbf{0.017} & 2.101 & \textbf{0.225} \\
& MC & 2.216 & -0.003 & 2.193 & -0.044 \\
& DRMC & 1.975 & 0.014 & 2.074 & 0.192 \\
& CMC & 2.797 & -0.015 & 2.693 & -0.157 \\
& DRCMC & 1.977 & 0.011 & 2.056 & 0.155 \\
& OTCMC & 1.934 & 0.013 & 2.023 & 0.178 \\
& GMM-NF & \textbf{1.772} & 0.013 & \textbf{1.864} & 0.201 \\
\midrule
\multirow{7}{*}{9}
& EW & 1.948 & 0.017 & 2.101 & 0.225 \\
& MC & 2.482 & 0.006 & 2.538 & 0.072 \\
& DRMC & 1.950 & 0.014 & 2.081 & 0.195 \\
& CMC & 3.138 & -0.024 & 2.921 & -0.239 \\
& DRCMC & 1.948 & 0.012 & 2.060 & 0.170 \\
& OTCMC & 1.904 & 0.014 & 2.030 & 0.194 \\
& GMM-NF & \textbf{1.586} & \textbf{0.026} & \textbf{1.819} & \textbf{0.389} \\
\midrule
\bottomrule
\end{tabular}
\caption{Performance metrics for different contextual portfolio optimization models under varying $\eta$ values.}
\label{tab:model_comparison}
\end{table}

\subsection{Multi-stage Wind Energy Optimization}

In our final experiment, we extend our study to the multistage setting. Following~\cite{park2022data} and \cite{kim2011optimal}, we investigate a multistage energy planning problem in the day-ahead market. At the beginning of day~$t$, the producer observes the day-ahead hourly prices $\bm{p}_t \in \mathbb{R}^{24}_+$ and determines the commitment levels $\bm{u}_t \in \mathbb{R}^{24}_+$ for the next day’s production. On day $t{+}1$, the commitments are fulfilled either through generation $\tilde \bxi_{t+1}$ or by discharging from the storage units. The unmet commitments incur a penalty of twice the day-ahead price, while the excess generations are used to recharge storage. When the storage capacity is saturated, the surplus energy is curtailed.

We obtain the hourly wind energy data from the North American Land Data Assimilation System~\citep{xia2013overview} from $2002$ to $2011$ at \textbf{Ohio} and \textbf{North Carolina} regions. The day-ahead prices are publicly available at the PJM market\footnote{Day-ahead price data can be downloaded from PJM: \url{http://dataminer2.pjm.com/feed/da_hrl_lmps/definition}}. We construct \emph{weekly trajectories} consisting of seven consecutive days (24 hours $\times$ 7 stages), yielding 520 trajectories in total. To capture seasonal variation, we partition the decision-making problem into quarterly subproblems; in each quarter, we take the first three years of data as in-sample trajectories and evaluate out-of-sample performance on the remainder. The implementation details are presented in Appendix~\ref{subsec:Wind_energy_appendix}.

We integrate our proposed GMM-DRO transition model into the stochastic dual dynamic programming (SDDP) framework. By iteratively constructing piecewise-linear cuts to approximate the expected cost-to-go function, SDDP is widely viewed as one of the most effective approaches for solving multistage stochastic programming problems~\citep{fullner2025stochastic}. To evaluate our method, we benchmark it against three prominent methods from the literature: stagewise independent SDDP (IND-SDDP)~\citep{shapiro2011analysis}, Neural SDDP (N-SDDP)~\citep{dai2021neural}, and distributionally robust Nadaraya-Watson SDDP (NW-SDDP)~\citep{park2022data}. The first method, IND-SDDP, assumes no temporal dependence between different stages and therefore ignores the Markovian structure of the uncertainty process. The second approach, N-SDDP, integrates deep learning with SDDP to effectively enhance solution efficiency. It utilizes neural networks to predict affine cuts that approximate the cost-to-go functions, serving as a warm-start mechanism that is subsequently embedded and refined within the standard SDDP iterations. The third benchmark, NW-SDDP, represents a recent data-driven advancement for handling Markovian uncertainty. This method utilizes Nadaraya-Watson kernel regression to estimate transition probabilities and applies a DRO scheme to robustify the policy. Our method instead fits a global joint GM model for consecutive-stage uncertainties and uses the induced likelihood-ratio weights, as described in~\eqref{eq:apprx_exptected_ctg_new}. 

\begin{table}[h]
\centering
\begin{tabular}{ l l c c c }
\toprule
 Data set & Model & mean $\uparrow$ & $10^{\text{th}}$ percentile $\uparrow$ & $90^{\text{th}}$ percentile $\uparrow$ \\
\midrule
 \multirow{4}{*}{Ohio} 
 & IND-SDDP & 5.699 & 1.074 & 11.563 \\
 & N-SDDP  & 5.593 & 0.735 & 11.671 \\
 & NW-SDDP  & 5.866 & 1.547 & 11.437 \\
 & GMM-SDDP & \textbf{5.927} & \textbf{1.650} & \textbf{11.677} \\
\midrule
 \multirow{4}{*}{North Carolina} 
 & IND-SDDP & 7.604 & 1.419 & 14.069 \\
 & N-SDDP  & 7.071 & 0.392 & 14.031 \\
 & NW-SDDP  & 8.011 & 1.935 & \textbf{15.167} \\
 & GMM-SDDP & \textbf{8.479} & \textbf{2.539} & 14.664 \\
\bottomrule
\end{tabular}
\caption{Out-of-sample performance of different methods in the multistage day-ahead wind energy planning problem (in \$100{,}000).}
\label{tab:comparison3}
\end{table}

Table~\ref{tab:comparison3} summarizes the out-of-sample performance of different methods. First, the NW-SDDP method outperforms IND-SDDP by effectively exploiting the Markovian structure of the underlying uncertainty. This observation underscores the importance of capturing temporal dependencies in multistage settings. However, because the NW-SDDP method relies on local non-parametric weights, it inevitably suffers from the curse of dimensionality, which limits its performance in our high-dimensional problem instances. In contrast, our proposed GMM-SDDP approach offers a structured parametric representation of the dependence between consecutive daily power-price states, which is more stable in this high-dimensional problem.  Benefited from this nice property, it dominates the benchmark methods on average and at the $10^{th}$-percentile across both datasets. Finally, although the N-SDDP model exhibits superior computational efficiency (with detailed running times reported in Table~\ref{tab:runtime_sddp} in Appendix \ref{subsec:Wind_energy_appendix}), its overall solution quality is inferior to that of the other benchmark methods.

\section{Concluding Remarks}

This paper introduced a novel framework for data-driven contextual optimization, designed to address a fundamental trade-off in existing solution schemes. Prevailing parametric methods often lack the flexibility to handle complex, multimodal uncertainties, while non-parametric approaches are notoriously hampered by the curse of dimensionality. Our GMM-based framework provides an easy-to-implement yet powerful methodology that elegantly bridges this gap. By offering a structured yet flexible estimate of the conditional distribution, this approach is highly suitable for high-dimensional mixture distributions where traditional non-parametric methods fail and standard parametric models lack sufficient expressive power. To relax strict distributional assumptions, we integrate our framework with NFs. This extension provides additional modeling flexibility that can improve the distributional approximation and the resulting optimized decisions. Additionally, we design a DRO variant that is particularly useful when historical data are limited and overfitting is a primary concern. 

Finally, we extended our framework to the domain of dynamic decision-making by designing a novel GMM-based solution scheme for multistage stochastic optimization problems under Markovian uncertainty. Our theoretical analysis confirms that this method achieves significantly better sample complexity than traditional approaches, providing a powerful new tool for solving long-horizon, high-dimensional multistage problems. Extensive numerical experiments on a series of financial and business decision problems demonstrate the practical superiority of our approach, as it consistently outperforms state-of-the-art benchmarks, especially in high-dimensional scenarios and when faced with complex data structures.

% \section*{Acknowledgments}
%G.~Hanasusanto are funded in part by the National Science Foundation under grants 2342505 and 2343869. 

% \newpage
\bibliographystyle{plainnat}
\bibliography{bibliography}

\clearpage
\begin{APPENDICES}
\setlength{\parskip}{1em}

\section{Proofs of Section~\ref{sec:contextual_GMM}}\label{sec:appendix_A}
\begin{proof}{Proof of Proposition~\ref{prop:conditional_param}}
We first consider the case where the random vectors $\tilde{\s}$ and $\txi$ are jointly Gaussian (i.e., $K=1$) with the density function $\mathcal{N}\left((\s,\bxi)| \bm{\mu},\bm{\Sigma} \right)$. The marginal densities of $\tilde{\bxi}$ and $\tilde{\s}$ are $\mathcal{N}(\bxi| \bm{\mu}_{\bxi}, \bm{\Sigma}_{\bxi\bxi})$ and $\mathcal{N}(\s|\bm{\mu}_{\s}, \bm{\Sigma}_{\s\s})$, respectively. The density function of the conditional distribution of $\txi$ given $\s$ is \cite[Section 2.3.2, Equations (2.94)-(2.98)]{bishop2006pattern}
\begin{equation*}
p(\bxi|\s)= \mathcal{N} \left(\bxi|\bm{\mu}_{\bxi|\s}, \bm{\Sigma}_{\bxi|\s}\right),
\end{equation*}
where $\bm{\mu}_{\bxi|\s} = \bm{\mu}_{\bxi} + \bm{\Sigma}_{\bxi\s}( \bm{\Sigma}_{\s\s})^{-1}(\s- \bm{\mu}_{\s})$ and $\bm{\Sigma}_{\bxi|\s} =\bm{\Sigma}_{\bxi\bxi}-\bm{\Sigma}_{\bxi\s}(\bm{\Sigma}_{\s\s})^{-1}\bm{\Sigma}_{\s\bxi}$.

We now extend this result to the case of a GM joint distribution, where the marginal density of $\tilde{\s}$ is
\begin{equation*}
p(\s)= \sum_{k=1}^K  p^k \mathcal{N}(\s|\bm{\mu}_{\s}^k,\bm{\Sigma}_{\s\s}^k). 
\end{equation*}
Here, the conditional density function becomes
\begin{equation*}
p(\bxi|\s) 
= \frac{p(\bxi,\s)}{p(\s)} 
= \sum_{k=1}^K\frac{  p^k \mathcal{N}\left( (\bxi,\s) | \bm{\mu}^k, \bm{\Sigma}^k \right)}{\sum_{j=1}^K p^j \mathcal{N}(\s| \bm{\mu}_{\s}^j, \bm{\Sigma}_{\s\s}^j)} 
= \sum_{k=1}^K \frac{  p^k \mathcal{N} \left(\s|\bm{\mu}_{\s}^k, \bm{\Sigma}_{\s\s}^k \right)}{\sum_{j=1}^K p^j \mathcal{N} (\s| \bm{\mu}_{\s}^j, \bm{\Sigma}_{\s\s}^j)} \mathcal{N} \left(\bxi\; | \; \bm{\mu}_{\bxi|\s}^k, \bm{\Sigma}_{\bxi|\s}^k\right) .
\end{equation*}
Thus, this is a GM distribution with components $\mathcal{N} (\bxi\; | \; \bm{\mu}_{\bxi|\s}^k, \bm{\Sigma}_{\bxi|\s}^k)$, $k\in[K]$, and mixture probabilities 
\begin{equation*}
\frac{  p^k \mathcal{N} \left(\s|\bm{\mu}_{\s}^k, \bm{\Sigma}_{\s\s}^k \right)}{\sum_{j=1}^K p^j \mathcal{N} (\s| \bm{\mu}_{\s}^j, \bm{\Sigma}_{\s\s}^j)}, ~ k\in[K],
\end{equation*}
which completes the proof. \qed
\end{proof}

\begin{proposition}
\label{prop:tv_imply_parameters}
Consider two Gaussians  $\N(\bmu,\bSigma)$ and $\Nhat(\hat\bmu,\hat\bSigma)$  in $\RR^D$. If their total variation distance $\mathds{TV}(\N(\bmu,\bSigma),\Nhat(\hat\bmu,\hat\bSigma)) \leq \eps$, and the covariance matrices $\bSigma, \hat{\bSigma} \preceq \beta \I$, then 
% \begin{equation*}
% \|\bSigma-\hat\bSigma\|^2\lessapprox \frac{16\beta^2}{D} \quad\textup{and}\quad \frac{1}{8\beta}\|\bmu-\hat\bmu\|\lessapprox{8\beta}.
% \end{equation*}
\begin{equation*}
\|\bSigma - \hat\bSigma\|_F \leq 4\beta \sqrt{\frac{\eps}{1-\eps}}, \quad \text{ and }  \quad \|\bmu-\hat\bmu\|_2 \leq \sqrt{\frac{8\beta\eps}{1-\eps}}.
\end{equation*}

\end{proposition}
\begin{proof}{Proof of \Cref{prop:tv_imply_parameters}.}
The total variation distance admits a lower bound in the squared Hellinger distance, which in the case of two Gaussians $\N(\bmu,\bSigma)$ and $\Nhat(\hat\bmu,\hat\bSigma)$ is given by 
\begin{equation}
\label{eq:Hellinger}
\mathds{H}^2(\N(\bmu,\bSigma),\Nhat(\hat\bmu,\hat\bSigma))\coloneqq 1-\frac{|\bSigma|^{\frac{1}{4}}|\hat\bSigma|^{\frac{1}{4}}}{\left|\bm\Omega\right|^{\frac{1}{2}}}\exp\left(-\frac{1}{8}(\bmu-\hat\bmu)^\top\bm\Omega^{-1}(\bmu-\hat\bmu)\right)
\leq \mathds{TV}(\N(\bmu,\bSigma),\Nhat(\hat\bmu,\hat\bSigma)),
\end{equation}
where $\bm\Omega\coloneqq \frac{\bSigma+\hat\bSigma}{2}$. 
We will further derive a lower bound on the squared Hellinger distance. To simplify the exposition, we denote $\bm\Delta\coloneqq \bSigma-\hat\bSigma$; hence $\bSigma=\bm\Omega+\bm\Delta/2$ and    $\hat\bSigma=\bm\Omega-\bm\Delta/2$. 
Define $\bm B\coloneq\frac{1}{2}\bm\Omega^{-\frac{1}{2}}\bm\Delta\bm\Omega^{-\frac{1}{2}}$ which is symmetrically normalized covariance difference, then $\bSigma = \bm\Omega+\bm\Omega^{\frac{1}{2}}\bm B\bm\Omega^{\frac{1}{2}}, \ \hat\bSigma = \bm\Omega-\bm\Omega^{\frac{1}{2}}\bm B\bm\Omega^{\frac{1}{2}}$. 
By the multiplicative property of the determinant, 
\begin{equation*}
\begin{aligned}
|\bSigma|=|\bm\Omega^{\frac{1}{2}}(\I+\bm B)\bm\Omega^{\frac{1}{2}}| = |\bm\Omega^{\frac{1}{2}}||\I +\bm B||\bm\Omega^{\frac{1}{2}}| = |\bm\Omega||\I + \bm B| \\
|\hat\bSigma|=|\bm\Omega^{\frac{1}{2}}(\I-\bm B)\bm\Omega^{\frac{1}{2}}| = |\bm\Omega^{\frac{1}{2}}||\I -\bm B||\bm\Omega^{\frac{1}{2}}| = |\bm\Omega||\I - \bm B|  
\end{aligned}
\end{equation*}
From these results, 
\begin{equation*}
    \frac{|\bSigma|^{\frac{1}{4}}|\hat\bSigma|^{\frac{1}{4}}}{\left|\bm\Omega\right|^{\frac{1}{2}}}=\frac{(|\bm \Omega||\I +\bm B|)^{\frac{1}{4}}(|\bm \Omega||\I -\bm B|)^\frac{1}{4}}{|\bm \Omega|^\frac{1}{2}} = (|\I+\bm B||\I-\bm B|)^\frac{1}{4} = (|\I - \bm B^2|)^\frac{1}{4} 
\end{equation*}
Substituting $\frac{|\bSigma|^{\frac{1}{4}}|\hat\bSigma|^{\frac{1}{4}}}{\left|\bm\Omega\right|^{\frac{1}{2}}}= R$, $R^4=|\bm I - \bm B^2|$.
Since $\bm B$ is a symmetric matrix and by using logarithm, $4\log R = \log (|\I - \bm B^2|)$, and $4\log R = \sum_{i=1}^d \log(1- \lambda_i(\bm B)^2)$, where $\lambda_i$ is a eigenvalue of $\bm B$. 
Also, $\bSigma \succ 0$ and $\hat\bSigma\succ0$, we have $\I+\bm B\succ 0$ and $\I-\bm B \succ 0$. Therefore, every eigenvalue $\lambda_i(\bm B)$ satisfies $-1<\lambda_i(\bm B)<1$.
Using $\log (1-x) \leq -x$ for $0\leq x<1$,
\begin{equation*}
\begin{aligned}
\sum_{i=1}^d \log\left(1- \lambda_i(\bm B)^2\right) \leq - \sum_{i=1}^d \lambda_i(\bm B)^2 \\ 
4\log R \leq -\sum_{i=1}^d \lambda_i(\bm B)^2 = -\|\bm B\|_F^2 \\
R \leq \exp\left(-\frac{1}{16}||\bm \Omega^{-\frac{1}{2}}\bm\Delta\bm\Omega^{-\frac{1}{2}}\|_F^2\right)
\end{aligned}
\end{equation*}
Then the lower bound of Hellinger distance is defined as
\begin{equation*}
\begin{aligned}
\mathds{H}^2(\N(\bmu,\bSigma),\Nhat(\hat\bmu,\hat\bSigma))\geq 1-\exp\left(-\frac{1}{16}\|\bm \Omega^{-\frac{1}{2}}\bm\Delta\bm\Omega^{-\frac{1}{2}}\|_F^2-\frac{1}{8}(\bmu-\hat\bmu)^\top\bm\Omega^{-1}(\bmu-\hat\bmu)\right)\\
\end{aligned}
\end{equation*}
From the assumption $\bSigma \preceq \beta\I$ and $\hat\bSigma \preceq \beta\I$, it follows that $\bm\Omega=\frac{\bSigma+\hat\bSigma}{2}\preceq \beta\I.$
Accordingly, we have $\bm\Omega^{-1}\succeq \frac{1}{\beta}\I$, and the smallest singular value of $\bm\Omega^{-1/2}$ satisfies
$\sigma_{\min}(\bm\Omega^{-1/2})\ge \beta^{-1/2}$. Moreover, we have $(\bmu-\hat\bmu)^\top \bm\Omega^{-1}(\bmu-\hat\bmu)\ge\frac{1}{\beta}\|\bmu-\hat\bmu\|_2^2.$
Using the above results, 
\begin{equation*}
\|\bm \Omega^{-\frac{1}{2}}\bm\Delta\bm\Omega^{-\frac{1}{2}}\|_F \geq \sigma_{\min}(\bm \Omega^{-\frac{1}{2}})^2||\bm\Delta||_F\geq\frac{1}{\beta}\|\bm\Delta\|_F.
\end{equation*}
Therefore, the lower bound of the squared Hellinger distance and the relation are satisfied as follows 
\begin{equation*}
1- \exp\left(-\frac{1}{16\beta^2}\|\bSigma - \hat\bSigma \|^2_F - \frac{1}{8\beta}\|\bmu-\hat\bmu\|^2_2\right) \leq \mathds{H}^2(\N(\bmu,\bSigma),\Nhat(\hat\bmu,\hat\bSigma)) \leq \mathds{TV}(\N(\bmu,\bSigma),\Nhat(\hat\bmu,\hat\bSigma)) \leq \eps.
\end{equation*}
Finally, by noting that $x/(1+x) \leq 1-\exp(-x)$ for all $x>-1$, it follows
\begin{equation*}
    \frac{1}{16\beta^2}\|\bSigma - \hat\bSigma \|^2_F + \frac{1}{8\beta}\|\bmu-\hat\bmu\|^2_2\leq \frac{\eps}{1-\eps},
\end{equation*}
which implies the conclusion and finishes the proof. \qed

\end{proof}

\begin{proof}{Proof of \Cref{thm:tv_conditional}}
    We define the function
    \begin{equation*}
\phi(\s)\;:=\;\int_{\R^{R}}\!\big|f_{\PP}(\s,\bxi)-f_{\mathbb Q}(\s,\bxi)\big|\, {\rm d}\bxi \qquad \forall \s\in\R^{Q}.
    \end{equation*}
     Then we show that the function $\phi(\s)$ is also H\"older continuous. Specifically, for all $\s_{1},\s_{2}\in\R^{Q}$:
\begin{equation*}
    \begin{array}{rl}
        |\phi(\s_{1})-\phi(\s_{2})|
&\le\!\displaystyle\int_{\R^{R}}\!\big|\,|f_{\PP}(\s_{1},\bxi)-f_{\mathbb Q}(\s_{1},\bxi)|
  -|f_{\PP}(\s_{2},\bxi)-f_{\mathbb Q}(\s_{2},\bxi)|\,\big|\,{\rm d}\bxi\\[10pt]
&\le\!\displaystyle\int_{\R^{R}}\!|f_{\PP}(\s_{1},\bxi)-f_{\PP}(\s_{2},\bxi)|\,{\rm d}\bxi
  +\!\int_{\R^{R}}\!|f_{\mathbb Q}(\s_{1},\bxi)-f_{\mathbb Q}(\s_{2},\bxi)|\,{\rm d}\bxi,
    \end{array}
\end{equation*}
where the last inequality is from the reverse triangle inequality $\big| |u| - |v|\big| \leq |u-v|$. Since both $f_{\PP}$ and $f_{\mathbb Q}$ are with H\"older continuous slice functions, it follows that $\phi(\s)$ is also H\"older continuous:
\begin{equation}\label{eqn:phi_holder}
    |\phi(\s_{1})-\phi(\s_{2})| \leq \left( K_{\PP} + K_{\mathbb Q}\right) \left\| \s_1 - \s_2\right\|^\alpha \quad \forall \s_1, \; \s_2 \in \R^Q.
\end{equation}

Now we fix any side information $\s\in \R^{Q}$ and bound $\mathds{TV}(\PP_{\bxi | \s},\mathbb Q_{\bxi | \s})$. We note that
\begin{equation*}
\mathds{TV}(\PP_{\bxi\mid\s},\mathbb Q_{\bxi\mid\s}) = \frac{1}{2}\!\int_{\R^{R}}\!\big|f_{\PP}(\bxi\mid\s)-f_{\mathbb Q}(\bxi\mid\s)\big|\,{\rm d}\bxi.
\end{equation*}
For any $a,b\ge 0$ and $x,y>0$, the algebraic identity
$\tfrac{a}{x}-\tfrac{b}{y}=\tfrac{a-b}{x}+\tfrac{b(y-x)}{xy}$ gives
$\big|\tfrac{a}{x}-\tfrac{b}{y}\big|\le\tfrac{|a-b|}{x}+\tfrac{b}{xy}|x-y|$.
Noting that the conditional density $f_{\PP}(\bxi \mid \s) = f_{\PP}(\s, \bxi) / f_{\PP}(\s)$ and similarly for $\mathbb Q$, applying the above algebraic inequality with $a=f_{\PP}(\s,\bxi)$, $b=f_{\mathbb Q}(\s,\bxi)$, $x=f_{\PP}(\s)$, $y=f_{\mathbb Q}(\s)$ and
integrating in $\bxi$:
\begin{equation*}
    \mathds{TV}(\PP_{\bxi\mid\s},\mathbb Q_{\bxi\mid\s}) \leq \frac{1}{2}\left(\frac{\phi(\s)}{f_{\PP}(\s)}\;+\;\frac{|f_{\PP}(\s)-f_{\mathbb Q}(\s)|}{f_{\PP}(\s)}\right)
\end{equation*}
     By the triangle inequality for the Lebesgue integral,
     \begin{equation*}
         |f_{\PP}(\s)-f_{\mathbb Q}(\s)| =\Big|\!\int_{\R^{R}}\!\big(f_{\PP}(\s,\bxi)-f_{\mathbb Q}(\s,\bxi)\big)\,{\rm d}\bxi\Big| \le\phi(\s).
     \end{equation*}
     Therefore, it implies that 
     \begin{equation}\label{eqn:tv_phi}
         \mathds{TV}(\PP_{\bxi\mid\s},\mathbb Q_{\bxi\mid\s}) \leq  \frac{\phi(\s)}{f_{\PP}(\s)}.
     \end{equation}
     Next, we bound $\phi(\s)$. Let $r>0$ be an index that is to be optimized later. Denote by $B_{r}=B_{r}(\s)\subset\R^{Q}$ the Euclidean ball of radius $r$ centred
at $\s$, with volume $|B_{r}|=V_{Q}\,r^{Q}$. Adding and subtracting
the average of $\phi$ over $B_{r}$,
\begin{equation}\label{eq:avg-trick}
\phi(\s)\;=\;\frac{1}{|B_{r}|}\!\int_{B_{r}}\!\phi(\s)\, {\rm d}\s^\prime
\;\le\;\frac{1}{|B_{r}|}\!\int_{B_{r}}\!\phi(\s^\prime)\,{\rm d}\s^\prime
\;+\;\frac{1}{|B_{r}|}\!\int_{B_{r}}\!|\phi(\s)-\phi(\s^\prime)|\,{\rm d}\s^\prime.
\end{equation}
The first term is bounded using $\phi\ge 0$ and Fubini's theorem:
\begin{equation*}
    \frac{1}{|B_{r}|}\!\int_{B_{r}}\!\phi(\s^\prime)\,{\rm d}\s^\prime
\;\le\;\frac{1}{|B_{r}|}\!\int_{\R^{Q}}\!\phi(\s^\prime)\,{\rm d}\s^\prime
\; = \;\frac{2\,\mathds{TV}(\PP,\mathbb Q)}{V_{Q}\,r^{Q}} \;\le\;\frac{\sqrt{2\epsilon}}{V_{Q}\,r^{Q}},
\end{equation*}
where the equality is due to $\int_{\R^{Q}}\phi(\s^\prime)\,{\rm d}\s^\prime =\int_{\R^{Q+R}}|f_{\PP}(\s^\prime,\bxi)-f_{\mathbb Q}(\s^\prime,\bxi)|\,{\rm d}\s^\prime\,{\rm d}\bxi =2\mathds{TV}(\PP,\mathbb Q)$, and the last inequality is from the consequence of Pinsker's inequality $\mathds{TV}(\PP,\mathbb Q) \leq \sqrt{\DD\left(f_{\PP} \big|\big|f_{\mathbb Q}\right)/2} \leq \sqrt{\eps/2}$.  

The second term in \eqref{eq:avg-trick}, according to \eqref{eqn:phi_holder}, is bounded by
\begin{equation*}
    |\phi(\s)-\phi(\s^\prime)| \le \left( K_{\PP} + K_{\mathbb Q}\right) \left\| \s - \s^\prime \right\|^\alpha \le \left( K_{\PP} + K_{\mathbb Q}\right) r^\alpha \qquad \forall \s^\prime \in B(r).
\end{equation*}

Combining these in~\eqref{eq:avg-trick} yields
\begin{equation}\label{eq:phi-bound}
\phi(\s)\; \le\;\frac{\sqrt{2\epsilon}}{V_{Q}\,r^{Q}}\;+\;\left( K_{\PP} + K_{\mathbb Q}\right) r^\alpha 
\qquad\text{for every }r>0.
\end{equation}

Noting that the above inequality is true for any $r>0$, we define the function
     \begin{equation*}
         h(r) := \frac{\sqrt{2\epsilon}}{V_{Q}\,r^{Q}}\;+\;\left( K_{\PP} + K_{\mathbb Q}\right) r^\alpha  \qquad \forall r>0.
     \end{equation*}
     By taking the derivative, it follows that the minimum of $h$ is attained at
     \begin{equation*}
         r^\star= \left(\frac{Q\sqrt{2\eps}}{V_{Q}\,\alpha\;\left( K_{\PP} + K_{\mathbb Q}\right)}\right)^{1/(Q+\alpha)}.
     \end{equation*}
     
     Substituting $r^\star$ to \eqref{eq:phi-bound}, we obtain
     \begin{equation*}
         \begin{array}{rl}
              \phi(\s) \leq & \displaystyle\frac{\sqrt{2\eps}}{V_Q}\left(\frac{V_{Q}\,\alpha\;\left( K_{\PP} + K_{\mathbb Q}\right)}{Q\sqrt{2\eps}}\right)^{\frac{Q}{Q+\alpha}} + \left( K_{\PP} + K_{\mathbb Q}\right) \left(\frac{Q\sqrt{2\eps}}{V_{Q}\,\alpha\;\left( K_{\PP} + K_{\mathbb Q}\right)}\right)^{\frac{\alpha}{Q+\alpha}}\\
              \leq & \left(\left(\displaystyle\frac{\alpha}{Q}\right)^{\frac{Q}{Q+\alpha}} + \left(\displaystyle\frac{Q}{\alpha}\right)^{\frac{\alpha}{Q+\alpha}}\right) \left( K_{\PP} + K_{\mathbb Q}\right)^{\frac{Q}{Q+\alpha}}\left(\displaystyle\frac{\sqrt{2\eps}}{V_Q}\right)^{\frac{\alpha}{Q+\alpha}}.
         \end{array}
     \end{equation*}
     The above inequality, together with \eqref{eqn:tv_phi}, implies the conclusion.
     \qed
\end{proof}

Next, we present the proof of Theorem~\ref{thm:approx_error}. The proof of Theorem~\ref{thm:approx_error} relies on several results.  
\begin{lemma}
\label{lem:ratio_bound}
 Let $a,a'\in\RR$, $b,b'\in\RR_+$ and suppose that $b,b'\geq\underline b>0$. If $|a-a'|\leq \epsilon_a$ and $|b-b'|\leq \epsilon_b$ then 
 \begin{equation*}
\left|\frac{a}{b}-\frac{a'}{b'}\right|\leq \frac{\epsilon_a}{b} + \frac{\epsilon_b|a'|}{bb'}\leq  \frac{\epsilon_a}{\underline b} + \frac{\epsilon_b|a'|}{\underline b^2}. 
\end{equation*}

\end{lemma}
\begin{proof}{Proof of Lemma~\ref{lem:ratio_bound}}
We have 
\begin{equation*}
\left|\frac{a}{b}-\frac{a'}{b'}\right|=\left|\frac{ab'-a'b}{bb'}\right|\leq \frac{|ab'-a'b'|+|a'b'-a'b|}{bb'}\leq  \frac{\epsilon_a}{b} + \frac{\epsilon_b|a'|}{bb'} \leq  \frac{\epsilon_a}{\underline b} + \frac{\epsilon_b|a'|}{\underline b^2} ,
\end{equation*}
which proves the claim. \qed
\end{proof}

\begin{lemma}\label{lem:matrixinv2}
Suppose the matrices $\bm A$ and $\bm B$ are strictly positive definite with $\alpha \bm I \preceq \bm A$, $\alpha \bm I \preceq \bm B$ for $\alpha\in\RR_{++}$. If $\|\bm B-\bm A\|\leq \eps$, then 
\[\| \bm A^{-1}-\bm B^{-1}\| \leq \frac{\epsilon}{\alpha^2}
\]
\end{lemma}
\begin{proof}{Proof of Lemma \ref{lem:matrixinv2}.}
  Since $\bm A^{-1} - \bm B^{-1} = \bm A^{-1}(\bm B - \bm A) \bm B^{-1}$,
we have
\begin{equation*}
\|\bm A^{-1} - \bm B^{-1} \| = \|\bm A^{-1}(\bm B - \bm A) \bm B^{-1} \|
\leq\|\bm A^{-1}\| \|(\bm B - \bm A)\| \| \bm B^{-1}  \|
\leq\frac{\epsilon}{\alpha^2}.
\end{equation*}
Thus, the claim follows.   
\end{proof}

\begin{lemma}\label{lem:normalpdfbound} 
Suppose $\bm \mu, \hat{\bm \mu} \in \R^D$ and $\bm \Sigma,  \hat\bSigma \in \mathbb S^D_{++}$ with $\alpha \I \preceq \bm \Sigma,  \hat\bSigma$. If 
\[\|{\bm \mu}-\hat{\bm\mu} \| \leq \epsilon_{\bmu}, \|{\bm\Sigma} - \hat{\bm\Sigma} \| \leq \epsilon_{\bSigma},
\]
then for every $\bm z\in\mathbb{R}^D$, we have 
\begin{align*}
    &\;\;\left|\mathcal{N}\left( \z| \bm{\mu}, \bm{\Sigma}\right)-\mathcal{N}\left( \z| \hat{\bm\mu}, \hat{\bm{\Sigma}}\right)\right| \\
    \leq&\;\;\mathcal{N}\left( \z| \bm{\mu}, \bm{\Sigma}\right)\left(\epsilon_\bmu \|\bm\Sigma^{-1}(\z-\bm\mu)\| + \frac{D\epsilon_\bSigma}{2} \left\|\bSigma^{-1}(\z-\bm\mu)(\z-\bm\mu)^\top\bSigma^{-1}-\bSigma^{-1}\right\| + + o(\epsilon_\bmu+\epsilon_\bSigma)\right) \\
    \leq&\;\;\mathcal{N}\left( \z| \bm{\mu}, \bm{\Sigma}\right)\left(\frac{\epsilon_\bmu}{\alpha} \left(\|\z\|+\gamma\right) + \frac{D\epsilon_\bSigma}{2\alpha^2}\left((\|\z\|+\gamma)^2+\alpha\right)  + o(\epsilon_\bmu+\epsilon_\bSigma)\right). 
\end{align*}
\end{lemma}
\begin{proof}{Proof of Lemma \ref{lem:normalpdfbound}.}
Recall that the normal density function is given by
\[\mathcal{N}\left( \z| \bm{\mu}, \bm{\Sigma}\right)=\frac{\exp\left(-\frac{1}{2}(\z-\bm\mu)^\top\bm\Sigma^{-1}(\z-\bm\mu)\right)}{\sqrt{(2\pi)^D|\bm \Sigma|}}.
\]
We first construct a first-order approximation around $(\bmu,\bSigma)$. The gradient with respect to the mean vector is
\begin{align*}
\nabla_\bmu \mathcal{N}\left( \z| \bm{\mu}, \bm{\Sigma}\right)&=\frac{\exp\left(-\frac{1}{2}(\z-\bm\mu)^\top\bm\Sigma^{-1}(\z-\bm\mu)\right)}{\sqrt{(2\pi)^D|\bm \Sigma|}} \bm\Sigma^{-1}(\z-\bm\mu)= \mathcal{N}\left( \z| \bm{\mu}, \bm{\Sigma}\right)\bm\Sigma^{-1}(\z-\bm\mu).
\end{align*}
To obtain the gradient with respect to the covariance matrix, we first take the logarithm: 
\begin{align*}
\log \mathcal{N}\left( \z| \bm{\mu}, \bm{\Sigma}\right)&=-\frac{1}{2}(\z-\bm\mu)^\top\bm\Sigma^{-1}(\z-\bm\mu) - \log \sqrt{(2\pi)^D} -\frac{1}{2} \log{|\bm \Sigma|}.
\end{align*}
Using \cite[Equation (63)]{petersen2008matrix} to compute the gradient of the quadratic term
\begin{equation*}
\nabla_\bSigma (\z-\bm\mu)^\top\bSigma^{-1}(\z-\bm\mu) = -\bSigma^{-1}(\bm{z}-\bm\mu)(\z-\bm\mu)^\top\bSigma^{-1}
\end{equation*}
and \cite[Equation (57)]{petersen2008matrix} to compute the gradient of the logarithm term 
\begin{equation*}
\nabla_\bSigma\log{|\bm \Sigma|} = \bSigma^{-1},
\end{equation*}
we find that 
\begin{align*}
\nabla_\bSigma \log \mathcal{N}\left( \z| \bm{\mu}, \bm{\Sigma}\right)&=\frac{1}{2}\bSigma^{-1}(\z-\bm\mu)(\z-\bm\mu)^\top\bSigma^{-1}-\frac{1}{2} \bSigma^{-1}.
\end{align*}
By the chain rule, we have 
\begin{equation*}
\nabla_\bSigma \mathcal{N}\left( \z| \bm{\mu}, \bm{\Sigma}\right) = \mathcal{N}\left( \z| \bm{\mu}, \bm{\Sigma}\right) \nabla_\bSigma \log \mathcal{N}\left( \z| \bm{\mu}, \bm{\Sigma}\right)= \mathcal{N}\left( \z| \bm{\mu}, \bm{\Sigma}\right) \left(\frac{1}{2}\bSigma^{-1}(\z-\bm\mu)(\z-\bm\mu)^\top\bSigma^{-1}-\frac{1}{2} \bSigma^{-1}\right).
\end{equation*}
 Hence, the Taylor expansion around $(\bmu,\bSigma)$, $$\mathcal{N}\left( \z| \hat{\bm\mu}, \hat{\bm{\Sigma}}\right)=\mathcal{N}\left( \z| \bm{\mu}, \bm{\Sigma}\right)- \nabla_\bmu \mathcal{N}\left( \z| \bm{\mu}, \bm{\Sigma}\right)^\top (\bmu-\hat\bmu) -\tr\left(\nabla_\bSigma \mathcal{N}\left( \z| \bm{\mu}, \bm{\Sigma}\right) (\bSigma-\hat\bSigma)\right) + R_1(\hat\bmu,\hat\bSigma)$$  yields the upper bound 
\begin{align*}
&\left|\mathcal{N}\left( \z| \bm{\mu}, \bm{\Sigma}\right)-\mathcal{N}\left( \z| \hat{\bm\mu}, \hat{\bm{\Sigma}}\right)\right|\\
&\leq  \left|\nabla_\bmu \mathcal{N}\left( \z| \bm{\mu}, \bm{\Sigma}\right)^\top (\bmu-\hat\bmu)\right| + \left|\tr\left(\nabla_\bSigma \mathcal{N}\left( \z| \bm{\mu}, \bm{\Sigma}\right) (\bSigma-\hat\bSigma)\right)\right| + R_1(\hat\bmu,\hat\bSigma)\\
&\leq \epsilon_\bmu \|\nabla_\bmu \mathcal{N}\left( \z| \bm{\mu}, \bm{\Sigma}\right)\| + D\epsilon_\bSigma \|\nabla_\bSigma \mathcal{N}\left( \z| \bm{\mu}, \bm{\Sigma}\right)\|
+ R_1(\hat\bmu,\hat\bSigma)\\
& =\mathcal{N}\left( \z| \bm{\mu}, \bm{\Sigma}\right)\left(\epsilon_\bmu \|\bm\Sigma^{-1}(\z-\bm\mu)\| + \frac{D\epsilon_\bSigma}{2} \left\|\bSigma^{-1}(\z-\bm\mu)(\z-\bm\mu)^\top\bSigma^{-1}-\bSigma^{-1}\right\| \right)+ R_1(\hat\bmu,\hat\bSigma)\\
&\leq \mathcal{N}\left( \z| \bm{\mu}, \bm{\Sigma}\right)\left(\frac{\epsilon_\bmu}{\alpha} \|\z-\bm\mu\| + \frac{D\epsilon_\bSigma}{2\alpha^2}\left(\left\|\z-\bm\mu\right\|^2+\alpha\right) \right) +  R_1(\hat\bmu,\hat\bSigma)\\
&\leq \mathcal{N}\left( \z| \bm{\mu}, \bm{\Sigma}\right)\left(\frac{\epsilon_\bmu}{\alpha} \left(\|\z\|+\gamma\right) + \frac{D\epsilon_\bSigma}{2\alpha^2}\left((\|\bm z\|+\gamma)^2+\alpha\right) \right) +  R_1(\hat\bmu,\hat\bSigma). 
% &\leq \mathcal{N}\left( \z| \bm{\mu}, \bm{\Sigma}\right)\left(\frac{\epsilon_\bmu}{\alpha} \left(\|\z\|+\gamma\right) + \frac{D\epsilon_\bSigma}{2\alpha^2}\left((\|\bm x\|+\gamma)^2+\alpha\right) \right) +  R_1(\hat\bmu,\hat\bSigma). 
\end{align*}
The third line follows from the assumptions $\|{\bm \mu}-\hat{\bm\mu} \| \leq \epsilon_{\bmu}$ and $\|{\bm\Sigma} - \hat{\bm\Sigma} \| \leq \epsilon_{\bSigma}$, as well as the  
inequality $\tr(\bm A\bm B)\leq D \|\bm A\|\|\bm B\|$ for $\bm A,\bm B\in\RR^{D\times D}$. 
Note that the remainder $R_1(\hat\bmu,\hat\bSigma)$ is bounded by $\mathcal{N}\left( \z| \bm{\mu}, \bm{\Sigma}\right)o(\|\bmu-\hat\bmu\| + \|\bSigma-\hat\bSigma\|)$.
Substituting this bound yields the desired result. \qed
\end{proof}

\begin{proposition}\label{prop:gmmpdfbound}
Let $\M=\sum_{k\in[K]}p^{k}\N(\bmu^{k},\bSigma^{k})$ and $\Mhat=\sum_{k\in[K]}\hat p^k\Nhat(\hat\bmu^k,\hat\bSigma^k)$ be two GM distributions in $\R^D$ with  densities $f$ and $\hat f$, respectively. 
If $|p^{k}-\hat p^k|\leq \epsilon_p$, $\|\bmu^{k}-\hat\bmu^k|\leq\epsilon_\bmu$, and $\|\bSigma^{k}-\hat\bSigma^k\|\leq\epsilon_\bSigma$ for all $k\in[K]$, then the densities satisfy
\begin{align*}
&|f(\z)-\fhat(\z)|\leq f(\z)\left(C_D\|\z\|^2+C'_D\|\z\|+C''_D\right),
\end{align*}
where $C_D=\left(1 + \frac{\epsilon_p}{\underline{p}} \right)  \frac{D \epsilon_\bSigma}{2 \alpha^2}$, $C'_D=\left(1 + \frac{\epsilon_p}{\underline{p}} \right)  \left( \frac{\epsilon_\bmu}{\alpha} + \frac{D \epsilon_\bSigma  \gamma}{\alpha^2} \right)$, and $C''_D=\frac{\epsilon_p}{\underline{p}} 
+ \left(1 + \frac{\epsilon_p}{\underline{p}} \right)
\left( \frac{\epsilon_\bmu}{\alpha} \, \gamma 
+ \frac{D \epsilon_\bSigma}{2 \alpha^2} (\gamma^2 + \alpha) 
+ o(\epsilon_\mu + \epsilon_\bSigma) \right)$.

\end{proposition}
\begin{proof}{Proof}
By the triangle inequality, we have
\begin{align*}
&|f(\z)-\fhat(\z)|\\
=&\left|\sum_{k\in[K]}p^k\mathcal{N}\left( \z| \bm{\mu}^k,\bm{\Sigma}^k\right)-\sum_{k\in[K]}\phat^k\mathcal{N}\left( \z| {\hat\bmu}^k,{\hat\bSigma}^k\right) \right|\\
\leq& \left|\sum_{k\in[K]}p^k\mathcal{N}\left( \z| \bm{\mu}^k,\bm{\Sigma}^k\right)-\sum_{k\in[K]}\hat p^k\mathcal{N}\left( \z| {\bmu}^k,{\bSigma}^k\right) \right|+\left|\sum_{k\in[K]}\hat p^k\mathcal{N}\left( \z| {\bmu}^k,{\bSigma}^k\right)-\sum_{k\in[K]}\phat^k\mathcal{N}\left( \z| {\hat\bmu}^k,{\hat\bSigma}^k\right) \right|\\
\leq& \sum_{k\in[K]}\epsilon_p\mathcal{N}\left( \z| {\bmu}^k,{\bSigma}^k\right) + \sum_{k\in[K]}\hat p^k\left|\mathcal{N}\left( \z| \bm{\mu}^k,\bm{\Sigma}^k\right)-\mathcal{N}\left( \z| {\hat\bmu}^k,{\hat\bSigma}^k\right) \right|.
\end{align*}

Since $\|\bmu^{k}-\hat\bmu^k|\leq\epsilon_\bmu$ and $\|\bSigma^{k}-\hat\bSigma^k\|\leq\epsilon_\bSigma$, Lemma \ref{lem:normalpdfbound} yields:
\begin{align*}
    &\left|\mathcal{N}\left( \z| \bm{\mu}^k, \bm{\Sigma}^k\right)-\mathcal{N}\left( \z| \hat{\bm\mu}^k, \hat{\bm{\Sigma}}^k\right)\right| \\
    \leq & \;\mathcal{N}\left( \z| \bm{\mu}^k, \bm{\Sigma}^k\right)\left(\frac{\epsilon_\bmu}{\alpha} \left(\|\z\|+\gamma\right) + \frac{D\epsilon_\bSigma}{2\alpha^2}\left((\|\bm z\|+\gamma)^2+\alpha\right)  + o(\epsilon_\bmu+\epsilon_\bSigma)\right).
    \end{align*}
We thus obtain
\begin{align*}
&|f(\z)-\fhat(\z)|\\
\leq&\sum_{k\in[K]}\epsilon_p\mathcal{N}\left( \z| {\bmu}^k,{\bSigma}^k\right) +\sum_{k\in[K]}\hat p^k\mathcal{N}\left( \z| \bm{\mu}^k, \bm{\Sigma}^k\right)\left(\frac{\epsilon_\bmu}{\alpha} \left(\|\z\|+\gamma\right) + \frac{D\epsilon_\bSigma}{2\alpha^2}\left((\|\bm z\|+\gamma)^2+\alpha\right)  + o(\epsilon_\bmu+\epsilon_\bSigma)\right)\\
\leq&\sum_{k \in [K]} p^k \, \mathcal{N}(\z | \bmu^k, \bSigma^k)
\Bigg[
\frac{\epsilon_p}{\underline{p}}
+ \left( 1 + \frac{\epsilon_p}{\underline{p}} \right)
\left(
\frac{\epsilon_\bmu}{\alpha} (\|\z\| + \gamma)
+ \frac{D \epsilon_\bSigma}{2 \alpha^2} \left((\|\z\| + \gamma)^2 + \alpha \right)
+ o(\epsilon_\bmu + \epsilon_\bSigma)
\right)\Bigg],
\end{align*}
where we have upper bounded $\hat p^k$ with $p^k+\epsilon_p$ and  $\epsilon_p$ with $\frac{p^k}{\pl}\epsilon_p$, and rearranged the terms. This completes the proof. 
\qed
\end{proof}

\begin{proposition}\label{prop:bounds_conditional_parameters}
Let $\M=\sum_{k\in[K]}p^{k}\N(\bmu^{k},\bSigma^{k})$ and $\Mhat=\sum_{k\in[K]}\hat p^k\Nhat(\hat\bmu^k,\hat\bSigma^k)$ be two GM distributions in $\R^D$.  If $|p^{k}-\hat p^k|\leq \epsilon_p$, $\|\bmu^{k}-\hat\bmu^k|\leq\epsilon_\bmu$, and $\|\bSigma^{k}-\hat\bSigma^k\|\leq\epsilon_\bSigma$ for all $k\in[K]$, then 
\begin{align*}
| p^k_{\xii|\s}-\hat{ p}^k_{\xii|\s} |
&\;\;\leq \left(\frac{p^k \mathcal{N}\left( \bm{s}| {\bm\mu}_s^k, {\bm{\Sigma}}_{\s\s}^k\right) }{f(\s)}
+ \frac{\phat^k\mathcal{N}\left( \bm{s}| \hat{\bm\mu}_s^k, \hat{\bm{\Sigma}}_{ss}^k\right)}{\hat f(\s)}\right) \left(C_{Q}\| \s\|^2+C'_{Q}\| \s\| + C''_{Q}\right)\\
&\;\;\leq2C_{Q}\| \s\|^2+2C'_{Q}\| \s\| +2 C''_{Q},\\
    \|\bm\mu^k_{\xii|\s}-\hat{\bm \mu}^k_{\xii|\s}\| &\leq \left(\frac{\beta}{\alpha}+1\right)\epsilon_\bmu+\frac{\alpha+\beta}{\alpha^2}\left(\|\bm s\|+\gamma \right)\epsilon_\bSigma,\\
    \|\bm\Sigma^k_{\bxi|\s}-\hat{\bm \Sigma}^k_{\bxi|\s}\| &\leq \left(\frac{\beta}{\alpha} \right)^2 \epsilon_\bSigma.
\end{align*}
\end{proposition}
\begin{proof}{Proof of Proposition~\ref{prop:bounds_conditional_parameters}}
Note that $f(\s)=\sum_{k\in[K]} p^k\mathcal{N}\left( \bm{s}| \bm{\mu}_\s^k, \bm{\Sigma}_{\s\s}^k\right)$
and $\fhat(\s)=\sum_{k\in[K]} \hat  p^k\mathcal{N}\left( \bm{s}| \hat{\bm\mu}_\s^k, \hat{\bm{\Sigma}}_{\s\s}^k\right)$. 
By Lemma \ref{lem:ratio_bound},
\begin{align*}
| p^k_{\xii|\s}-\hat{ p}^k_{\xii|\s} | &=\left| \frac{  p^k\mathcal{N}\left( \bm{s}| \bm{\mu}_s^k, \bm{\Sigma}_{ss}^k\right)}{f(\s) } - \frac{ \hat p^k\mathcal{N}\left( \bm{s}| \hat{\bm\mu}_s^k, \hat{\bm{\Sigma}}_{ss}^k\right)}{\fhat(\s)} \right|\\
&\leq \frac{\left|p^k\mathcal{N}\left( \bm{s}| \bm{\mu}_\s^k, \bm{\Sigma}_{\s\s}^k\right)-\hat  p^k\mathcal{N}\left( \bm{s}| \hat{\bm\mu}_\s^k, \hat{\bm{\Sigma}}_{\s\s}^k\right)\right|}{f(\s)}+\frac{\phat^k\mathcal{N}\left( \bm{s}| \hat{\bm\mu}_s^k, \hat{\bm{\Sigma}}_{ss}^k\right)\left|f(\s)-\fhat(\s)\right|}{f(\s)\fhat(\s)}\\
&\leq \underbrace{\frac{\left|p^k\mathcal{N}\left( \bm{s}| \bm{\mu}_\s^k, \bm{\Sigma}_{\s\s}^k\right)-\hat  p^k\mathcal{N}\left( \bm{s}| \hat{\bm\mu}_\s^k, \hat{\bm{\Sigma}}_{\s\s}^k\right)\right|}{f(\s)}}_{\rm (i)}+\underbrace{\frac{|f(\s)-\fhat(\s)|}{f(\s)}}_{\rm (ii)},
\end{align*}
where we use the fact that $\hat p^k\mathcal{N}\left( \bm{s}| \hat{\bm\mu}_\s^k, \hat{\bm{\Sigma}}_{\s\s}^k\right)/\fhat(\s)\leq 1$ in the second line. 
We now bound each term separately:
\begin{enumerate}[(i)]
\item Using the same derivation as in the proof of Proposition \ref{prop:gmmpdfbound}, we find that 
\begin{align*}
\left|p^k\mathcal{N}\left( \bm{s}| \bm{\mu}_\s^k, \bm{\Sigma}_{\s\s}^k\right)-\hat  p^k\mathcal{N}\left( \bm{s}| \hat{\bm\mu}_\s^k, \hat{\bm{\Sigma}}_{\s\s}^k\right)\right|
\leq&p^k \mathcal{N}\left( \bm{s}| {\bm\mu}_s^k, {\bm{\Sigma}}_{\s\s}^k\right) 
\left(C_{Q}\| \s\|^2+C'_{Q}\| \s\| + C''_{Q}\right).
\end{align*}
Hence, 
\begin{align*}
\frac{\left|p^k\mathcal{N}\left( \bm{s}| \bm{\mu}_\s^k, \bm{\Sigma}_{\s\s}^k\right)-\hat  p^k\mathcal{N}\left( \bm{s}| \hat{\bm\mu}_\s^k, \hat{\bm{\Sigma}}_{\s\s}^k\right)\right|}{f(\s)}
\leq &\;\frac{\left|p^k\mathcal{N}\left( \bm{s}| \bm{\mu}_\s^k, \bm{\Sigma}_{\s\s}^k\right)-\hat  p^k\mathcal{N}\left( \bm{s}| \hat{\bm\mu}_\s^k, \hat{\bm{\Sigma}}_{\s\s}^k\right)\right|}{ p^k  \mathcal{N}\left( \bm{s}| {\bm\mu}_s^k, {\bm{\Sigma}}_{\s\s}^k\right) }\\
\leq & \; C_{Q}\| \s\|^2+C'_{Q}\| \s\| + C''_{Q},
\end{align*}
where the first inequality follows from $f(\s)\geq p^k\mathcal{N}\left( \bm{s}| \bm{\mu}_\s^k, \bm{\Sigma}_{\s\s}^k\right)$.
\item By Proposition \ref{prop:gmmpdfbound}, we have 
\begin{align*}
&\frac{\left|f(\s)-\fhat(\s)\right|}{f(\s)}\leq C_{Q}\| \s\|^2+C'_{Q}\| \s\| + C''_{Q}. 
\end{align*}
\end{enumerate}
Thus, combining both bounds, we obtain 
\begin{align*}
| p^k_{\xii|\s}-\hat{ p}^k_{\xii|\s} | &\leq 2C_{Q}\| \s\|^2+2C'_{Q}\| \s\| +2 C''_{Q}. 
\end{align*}

We then derive the error bound for the conditional mean. 
\begin{align*}
    \|\bm\mu^k_{\bxi|\s}-\hat{\bm \mu}^k_{\bxi|\s}\| 
    =& \|\bm{\mu}^k_\bxi + \bm{\Sigma}_{\bxi\s}^k( \bm{\Sigma}_{\s\s}^{k})^{-1}(\s-\bm{\mu}_{\s}^k) - \hat{\bm\mu}^k_\bxi - \Sigmat_{\bxi\s}^k( \bm{\Sigmat}_{\s\s}^{k})^{-1}(\s-\hat{\bm\mu}_s^k)\|\\
    \leq& \|\bm{\mu}^k_\bxi- \hat{\bm\mu}^k_\bxi\|+ \|\bm{\Sigma}_{\bxi\s}^k( \bm{\Sigma}_{\s\s}^{k})^{-1}(\s-\bm{\mu}_{\s}^k)- \Sigmat_{\bxi\s}^k( \bm{\Sigmat}_{\s\s}^{k})^{-1}(\s-\hat{\bm\mu}_s^k)\|\\
    \leq&\epsilon_\mu+\|\bm s\| \cdot \|\bm{\Sigma}_{\bxi\s}^k( \bm{\Sigma}_{\s\s}^{k})^{-1} -  \Sigmat_{\bxi\s}^k( \bm{\Sigmat}_{\s\s}^{k})^{-1}\| + \|\bm{\Sigma}_{\bxi\s}^k( \bm{\Sigma}_{\s\s}^{k})^{-1}\bm\mu_s^k -  \Sigmat_{\bxi\s}^k( \bm{\Sigmat}_{\s\s}^{k})^{-1}\hat{\bm\mu}_s^k\| 
\end{align*}
where the last inequality comes from the fact that the norm of the sub-vector is less than the norm of the whole vector. For the second term, we have
\begin{align*}
    \|\bm{\Sigma}_{\bxi\s}^k( \bm{\Sigma}_{\s\s}^{k})^{-1} -  \Sigmat_{\bxi\s}^k( \bm{\Sigmat}_{\s\s}^{k})^{-1}\|
    \leq & \|\bm{\Sigma}_{\bxi\s}^k( \bm{\Sigma}_{\s\s}^{k})^{-1}- \hat{\bm{\Sigma}}_{\bxi\s}^k( \bm{\Sigma}_{\s\s}^{k})^{-1} \|+\|\hat{\bm{\Sigma}}_{\bxi\s}^k( \bm{\Sigma}_{\s\s}^{k})^{-1} - \Sigmat_{\bxi\s}^k( \bm{\Sigmat}_{\s\s}^{k})^{-1} \|\\
    \leq & \|( \bm{\Sigma}_{\s\s}^{k})^{-1} \|\cdot \|\bm\Sigma_{\bxi\s}^k - \Sigmat_{\bxi\s}^k \|+ \|\Sigmat_{\bxi\s}^k \|\cdot \|( \bm{\Sigma}_{\s\s}^{k})^{-1} - ( \bm{\Sigmat}_{\s\s}^{k})^{-1}\|\\
    \leq & \frac{\epsilon_\Sigma}{\alpha}+\frac{\beta\epsilon_\Sigma}{\alpha^2}\\
    =&\frac{\alpha+\beta}{\alpha^2}\epsilon_\Sigma,
\end{align*}
where the third inequality comes from Lemma~\ref{lem:matrixinv2}. Next, we employ this result to the third term and obtain
\begin{align*}
    \|\bm{\Sigma}_{\bxi\s}^k( \bm{\Sigma}_{\s\s}^{k})^{-1}\bm\mu_s^k -  \Sigmat_{\bxi\s}^k( \bm{\Sigmat}_{\s\s}^{k})^{-1}\hat{\bm\mu}_s^k\|\leq&\|\bm{\Sigma}_{\bxi\s}^k( \bm{\Sigma}_{\s\s}^{k})^{-1}\bm\mu_s^k -  \bm\Sigma_{\bxi\s}^k( \bm{\Sigma}_{\s\s}^{k})^{-1}\hat{\bm\mu}_s^k\|+ \|\bm\Sigma_{\bxi\s}^k( \bm{\Sigma}_{\s\s}^{k})^{-1}\hat{\bm\mu}_s^k- \Sigmat_{\bxi\s}^k( \bm{\Sigmat}_{\s\s}^{k})^{-1}\hat{\bm\mu}_s^k\|\\
    \leq&\|\bm{\Sigma}_{\bxi\s}^k( \bm{\Sigma}_{\s\s}^{k})^{-1} \|\cdot\|\bm \mu_s^k -\hat{\bm\mu}_s^k \|+ \|\bm{\Sigma}_{\bxi\s}^k( \bm{\Sigma}_{\s\s}^{k})^{-1} -  \Sigmat_{\bxi\s}^k( \bm{\Sigmat}_{\s\s}^{k})^{-1} \| \cdot \| \hat{\bm\mu}_s^k\|\\
    \leq& \frac{\beta}{\alpha}\epsilon_\mu+\frac{\alpha+\beta}{\alpha^2}\epsilon_\Sigma\|\hat{\bm\mu}_s^k\|.
\end{align*}
Combining the results for the second and third terms yields
\[
\|\bm\mu^k_{\bxi|\s}-\hat{\bm \mu}^k_{\bxi|\s}\| \leq \left(\frac{\beta}{\alpha}+1\right)\epsilon_\mu+\frac{\alpha+\beta}{\alpha^2}\left(\gamma+\|\bm s\|\right)\epsilon_\Sigma.
\]
Finally, we derive the error bound for the conditional covariance. We first define the $k$-th precision matrix, which is the inverse of the $k$-th covariance matrix, as 
\begin{equation*}
\bm{\Psi}^k = (\bm{\Sigma}^k)^{-1} 
= 
\left[\begin{array}{cc}     
\bm{\Psi}^k_{\s\s} & \bm{\Psi}^k_{\s\bxi} \\
\bm{\Psi}_{\bxi\s}^k & \bm{\Psi}^k_{\bxi\bxi}
\end{array} \right] \in \mathbb{S}^{Q+R}_{++}.
\end{equation*}

Taking $\bm\Psi^k$ and $\hat{\bm\Psi}^k$ into Lemma~\ref{lem:matrixinv2}, we obtain
\[\|\bm\Psi^k- \hat{\bm\Psi}^k\| \leq \frac{\epsilon_\Sigma}{\alpha^2}.
\]
Note that the norm of a submatrix is less than the whole matrix, i.e.,
\[\|\bm\Psi^k_{\bxi\bxi}- \hat{\bm\Psi}^k_{\bxi\bxi}\| \leq \frac{\epsilon_\Sigma}{\alpha^2}.
\]
Since the conditional covariance matrix is equal to the inverse of $\bm\Psi^k_{\bxi\bxi}$, we apply Lemma~\ref{lem:matrixinv2} again and obtain
\[\|\bm\Sigma^k_{\bxi|\s}-\hat{\bm \Sigma}^k_{\bxi|\s}\| \leq \left(\frac{\beta}{\alpha} \right)^2 \epsilon_\Sigma.
\]
This completes the proof. \qed

\end{proof}

\begin{lemma}
\label{lem:norm_xi_bound}
If $\tilde\z\sim\N(\bmu,\bSigma)$ is a Gaussian vector in $\RR^D$ then $\EE_{\N(\bmu,\bSigma)}[\|\bSigma^{-1}(\tilde\z-\bmu)\|]\leq \sqrt{D}\|\bSigma^{-1}\|$. 
\end{lemma}
\begin{proof}{Proof of Lemma~\ref{lem:norm_xi_bound}}
Let $\tilde\z\sim\N(\bmu,\bSigma)$. Then, $\bSigma^{-\frac{1}{2}}(\tilde\z-\bmu)\sim\N(\bm 0,\I)$, and thus $\bSigma^{-1}(\tilde\z-\bmu)\sim\N(\bm 0,\bSigma^{-1})$. We have 
\begin{equation*}
\EE_{\N(\bmu,\bSigma)}[\|\bSigma^{-1}(\tilde\z-\bmu)\|] = \EE_{\N(\bm 0,\bSigma^{-1})}[\|\tilde\z\|]\leq \sqrt{\EE_{\N(\bm 0,\bSigma^{-1})}[\|\tilde\z\|^2]}=\sqrt{\tr(\bSigma^{-1})}\leq \sqrt{D}\|\bSigma^{-1}\|. 
\end{equation*}
where we use Jensen's inequality to upper bound the expectation of a norm. This completes the proof. \qed
\end{proof}

\begin{lemma}
\label{lem:norm_outer_xi_bound}
If $\tilde\z\sim\N(\bmu,\bSigma)$ is a Gaussian vector in $\RR^D$ then 
\begin{align*}
\EE_{\N(\bmu,\bSigma)}\left[\left\|\bSigma^{-1}(\tilde\z-\bm\mu)(\tilde\z-\bm\mu)^\top\bSigma^{-1}-\bSigma^{-1}\right\|\right]\leq D\|\bSigma^{-1}\|.
\end{align*}
\end{lemma}
\begin{proof}{Proof of Lemma~\ref{lem:norm_outer_xi_bound}}
Let $\tilde\z\sim\N(\bmu,\bSigma)$. Then, $\bSigma^{-1}(\tilde\z-\bmu)\sim\N(\bm 0,\bSigma^{-1})$. By \citep[Theorem 4]{koltchinskii2017concentration} we have:
\begin{align*}
\EE_{\N(\bmu,\bSigma)}\left[\left\|\bSigma^{-1}(\tilde\z-\bm\mu)(\tilde\z-\bm\mu)^\top\bSigma^{-1}-\bSigma^{-1}\right\|\right]\leq&\;\; \|\bSigma^{-1}\|\max\left\{\sqrt{\frac{\tr(\bSigma^{-1})}{\|\bSigma^{-1}\|}},\frac{\tr(\bSigma^{-1})}{\|\bSigma^{-1}\|}\right\}\\
=&\;\; \max\left\{\sqrt{\|\bSigma^{-1}\|\tr(\bSigma^{-1})},{\tr(\bSigma^{-1})}\right\}\\
\leq&\;\; \max\left\{\sqrt{D}\|\bSigma^{-1}\|,{D\|\bSigma^{-1}\|}\right\}\\
=&\;\; D\|\bSigma^{-1}\|. 
\end{align*}
Thus, the claim follows. \qed
\end{proof}

\begin{proof}{Proof of Theorem \ref{thm:approx_error}}
Let $f$ and $\fhat$ be the density functions of $\M$ and $\Mhat$, respectively. We have 
\begin{align*}
\left|\EE_{\Ms}\left[\ell(\x,\txi)\right]- \EE_{\Mhats}\left[\ell(\x,\txi)\right]\right|
=&\;\;\left|\int \ell(\x,\xii){f(\xii|\s)}{\rm d}\xii-\int \ell(\x,\xii){\fhat(\xii|\s)}{\rm d}\xii\right|\\
=&\;\;\left|\int \ell(\x,\xii)\left({f(\xii|\s)}-{\fhat(\xii|\s)}\right){\rm d}\xii\right|\\
\leq &\;\;\lbar \int \left|{f(\xii|\s)}-{\fhat(\xii|\s)}\right|{\rm d}\xii. 
\end{align*}
We next bound the difference of conditional densities:
\begin{align*}
&\left|{f(\xii|\s)}-{\fhat(\xii|\s)}\right|\\
= & \left|\sum_{k\in[K]} p_\xis^k\mathcal{N}\left( \bxi| \bm{\mu}_\xis^k, \bm{\Sigma}_{\xis}^k\right)-\sum_{k\in[K]} \hat  p_\xis^k\mathcal{N}\left( \bxi| \hat{\bm\mu}_\xis^k, \hat{\bm{\Sigma}}_{\xis}^k\right)\right|\\
\leq & \left|\sum_{k\in[K]} p_\xis^k\mathcal{N}\left( \bxi| \bm{\mu}_\xis^k, \bm{\Sigma}_{\xis}^k\right)-\sum_{k\in[K]} \hat  p_\xis^k\mathcal{N}\left( \bxi| {\bm\mu}_\xis^k, {\bm{\Sigma}}_{\xis}^k\right)\right|\\
&\qquad+  \left|\sum_{k\in[K]}\hat p_\xis^k\mathcal{N}\left( \bxi | \bm{\mu}_\xis^k, \bm{\Sigma}_{\xis}^k\right)-\sum_{k\in[K]} \hat  p_\xis^k\mathcal{N}\left( \bxi| \hat{\bm\mu}_\xis^k, \hat{\bm{\Sigma}}_{\xis}^k\right)\right|\\
\leq &{\sum_{k\in[K]} \left| p_\xis^k- \hat  p_\xis^k\right|\mathcal{N}\left( \bxi| {\bm\mu}_\xis^k, {\bm{\Sigma}}_{\xis}^k\right)}
+ {\sum_{k\in[K]}\hat p_\xis^k \left|\mathcal{N}\left( \bxi| \bm{\mu}_\xis^k, \bm{\Sigma}_{\xis}^k\right)-\mathcal{N}\left( \bxi| \hat{\bm\mu}_\xis^k, \hat{\bm{\Sigma}}_{\xis}^k\right)\right|}.
\end{align*}
Hence, 
\begin{align*}
&\left|\EE_{\Ms}\left[\ell(\x,\txi)\right]- \EE_{\Mhats}\left[\ell(\x,\txi)\right]\right|\\
\leq&\lbar \underbrace{\int \sum_{k\in[K]} \left| p_\xis^k- \hat  p_\xis^k\right|\mathcal{N}\left( \bxi| {\bm\mu}_\xis^k, {\bm{\Sigma}}_{\xis}^k\right){\rm d}\xii}_{\rm (i)}+\lbar
\underbrace{\int  \sum_{k\in[K]}\hat p_\xis^k \left|\mathcal{N}\left( \bxi| \bm{\mu}_\xis^k, \bm{\Sigma}_{\xis}^k\right)-\mathcal{N}\left( \bxi| \hat{\bm\mu}_\xis^k, \hat{\bm{\Sigma}}_{\xis}^k\right)\right| {\rm d}\xii}_{\rm (ii)}. 
\end{align*}

We analyze each integral separately:
\begin{enumerate}[(i)]
\item Using Proposition \ref{prop:bounds_conditional_parameters}, we obtain
\begin{align*}
 &\;\;\int\sum_{k\in[K]} \left| p_\xis^k- \hat  p_\xis^k\right|\mathcal{N}\left( \bxi| {\bm\mu}_\xis^k, {\bm{\Sigma}}_{\xis}^k\right){\rm d}\xii   \\
 \leq&\;\;\left(C_{Q}\| \s\|^2+C'_{Q}\| \s\| + C''_{Q}\right) \sum_{k\in[K]}\left(\frac{p^k \mathcal{N}\left( \bm{s}| {\bm\mu}_s^k, {\bm{\Sigma}}_{\s\s}^k\right) }{f(\s)}
+ \frac{\phat^k\mathcal{N}\left( \bm{s}| \hat{\bm\mu}_s^k, \hat{\bm{\Sigma}}_{ss}^k\right)}{\hat f(\s)}\right) \int \mathcal{N}\left( \bxi| {\bm\mu}_\xis^k, {\bm{\Sigma}}_{\xis}^k\right){\rm d}\xii\\
=&\;\;2\left(C_{Q}\| \s\|^2+C'_{Q}\| \s\| + C''_{Q}\right). 
\end{align*}
\item Using Lemma \ref{lem:normalpdfbound}, we have
\begin{align*}
&\;\; \left|\mathcal{N}\left( \bxi| \bm{\mu}_\xis^k, \bm{\Sigma}_{\xis}^k\right)-\mathcal{N}\left( \bxi| \hat{\bm\mu}_\xis^k, \hat{\bm{\Sigma}}_{\xis}^k\right)\right|\\
\leq &\;\; \mathcal{N}\left( \bxi| \bm{\mu}_\xis^k, \bm{\Sigma}_{\xis}^k\right)\Bigg(\epsilon'_\bmu \|(\bSigma^{k}_\xis)^{-1}(\bxi-\bm\mu^k_\xis)\| \\
&\qquad+ \frac{R\epsilon'_\bSigma}{2} \left\|(\bSigma^k_\xis)^{-1}(\bxi-\bm\mu^k_\xis)(\bxi-\bm\mu^k_\xis)^\top(\bSigma^k_\xis)^{-1}-(\bSigma^k_\xis)^{-1}\right\| +o(\epsilon_\bmu + \epsilon_\bSigma) \Bigg),
\end{align*}
where $\epsilon'_\bmu=\left(\frac{\beta}{\alpha}+1\right)\epsilon_\bmu+\frac{\alpha+\beta}{\alpha^2}\left(\|\bm s\|+\gamma \right)\epsilon_\bSigma$ and $\epsilon'_\bSigma=\left(\frac{\beta}{\alpha} \right)^2 \epsilon_\bSigma$ by Proposition \ref{prop:bounds_conditional_parameters}. 
Hence, 
\begin{align*}
&\;\;\int  \sum_{k\in[K]}\hat p_\xis^k \left|\mathcal{N}\left( \bxi| \bm{\mu}_\xis^k, \bm{\Sigma}_{\xis}^k\right)-\mathcal{N}\left( \bxi| \hat{\bm\mu}_\xis^k, \hat{\bm{\Sigma}}_{\xis}^k\right)\right| {\rm d}\xii\\
\leq &\;\; \sum_{k\in[K]}\hat p_\xis^k \EE_{\N(\bm{\mu}_{\xii|\s}^k,\bSigma_{\xii|\s}^k)}\left[\epsilon'_\bmu \|(\bSigma^{k}_\xis)^{-1}(\txi-\bm\mu^k_\xis)\|\right] + o(\epsilon_\bmu + \epsilon_\bSigma)\\
&\;\;\quad + \sum_{k\in[K]}\hat p_\xis^k \EE_{\N(\bm{\mu}_{\xii|\s}^k,\bSigma_{\xii|\s}^k)}\left[ \frac{R\epsilon'_\bSigma}{2} \left\|(\bSigma^k_\xis)^{-1}(\txi-\bm\mu^k_\xis)(\txi-\bm\mu^k_\xis)^\top(\bSigma^k_\xis)^{-1}-(\bSigma^k_\xis)^{-1}\right\|\right] . 
\end{align*}
By Lemma \ref{lem:norm_xi_bound}, we have
\begin{equation*}
\EE_{\N(\bm{\mu}_{\xii|\s}^k,\bSigma_{\xii|\s}^k)}\left[ \|(\bSigma^{k}_\xis)^{-1}(\txi-\bm\mu^k_\xis)\|\right] \leq \sqrt{R}\left\|(\bSigma_{\xii|\s}^k)^{-1}\right\|,
\end{equation*}
and by Lemma \ref{lem:norm_outer_xi_bound}, we have 
\begin{align*}
\EE_{\N(\bm{\mu}_{\xii|\s}^k,\bSigma_{\xii|\s}^k)}\left[  \left\|(\bSigma^k_\xis)^{-1}(\txi-\bm\mu^k_\xis)(\txi-\bm\mu^k_\xis)^\top(\bSigma^k_\xis)^{-1}-(\bSigma^k_\xis)^{-1}\right\|\right]\leq R\left\|(\bSigma_{\xii|\s}^k)^{-1}\right\|. 
\end{align*}
In addition, since $\bm{\Sigma}^k_{\xii|\s} \displaystyle = \;\; \bm{\Sigma}_{\bxi\bxi}^k- \bm{\Sigma}_{\xii\s}^k( \bm{\Sigma}_{\s\s}^{k})^{-1} \bm{\Sigma}_{\s\xii}^k$ constitutes a Schur complement of $\bSigma^k$, we must have $\lambda_{\text{max}}\left((\bm{\Sigma}^k_{\xii|\s})^{-1}\right)\leq \lambda_{\text{max}}((\bSigma^k)^{-1})\leq \frac{1}{\alpha}$   \cite[Lemma 2.3]{zhang2006schur}. 
This yields:
\begin{align*}
&\;\;\int  \sum_{k\in[K]}\hat p_\xis^k \left|\mathcal{N}\left( \bxi| \bm{\mu}_\xis^k, \bm{\Sigma}_{\xis}^k\right)-\mathcal{N}\left( \bxi| \hat{\bm\mu}_\xis^k, \hat{\bm{\Sigma}}_{\xis}^k\right)\right| {\rm d}\xii\\
\leq &\;\; \sum_{k\in[K]}\hat p_\xis^k \left(\sqrt{R}\epsilon'_\bmu+\frac{1}{2}R^2\epsilon'_\bSigma \right)\left\|(\bSigma_{\xii|\s}^k)^{-1}\right\| + o(\epsilon_\bmu + \epsilon_\bSigma)\\
\leq &\;\; \sum_{k\in[K]}\hat p_\xis^k \left(\sqrt{R}\epsilon'_\bmu+\frac{1}{2}R^2\epsilon'_\bSigma\right)\frac{1}{\alpha} + o(\epsilon_\bmu + \epsilon_\bSigma)\\
=&\;\;\left(\sqrt{R}\epsilon'_\bmu+\frac{1}{2}R^2\epsilon'_\bSigma\right)\frac{1}{\alpha} + o(\epsilon_\bmu + \epsilon_\bSigma). 
\end{align*}
\end{enumerate}

Combining both bounds, we obtain
\begin{align*}
&\left|\EE_{\Ms}\left[\ell(\x,\txi)\right]- \EE_{\Mhats}\left[\ell(\x,\txi)\right]\right|\\
\leq&\lbar \left(2\left(C_{Q}\| \s\|^2+C'_{Q}\| \s\| + C''_{Q}\right)+\left(\sqrt{R}\epsilon'_\bmu+\frac{1}{2}R^2\epsilon'_\bSigma\right)\frac{1}{\alpha}\right),
\end{align*}
which proves the claim.  
\qed
\end{proof}

\begin{proof}{Proof of Proposition~\ref{prop:normalizing_flow}.}
Let \(\Q^\btheta\) denote the pushforward of the latent distribution \(\M\) under the separable normalizing flow
\[
T_\btheta(\s,\bxi)=\big(H_\btheta(\s),F_\btheta(\bxi;\s)\big).
\]
Fix \((\s',\bxi')\), and define
\[
\s=H_\btheta^{-1}(\s'),
\qquad
\bxi=F_\btheta^{-1}(\bxi';\s).
\]
Then
\[
T_\btheta^{-1}(\s',\bxi')
=
\left(
H_\btheta^{-1}(\s'),
F_\btheta^{-1}\big(\bxi';H_\btheta^{-1}(\s')\big)
\right)
=
(\s,\bxi).
\]

By the change-of-variables formula, the joint density of \((\ts',\txi')\) under \(\Q^\btheta\) satisfies
\[
f_{\Q^\btheta}(\s',\bxi')
=
f_{\M}(\s,\bxi)
\left|
\det \mathbf J_{T_\btheta^{-1}}(\s',\bxi')
\right|.
\]
Similarly, since \(\s'=H_\btheta(\s)\), the marginal density of \(\ts'\) is
\[
f_{\Q^\btheta_{\s'}}(\s')
=
f_{\M_{\s}}(\s)
\left|
\det \mathbf J_{H_\btheta^{-1}}(\s')
\right|.
\]
Therefore,
\[
f_{\Q^\btheta_{\bxi'|\s'}}(\bxi'|\s')
=
\frac{
f_{\Q^\btheta}(\s',\bxi')
}{
f_{\Q^\btheta_{\s'}}(\s')
}
=
\frac{
f_{\M}(\s,\bxi)
}{
f_{\M_{\s}}(\s)
}
\cdot
\frac{
\left|\det \mathbf J_{T_\btheta^{-1}}(\s',\bxi')\right|
}{
\left|\det \mathbf J_{H_\btheta^{-1}}(\s')\right|
}.
\]
The first ratio is the latent conditional density:
\[
\frac{
f_{\M}(\s,\bxi)
}{
f_{\M_{\s}}(\s)
}
=
f_{\M_{\bxi|\s}}(\bxi|\s).
\]

It remains to simplify the Jacobian determinant. Since the inverse map has the form
\[
T_\btheta^{-1}(\s',\bxi')
=
\left(
H_\btheta^{-1}(\s'),
F_\btheta^{-1}\big(\bxi';H_\btheta^{-1}(\s')\big)
\right),
\]
its Jacobian is block lower triangular:
\[
\mathbf J_{T_\btheta^{-1}}(\s',\bxi')
=
\begin{bmatrix}
\dfrac{\partial H_\btheta^{-1}(\s')}{\partial \s'} & \bm 0 \\[1.2em]
\dfrac{\partial F_\btheta^{-1}(\bxi';\s)}{\partial \s'} &
\dfrac{\partial F_\btheta^{-1}(\bxi';\s)}{\partial \bxi'}
\end{bmatrix}.
\]
Hence,
\[
\left|
\det \mathbf J_{T_\btheta^{-1}}(\s',\bxi')
\right|
=
\left|
\det \mathbf J_{H_\btheta^{-1}}(\s')
\right|
\cdot
\left|
\det \mathbf J_{F_\btheta^{-1}}(\bxi';\s)
\right|,
\]
where the second Jacobian is taken with respect to \(\bxi'\), holding \(\s\) fixed. Substituting this into the previous display yields
\[
f_{\Q^\btheta_{\bxi'|\s'}}(\bxi'|\s')
=
f_{\M_{\bxi|\s}}(\bxi|\s)
\left|
\det \mathbf J_{F_\btheta^{-1}}(\bxi';\s)
\right|,
\]
where
\[
\s=H_\btheta^{-1}(\s'),
\qquad
\bxi=F_\btheta^{-1}(\bxi';\s).
\]
This proves the claim. \qed
\end{proof}

\section{Proofs of Section~\ref{sec:DRO}}

\begin{proof}{Proof of Proposition~\ref{prop:dro_reformulation}}
Under the prescribed piecewise affine loss function, the formulation \eqref{eq:cond_DRO_min_SAA} becomes 
\begin{align*}
\min_{\x\in\X,\lambda\in\RR_+} \varepsilon^2\lambda+\frac{1}{M}\sum_{m\in[M]}\sup_{\bm\omega\in\RR^R}\max_{j\in[J]}\bm a_j(\x)^\top\bm\omega+b_j(\x)-\lambda \|\bm\omega-\xii_m\|^2. 
\end{align*}
Introducing epigraphical variables $\bm\gamma\in\RR^M$ to bring the summands into the constraint system, we get:
\begin{equation}
\label{eq:semiinfinite_program}
\begin{array}{rl}
\min_{}\;\;\; &\displaystyle\varepsilon^2\lambda+\frac{1}{M}\sum_{m\in[M]}\gamma_m\\
\st\;\;& \displaystyle\x\in\X,\;\lambda\in\RR_+,\;\bm \gamma\in\RR^M\\
&\displaystyle\bm a_j(\x)^\top\bm\omega+b_j(\x)-\lambda \|\bm\omega-\xii_m\|^2\leq\gamma_m\quad\forall \bm\omega\in\RR^R\quad\forall j\in[J] \quad\forall m\in[M] 
\end{array}
\end{equation}
Each semi-infinite constraint is equivalent to a semidefinite constraint:
\begin{align*}
&\;\bm a_j(\x)^\top\xii+b_j(\x)-\lambda \|\bm\omega-\xii_m\|^2\leq\gamma_m\quad\forall \bm\omega\in\RR^R\\
\Longleftrightarrow\;&\;0 \leq -\bm a_j(\x)^\top\bm\omega-b_j(\x)+\lambda \|\bm\omega-\xii_m\|^2 + \gamma_m\quad\forall \bm\omega\in\RR^R\\
\Longleftrightarrow\;&\;0 \leq -\bm a_j(\x)^\top\bm\omega-b_j(\x)+\lambda\bm\omega^\top\bm\omega-2\lambda\bxi_m^\top\bm\omega+ \lambda\xii_m^\top\xii_m + \gamma_m\quad\forall \bm\omega\in\RR^R\\
\Longleftrightarrow&\bm 0\preceq \begin{bmatrix}
\lambda\I &\; -\left(\lambda\bxi_m + \frac{\bm a_j(\x)}{2}\right) \\ 
-\left(\lambda\bxi_m + \frac{\bm a_j(\x)}{2}\right)^\top & \lambda\xii_m^\top\xii_m + \gamma_m-b_j(\x)
\end{bmatrix}.
\end{align*}
By the Schur complement, this constraint is equivalent to the hyperbolic constraint:
\begin{align*}
 &\;\;\lambda\xii_m^\top\xii_m + \gamma_m-b_j(\x) \geq \frac{\left(\lambda\bxi_m + \frac{\bm a_j(\x)}{2}\right)^\top\left(\lambda\bxi_m + \frac{\bm a_j(\x)}{2}\right)}{\lambda},\quad \lambda\xii_m^\top\xii_m + \gamma_m-b_j(\x)\geq 0\\
\Longleftrightarrow &\;\; \left\|\begin{bmatrix} 2\lambda\bxi_m + \bm a_j(\x)  \\
\lambda\xii_m^\top\xii_m + \gamma_m-b_j(\x)-\lambda \end{bmatrix}\right\| \leq \lambda\xii_m^\top\xii_m + \gamma_m-b_j(\x) + \lambda, \quad \lambda\xii_m^\top\xii_m + \gamma_m-b_j(\x)\geq 0.
\end{align*}
Substituting this constraint back into \eqref{eq:semiinfinite_program} yields the desired formulation. \qed
\end{proof}

We then present the proof of Theorem \ref{thm:coverage_radius}. To this end, we first derive several auxiliary results. 
\begin{proposition}
\label{prop:Wasserstein_bound_Gaussian}
 The 2-Wasserstein distance between two Gaussians $\N(\bmu,\bSigma)$ and $\Nhat(\hat\bmu,\hat\bSigma)$  in $\RR^D$ admits the following upper bound:
\begin{equation*}
\Wass_2^2(\N(\bmu,\bSigma),\Nhat(\hat\bmu,\hat\bSigma))\leq \|\bmu-\hat\bmu\|^2 + D  \|\bSigma-\hat\bSigma\|.%\leq  \eps_\bmu^2 + D\epsilon_\bSigma. 
\end{equation*}
\end{proposition}
\begin{proof}{Proof of Proposition~\ref{prop:Wasserstein_bound_Gaussian}}
The 2-Wasserstein distance between two Gaussians $\N$ and $\Nhat$ \cite[Proposition 7]{givens1984class} has the closed form formula given by  
\begin{equation*}
\Wass_2^2(\N(\bmu,\bSigma),\Nhat(\hat\bmu,\hat\bSigma))=\|\bmu-\hat\bmu\|^2+\tr\left(\bSigma+\hat\bSigma-2((\bSigma^{\frac{1}{2}}\hat\bSigma\bSigma^{\frac{1}{2}})^{\frac{1}{2}}\right).
\end{equation*}
%By our assumption, we have $\|\bmu-\hat\bmu\|^2\leq \epsilon_\bmu^2$. 
It remains to upper bound the trace term. From \cite[Theorem 1]{bhatia2019bures}, we have  $\tr\left(\bSigma+\hat\bSigma-2((\bSigma^{\frac{1}{2}}\hat\bSigma\bSigma^{\frac{1}{2}})^{\frac{1}{2}}\right)^{\frac{1}{2}}=\min_{\bm U\in \mathcal U(d)}\|\bSigma^{\frac{1}{2}}-\hat\bSigma^{\frac{1}{2}}\bm U\|_F$, where $\mathcal U(D)$ is the set of all $D\times D$ unitary matrices. Setting $\bm U=\I$, we get 
\begin{equation*}
\tr\left(\bSigma+\hat\bSigma-2((\bSigma^{\frac{1}{2}}\hat\bSigma\bSigma^{\frac{1}{2}})^{\frac{1}{2}}\right)^2\leq \|\bSigma^{\frac{1}{2}}-\hat\bSigma^{\frac{1}{2}}\|^2_F\leq{D}\|\bSigma^{\frac{1}{2}}-\hat\bSigma^{\frac{1}{2}}\|^2.
\end{equation*}
Since $\|\bSigma^{\frac{1}{2}}-\hat\bSigma^{\frac{1}{2}}\|^2 \leq \|\bSigma-\hat\bSigma\|$ \citep{phillips1987uniform}, we thus obtain 
\begin{equation*}
\tr\left(\bSigma+\hat\bSigma-2((\bSigma^{\frac{1}{2}}\hat\bSigma\bSigma^{\frac{1}{2}})^{\frac{1}{2}}\right)^2\leq D \|\bSigma-\hat\bSigma\|. 
\end{equation*}
Combining the bounds yields the desired result. \qed
\end{proof}

\begin{proposition}
\label{prop:Wasserstein_GMMs}
Let $\M=\sum_{k\in[K]}p^{k}\N(\bmu^{k},\bSigma^{k})$ and $\Mhat=\sum_{k\in[K]}\hat p^k\Nhat(\hat\bmu^k,\hat\bSigma^k)$ be two GM distributions in $\R^D$.  If $\|\bmu^{k}-\hat\bmu^k|\leq\epsilon_\bmu$ and $\|\bSigma^{k}-\hat\bSigma^k\|\leq\epsilon_\bSigma$ for all $k\in[K]$, then 
\begin{align*}
\Wass_2^2(\M,\Mhat)\leq \epsilon_\bmu^2 + D\epsilon_\bSigma + \left(4\gamma^2 +2D\beta\right)\sum_{i\in[K]}{|p^{i}-\hat p^i|}.
\end{align*}
\end{proposition}
\begin{proof}{Proof of Proposition~\ref{prop:Wasserstein_GMMs}}
From \cite[Section 4]{delon2020wasserstein} and \cite[Section 3]{chen2018optimal}, we find that for any two GM distributions $\M$ and $\Mhat$ their 2-Wasserstein distance is upper bounded by:
\begin{equation}
    \label{eq:bound_Wasserstein_GM}
    \begin{array}{rl}
 \displaystyle \Wass_2^2(\M,\Mhat)\leq \min& \displaystyle \sum_{i,j\in[K]}\pi_{ij} \Wass_2^2(\N(\bmu^{i},\bSigma^{i}),\Nhat(\hat\bmu^{j},\hat\bSigma^{j}))\\
 \displaystyle \st& \displaystyle\bm\pi\in\Delta^K\times\Delta^K\\
  & \displaystyle\sum_{i\in[K]}\pi_{ij}=\hat p^j\quad \forall j\in[K]\\
  &  \displaystyle \sum_{j\in[K]}\pi_{ij}= p^i\quad \forall i\in[K].
  \end{array}
\end{equation}
Suppose that, for every $k\in[K]$, $p^k=\hat p^k + \delta_k$ for some $\delta_k\in[-1,1]$. By construction, we must have $\sum_{k\in[K]}\delta_k=0$ in order for both $\bm p$ and $\hat{\bm p}$ to be probability mass functions. Define
\begin{align*}
q^k:=\min(p^k,\hat p^k), \qquad
r^k:=(p^k-\hat p^k)_+, \qquad
s^k:=(\hat p^k-p^k)_+,
\end{align*}
so that
\begin{align*}
p^k=q^k+r^k, \qquad \hat p^k=q^k+s^k,
\end{align*}
and
\begin{align*}
m:=\sum_{k\in[K]}r^k=\sum_{k\in[K]}s^k=\frac{1}{2}\sum_{k\in[K]}|p^k-\hat p^k|
=\frac{1}{2}\sum_{k\in[K]}|\delta_k|.
\end{align*}
If $m=0$, then $p^k=\hat p^k$ for all $k\in[K]$, and choosing $\pi_{kk}=p^k$ and $\pi_{ij}=0$ for $i\neq j$ immediately yields the claim. Suppose now that $m>0$. Consider the solution
\begin{align*}
\pi_{ij}=\begin{cases}\displaystyle 
q^i+\frac{r^is^j}{m} &\text{ if }\; i=j\\
\displaystyle \frac{r^is^j}{m}  &\text{ if }\; i\neq j.
\end{cases}
\end{align*}
Since $\pi_{ij}\geq 0$ for all $i,j\in[K]$, and
\begin{align*}
\sum_{j\in[K]}\pi_{ij}
= q^i + \frac{r^i}{m}\sum_{j\in[K]} s^j
= q^i+r^i
= p^i,
\end{align*}
and similarly
\begin{align*}
\sum_{i\in[K]}\pi_{ij}
= q^j + \frac{s^j}{m}\sum_{i\in[K]} r^i
= q^j+s^j
= \hat p^j,
\end{align*}
 this solution is feasible to \eqref{eq:bound_Wasserstein_GM}, which implies that
\begin{align*}
\Wass_2^2(\M,\Mhat)\leq&\; \sum_{k\in[K]} q^k \Wass_2^2(\N(\bmu^{k},\bSigma^{k}),\Nhat(\hat\bmu^{k},\hat\bSigma^{k})) \\
&\; + \frac{1}{m}\sum_{i,j\in[K]} r^i s^j \Wass_2^2(\N(\bmu^{i},\bSigma^{i}),\Nhat(\hat\bmu^{j},\hat\bSigma^{j})).
\end{align*}
Using Proposition \ref{prop:Wasserstein_bound_Gaussian} and by our assumptions, we thus arrive at
\begin{align*}
\Wass_2^2(\M,\Mhat)\leq& \sum_{k\in[K]} q^k  \left(\|\bmu^k-\hat\bmu^k\|^2 +D\|\bSigma^k-\hat\bSigma^k\|\right) \\
&\; + \frac{1}{m}\sum_{i,j\in[K]} r^i s^j \left(\|\bmu^i-\hat\bmu^j\|^2 +D\|\bSigma^i-\hat\bSigma^j\|\right)\\
\leq &\sum_{k\in[K]} q^k \left(\epsilon_\bmu^2 + D\epsilon_\bSigma \right)
+ \frac{1}{m}\sum_{i,j\in[K]} r^i s^j \left(4\gamma^2 +2D\beta\right)\\
\leq& \epsilon_\bmu^2 + D\epsilon_\bSigma
+ \frac{4\gamma^2 +2D\beta}{m}\left(\sum_{i\in[K]} r^i\right)\left(\sum_{j\in[K]} s^j\right)\\
=& \epsilon_\bmu^2 + D\epsilon_\bSigma
+ \left(4\gamma^2 +2D\beta\right)m\\
=& \epsilon_\bmu^2 + D\epsilon_\bSigma
+ \frac{1}{2}\left(4\gamma^2 +2D\beta\right)\sum_{i\in[K]}|\delta_i|\\
\leq& \epsilon_\bmu^2 + D\epsilon_\bSigma
+ \left(4\gamma^2 +2D\beta\right)\sum_{i\in[K]}|\delta_i|.
\end{align*}
This completes the proof. \qed
\end{proof}
Equipped with these preliminary results, we now present the proof of Theorem \ref{thm:coverage_radius}.

\begin{proof}{Proof of Theorem \ref{thm:coverage_radius}}
By Proposition \ref{prop:bounds_conditional_parameters}, we have 
\begin{align*}
| p^k_{\xii|\s}-\hat{ p}^k_{\xii|\s} | &\leq\left(\frac{p^k \mathcal{N}\left( \bm{s}| {\bm\mu}_s^k, {\bm{\Sigma}}_{\s\s}^k\right) }{f(\s)}
+ \frac{\phat^k\mathcal{N}\left( \bm{s}| \hat{\bm\mu}_s^k, \hat{\bm{\Sigma}}_{ss}^k\right)}{\hat f(\s)}\right) \left(C_{Q}\| \s\|^2+C'_{Q}\| \s\| + C''_{Q}\right)\\
    \|\bm\mu^k_{\xii|\s}-\hat{\bm \mu}^k_{\xii|\s}\| &\leq \left(\frac{\beta}{\alpha}+1\right)\epsilon_\bmu+\frac{\alpha+\beta}{\alpha^2}\left(\|\bm s\|+\gamma \right)\epsilon_\bSigma\\
    \|\bm\Sigma^k_{\bxi|\s}-\hat{\bm \Sigma}^k_{\bxi|\s}\| &\leq \left(\frac{\beta}{\alpha} \right)^2 \epsilon_\bSigma.
\end{align*}
% \begin{proof}{Proof of Theorem \ref{thm:coverage_radius}}
% By Proposition \ref{prop:bounds_conditional_parameters}, we have 
% \begin{align*}
% | p^k_{\xii|\s}-\hat{ p}^k_{\xii|\s} | &\leq\left(\frac{p^k \mathcal{N}\left( \bm{s}| {\bm\mu}_s^k, {\bm{\Sigma}}_{\s\s}^k\right) }{f(\s)}
% + \frac{\phat^k\mathcal{N}\left( \bm{s}| \hat{\bm\mu}_s^k, \hat{\bm{\Sigma}}_{ss}^k\right)}{\hat f(\s)}\right) \left(C_{Q}\| \s\|^2+C'_{Q}\| \s\| + C''_{Q}\right)\\
%     \|\bm\mu^k_{\xii|\s}-\hat{\bm \mu}^k_{\xii|\s}\| &\leq \left(\frac{\beta}{\alpha}+1\right)\epsilon_\bmu+\frac{\alpha+\beta}{\alpha^2}\left(\|\bm s\|+\gamma \right)\epsilon_\bSigma\\
%     \|\bm\Sigma^k_{r|s}-\hat{\bm \Sigma}^k_{r|s}\| &\leq \left(\frac{\beta}{\alpha} \right)^2 \epsilon_\bSigma.
% \end{align*}
Hence, Proposition \ref{prop:Wasserstein_GMMs} implies that 
\begin{align*}
&\;\Wass_2^2(\Ms,\Mhats)\\
\leq&\;\left(\left(\frac{\beta}{\alpha}+1\right)\epsilon_\bmu+\frac{\alpha+\beta}{\alpha^2}\left(\|\bm s\|+\gamma \right)\epsilon_\bSigma\right)^2 + R\left(\frac{\beta}{\alpha} \right)^2 \epsilon_\bSigma\\
&\;\quad+ \left(4\gamma^2 +2R\beta\right)\sum_{i\in[K]} \left(\frac{p^i \mathcal{N}\left( \bm{s}| {\bm\mu}_s^i, {\bm{\Sigma}}_{\s\s}^i\right) }{f(\s)}
+ \frac{\phat^i\mathcal{N}\left( \bm{s}| \hat{\bm\mu}_s^i, \hat{\bm{\Sigma}}_{ss}^i\right)}{\hat f(\s)}\right) \left(C_{Q}\| \s\|^2+C'_{Q}\| \s\| + C''_{Q}\right)\\
\leq&\;\left(\left(\frac{\beta}{\alpha}+1\right)\epsilon_\bmu+\frac{\alpha+\beta}{\alpha^2}\left(\|\bm s\|+\gamma \right)\epsilon_\bSigma\right)^2 + R\left(\frac{\beta}{\alpha} \right)^2 \epsilon_\bSigma\\
&\;\quad+ 2 \left(4\gamma^2 +2R\beta\right)\left(C_{Q}\| \s\|^2+C'_{Q}\| \s\| + C''_{Q}\right).
\end{align*}
This completes the proof. 
\qed
\end{proof}

\begin{proof}{Proof of Corollary~\ref{coro:performance_bound}}
Since $\Ms\in\mathcal P_\varepsilon$, we must have 
\begin{equation*}
\EE_{\Ms}[\ell(\x,\txi)]\leq \sup_{\Q\in\mP_\varepsilon}\EE_{\Q}\left[\ell(\x,\txi)\right]\qquad\forall \x\in\X.
\end{equation*}
In particular, substituting the solution $\hat\x $ into both sides of the inequalities yields the desired guarantee. 
\qed
\end{proof}

\begin{proof}{Proof of Proposition~\ref{prop:dro_k}}
By the triangle inequality for the Wasserstein distance, we obtain
\begin{align*}
\Wass_2(\Ms^K,\Mhats^{K'})&\leq \Wass_2(\Ms^K,\Mhats^K)+ \Wass_2(\Mhats^K,\Mhats^{K'}).
\end{align*}
From Theorem \ref{thm:coverage_radius}, we have $ \Wass_2(\Ms^K,\Mhats^K)\leq \varepsilon$. Furthermore, from \cite[Section 4]{delon2020wasserstein} and \cite[Section 3]{chen2018optimal}, we obtain that $\Wass_2(\Mhats^K,\Mhats^{K'})$ is upper bounded by the optimal value of the following linear program: 
\begin{equation*}
    \begin{array}{rl}
 \displaystyle \Wass^2_2(\Mhats^K,\Mhats^{K'})\leq\overline{\Wass}^2_2(\Mhats^K,\Mhats^{K'}) \coloneqq\min & \displaystyle \sum_{i\in[K]}\sum_{j\in[K']}\pi_{ij} \Wass_2^2(\Nhat(\hat\bmu_{\bxi|\s}^{i},\hat\bSigma_{\bxi|\s}^{i}),\Nhat(\hat\bmu_{\bxi|\s}^{j},\hat\bSigma_{\bxi|\s}^{j}))\\
 \displaystyle \st& \displaystyle\bm\pi\in\Delta^K\times\Delta^{K'}\\
  & \displaystyle\sum_{i\in[K]}\pi_{ij}=\hat p^j\quad \forall j\in[K']\\
  &  \displaystyle \sum_{j\in[K']}\pi_{ij}=\hat p^i\quad \forall i\in[K],
  \end{array}
\end{equation*}
where $\Wass_2^2(\Nhat(\bmu_{\bxi|\s}^{i},\bSigma_{\bxi|\s}^{i}),\Nhat(\hat\bmu_{\bxi|\s}^{j},\hat\bSigma_{\bxi|\s}^{j}))=|\hat\bmu^i_{\bxi|\s}-\hat\bmu^j_{\bxi|\s}\|^2+\tr\left(\hat\bSigma_{\bxi|\s}^i+\hat\bSigma_{\bxi|\s}^j-2(((\hat\bSigma_{\bxi|\s}^i)^{\frac{1}{2}}\hat\bSigma_{\bxi|\s}^j(\hat\bSigma_{\bxi|\s}^i)^{\frac{1}{2}})^{\frac{1}{2}}\right)$  \cite[Proposition 7]{givens1984class}. 
Since $K\in\K$, we must have $\Wass_2(\Ms^K,\Mhats^{K'})\leq \varepsilon +\max_{L\in\K}{\overline{\Wass}_2(\Mhats^L,\Mhats^{K'})}$. The result then follows from optimizing for the best mixture size $K'\in\mathcal K$ that minimizes the bound. \qed
\end{proof}

\section{Proofs of Section~\ref{sec:observational}}

\begin{proof}{Proof of Theorem \ref{thm:error_observe_fix_x}}
We decompose the error as follows:
\begin{align*}
&\displaystyle\;\;\;\;\left|\EE_{\M}\left[\ell(\x,\txi)\big|\ts=\s\right]-\frac{\frac{1}{N}\sum_{n\in[N]}\ell(\x,\bxi_n)\hat f(\s|\bxi_n) }{\fhat(\s) }\right|\\
&\displaystyle=\left|\frac{\int \ell(\x,\xii) f(\s,\xii) {\rm d}\xii}{f(\s)}- \frac{\frac{1}{N}\sum_{n\in[N]} \ell(\x,\bxi_n)\hat f(\s|\bxi_n)}{\fhat(\s)}\right|\\
&\displaystyle\leq\underbrace{\left|\frac{\int \ell(\x,\xii) f(\s,\xii) {\rm d}\xii}{ f(\s)}-\frac{\int \ell(\x,\xii) \hat f(\s|\xii) f(\xii){\rm d}\xii}{ \hat f(\s)}\right|}_{\rm (a)}+\underbrace{\left|\frac{\int \ell(\x,\xii) \hat f(\s|\xii) f(\xii){\rm d}\xii}{\fhat(\s)}- \frac{\frac{1}{N}\sum_{n\in[N]}\ell(\x,\bxi_n)\hat f(\s|\bxi_n) }{\fhat(\s) }\right|}_{\rm (b)}. 
\end{align*}
We analyze each term separately:
\begin{enumerate}[(a):]
\item Using Lemma \ref{lem:ratio_bound}, we have:
\begin{align*}
&\;\left|\frac{\int \ell(\x,\xii) f(\s,\xii) {\rm d}\xii}{f(\s)}-\frac{\int \ell(\x,\xii) \hat f(\s|\xii) f(\xii){\rm d}\xii}{\fhat(\s)}\right|\\
\leq & \;\frac{\left|\int \ell(\x,\xii) f(\s,\xii) {\rm d}\xii-\int \ell(\x,\xii) \hat f(\s|\xii) f(\xii){\rm d}\xii\right|}{f(\s)} + \frac{|f(\s)-\fhat(\s)|\left|\int \ell(\x,\xii) \hat f(\s|\xii) f(\xii){\rm d}\xii\right|}{f(\s)\fhat(\s)}\\
\leq & \;\underbrace{\frac{\lbar\int\left| f(\s,\xii) -\hat f(\s|\xii) f(\xii)\right|{\rm d}\xii}{f(\s)}}_{\textup{(i)}}  + \underbrace{\frac{|f(\s)-\fhat(\s)|\cdot\lbar \int  \hat f(\s|\xii) f(\xii){\rm d}\xii}{f(\s)\fhat(\s)}}_{\textup{(ii)}}.
\end{align*}
We now bound each term:
\begin{enumerate}[(i)]
\item  We use Proposition \ref{prop:gmmpdfbound} to write
\begin{align*}
|f(\s,\xii) -\hat f(\s|\xii) f(\xii)|
\leq&\;\;|f(\s,\xii) -\fhat(\s,\xii)|+\hat f(\s|\xii) |\fhat(\xii)-f(\xii)|\\
\leq &\;\;f(\s,\xii)\left(C_{Q+R}\|(\s,\xii)\|^2 + C'_{Q+R}\|(\s,\xii)\| + C''_{Q+R}\right)\\
&\;\;+\fu f(\xii)\left(C_{R}\|\xii\|^2 + C'_{R}\|\xii\| + C''_{R}\right),
\end{align*} 
Integrating and applying the assumptions yields:
\begin{align*}
\textup{(i)}
\leq &\;\;\lbar\int f(\xii|\s) \left(C_{Q+R}\|(\s,\xii)\|^2 + C'_{Q+R}\|(\s,\xii)\| + C''_{Q+R}\right){\rm d}\xii\\
&\;\;+\frac{\lbar\fu}{\fl} \int f(\xii)\left(C_{R}\|\xii\|^2 + C'_{R}\|\xii\| + C''_{R}\right) {\rm d}\xii\\
\leq &\;\;\lbar\EE_{\Ms}\left[C_{Q+R}\|(\s,\xii)\|^2 + C'_{Q+R}\|(\s,\xii)\| + C''_{Q+R}\right]\\
&\;\;+\frac{\lbar\fu}{\fl} \EE_{\M}\left[C_{R}\|\xii\|^2 + C'_{R}\|\xii\| + C''_{R}\right].
\end{align*}
Now note:
\begin{align*}
 \EE_{\Ms} \left[\| (\s,\txi)\|^2 \right] &\leq  \|\s\|^2  + \EE_{\Ms} \left[ \|\txi\|^2\right]\\
 &=  \|\s\|^2 + \sum_{k\in[K]}p^k_{\xii|\s}\EE_{\N(\bm{\mu}_{\xii|\s}^k,\bSigma_{\xii|\s}^k)} \left[ \|\txi\|^2\right]\\
 & =\|\s\|^2 + \sum_{k\in[K]}p^k_{\xii|\s}\EE_{\N(\bm{\mu}_{\xii|\s}^k,\bSigma_{\xii|\s}^k)} \left[ \tr(\txi\txi^\top)\right]\\
  & =\|\s\|^2 + \sum_{k\in[K]}p^k_{\xii|\s}\left(\tr\left(\bSigma_{\xii|\s}^k\right)+ \|\bm{\mu}_{\xii|\s}^k\|^2\right).
\end{align*}
In addition, since\begin{align*}
\tr\left(\bm{\Sigma}^k_{\xii|\s}\right)&=\tr\left(  \bm{\Sigma}_{\bxi\bxi}^k- \bm{\Sigma}_{\xii\s}^k( \bm{\Sigma}_{\s\s}^{k})^{-1} \bm{\Sigma}_{\s\xii}^k\right)\leq \beta
\end{align*}
and
\begin{align*}
\|\bm{\mu}_{\xii|\s}^k\|^2  =& \;\; \|\bm{\mu}^k_{\xii} + \bm{\Sigma}_{\xii\s}^k( \bm{\Sigma}_{\s\s}^{k})^{-1}(\bm{s}-\bm{\mu}_\s^k)\|^2\\
\leq &\;\;\left(\gamma+\frac{\beta}{\alpha}(\|\s\|+\gamma)\right)^2,
\end{align*}
we obtain
\begin{align*}
 \EE_{\Ms} \left[\| (\s,\txi)\|^2 \right] &\leq  \|\s\|^2  + \beta + \left(\gamma+\frac{\beta}{\alpha}(\|\s\|+\gamma)\right)^2. 
\end{align*}
Next, by Jensen's inequality
\begin{align*}
 \EE_{\Ms} \left[\| (\s,\txi)\| \right] &\leq \sqrt{\EE_{\Ms} \left[\| (\s,\txi)\|^2 \right]}\\
  & =\sqrt{ \|\s\|^2  + \beta + \left(\gamma+\frac{\beta}{\alpha}(\|\s\|+\gamma)\right)^2}.
\end{align*}
Hence,
\begin{align*}
&\;\;\EE_{\Ms}\left[C_{Q+R}\|(\s,\txi)\|^2 + C'_{Q+R}\|(\s,\txi)\| + C''_{Q+R}\right]\\
\leq &\;\;C_{Q+R}\left( \|\s\|^2  + \beta + \left(\gamma+\frac{\beta}{\alpha}(\|\s\|+\gamma)\right)^2\right)+C'_{Q+R}\sqrt{ \|\s\|^2  + \beta + \left(\gamma+\frac{\beta}{\alpha}(\|\s\|+\gamma)\right)^2} + C''_{Q+R}.
\end{align*}

For the marginal expectation:
\begin{align*}
&\;\;\EE_{\M}\left[C_{R}\|\xii\|^2 + C'_{R}\|\xii\| + C''_{R}\right]\\
&\;\;\leq C_R \left(\sum_{k\in[K]}p^k\left(\tr(\bSigma^k)+\|\bmu\|^2\right)\right)+C'_R\sqrt{\sum_{k\in[K]}p^k\left(\tr(\bSigma^k)+\|\bmu\|^2\right)}+C''_R\\
&\;\; \leq C_R \left(\beta+\gamma^2\right)+C'_R\sqrt{\beta+\gamma^2}+C''_R.
\end{align*}

\item Again using Proposition  \ref{prop:gmmpdfbound}, we obtain 
\begin{align*}
&\frac{\left|f(\s)-\fhat(\s)\right|}{f(\s)}\leq C_{Q}\| \s\|^2+C'_{Q}\| \s\| + C''_{Q},
\end{align*}
and so
\begin{align*}
\frac{|f(\s)-\fhat(\s)|\cdot\lbar \int  \hat f(\s|\xii) f(\xii){\rm d}\xii}{f(\s)\fhat(\s)}
\leq  \frac{\lbar\fbar}{\fl}\left(C_{Q}\| \s\|^2+C'_{Q}\| \s\| + C''_{Q}\right). 
\end{align*}
Combining both bounds gives:
\begin{equation*}
\begin{array}{rl}
\;\;{\rm(a)}
\leq&\;\; \lbar\Bigg[\;C_{Q+R}\left( \|\s\|^2  + \beta + \left(\gamma+\frac{\beta}{\alpha}(\|\s\|+\gamma)\right)^2\right)\\
&\;\;\qquad+C'_{Q+R}\sqrt{ \|\s\|^2  + \beta + \left(\gamma+\frac{\beta}{\alpha}(\|\s\|+\gamma)\right)^2} + C''_{Q+R}\Bigg] \\
&\;\; + \frac{\lbar\fbar}{\fl} \left(C_R \left(\beta+\gamma^2\right)+C'_R\sqrt{\beta+\gamma^2}+C''_R+C_{Q}\| \s\|^2+C'_{Q}\| \s\| + C''_{Q}\right). 
\end{array}
\end{equation*}

\end{enumerate}
\item By Hoeffding's inequality, we have with probability at least $1-\delta$:
\begin{align*}
&\left|\int \ell(\x,\xii) \hat f(\s|\xii) f(\xii){\rm d}\xii-\frac{1}{N}\sum_{n\in[N]} \ell(\x,\bxi_n)\hat f(\s|\bxi_n)\right|\leq \lbar\fbar\sqrt{\frac{1}{2N}\log\left(\frac{4}{\delta}\right)}.
\end{align*}
Thus, 
\begin{align*}
&\displaystyle\;\;\;\;\left|\frac{\int \ell(\x,\xii) \hat f(\s|\xii) f(\xii){\rm d}\xii}{\hat f(\s)}- \frac{\frac{1}{N}\sum_{n\in[N]}\ell(\x,\bxi_n)\hat f(\s|\bxi_n) }{\fhat(\s) }\right|\leq \frac{\lbar\fbar}{\fl}\sqrt{\frac{1}{2N}\log\left(\frac{4}{\delta}\right)}.
\end{align*}

\end{enumerate}
Finally, combining (a), (b), and applying the  union bound with Assumption (D), we have
\begin{equation*}
\begin{array}{rl}
&\displaystyle\;\;\;\;\left|\EE_{\M}\left[\ell(\x,\txi)\big|\ts=\s\right]-\frac{1}{\fhat(\s)N}\sum_{n\in[N]}\ell(\x,\bxi_n)\hat f(\s|\bxi_n)\right|\\
\leq & \lbar\left[\;C_{Q+R}\left( \|\s\|^2  + \beta + \left(\gamma+\frac{\beta}{\alpha}(\|\s\|+\gamma)\right)^2\right)+C'_{Q+R}\sqrt{ \|\s\|^2  + \beta + \left(\gamma+\frac{\beta}{\alpha}(\|\s\|+\gamma)\right)^2} + C''_{Q+R}\right] \\
&\;\; + \frac{\lbar\fbar}{\fl} \left(C_R \left(\beta+\gamma^2\right)+C'_R\sqrt{\beta+\gamma^2}+C''_R+C_{Q}\| \s\|^2+C'_{Q}\| \s\| + C''_{Q}\right)\\
&\quad+\frac{\lbar\fbar}{\fl}\sqrt{\frac{1}{2N}\log\left(\frac{4}{\delta}\right)}.
\end{array}
\end{equation*}
holds with probability at least 1-2$\delta$. Setting $2\delta \rightarrow \delta$ completes the proof.  \qed

\end{proof}

The proof of Theorem~\ref{thm:out_of_sample} relies on the following preliminary result.
\begin{lemma}
\label{lem:gmm_proxy}
Let $\txi \in \RR^d$ be a random vector drawn from $\M=\sum_{k\in[K]}p^{k}\N(\bmu^{k},\bSigma^{k})$. Then, the centered vector $\txi - \EE[\txi]$ is sub-Gaussian, and its squared sub-Gaussian parameter $\sigma^2$ satisfies:
\begin{equation}
    \sigma^2 \le \max_{j \in [K]} \left\lVert\bm\Sigma^j\right\rVert_2 + \max_{k \in [K]} \left\lVert\bm\mu^k - \bar{\bm\mu}\right\rVert_2^2
\end{equation}
where $\bar{\bm\mu} = \EE[\txi] = \sum_{k=1}^K p^k \bm\mu^k$ is the global mean of the mixture. 
\end{lemma}

\begin{proof}{Proof of Lemma~\ref{lem:gmm_proxy}}
To establish the sub-Gaussian property, we bound the moment generating function (MGF) of the projection $\bm u^\top(\txi - \bar{\bm\mu})$ for any unit vector $\bm u \in \RR^d$. By the law of total expectation, the MGF can be written as a convex combination of the component MGFs:
\begin{align*}
    \EE\left[\exp\left(t \bm u^\top (\txi - \bar{\bm\mu})\right)\right] &= \sum_{k=1}^K p^k \EE\left[\exp\left(t \bm u^\top (\txi^k - \bar{\bm\mu})\right)\right] \\
    &= \sum_{k=1}^K p^k \exp\left(t \bm u^\top (\bm\mu^k - \bar{\bm\mu}) + \frac{t^2}{2}(\bm u^\top \bm\Sigma^k \bm u)\right),
\end{align*}
where $\txi^k \sim \mathcal{N}(\bm\mu^k, \bm\Sigma^k)$. We can bound the quadratic form $\bm u^\top \bm\Sigma^k \bm u \le \left\lVert\bm\Sigma^k\right\rVert_2$, where $\left\lVert\cdot\right\rVert_2$ is the spectral norm. This allows us to factor out a uniform bound for the variance terms:
\begin{equation}
\label{eq:proof_mid_step}
    \EE\left[\exp\left(t \bm u^\top (\txi - \bar{\bm\mu})\right)\right] \le \exp\left(\frac{t^2}{2} \max_{j \in [K]} \left\lVert\bm\Sigma^j\right\rVert_2\right) \left( \sum_{k=1}^K p^k \exp\left(t \bm u^\top (\bm\mu^k - \bar{\bm\mu})\right) \right).
\end{equation}
The remaining sum is the MGF of $\bm u^\top(\tilde{\bm \zeta} - \bar{\bm\mu})$, where $\tilde{\bm \zeta}$ is an auxiliary random vector that takes value $\bm\mu^k$ with probability $p^k$. 

By Hoeffding's lemma, the MGF of a zero-mean random variable bounded in $[-R, R]$ is bounded by $\exp(t^2 R^2 / 2)$. Therefore, we have:
\begin{align*}
    \sum_{k=1}^K p^k \exp\left(t \bm u^\top (\bm\mu^k - \bar{\bm\mu})\right)
    &\le \exp\left(\frac{t^2}{2} R^2\right) \\
    &= \exp\left(\frac{t^2}{2} \max_{k \in [K]} \left\lVert\bm\mu^k - \bar{\bm\mu}\right\rVert_2^2\right).
\end{align*}

Substituting this back into inequality \eqref{eq:proof_mid_step} yields:
\begin{align*}
    \EE\left[\exp\left(t \bm u^\top (\txi - \bar{\bm\mu})\right)\right] &\le \exp\left(\frac{t^2}{2} \max_{j\in [K]} \left\lVert\bm\Sigma^j\right\rVert_2\right) \exp\left(\frac{t^2}{2} \max_{k\in [K]} \left\lVert\bm\mu^k - \bar{\bm\mu}\right\rVert_2^2\right) \\
    &= \exp\left( \frac{t^2}{2} \left( \max_{j\in [K]} \left\lVert\bm\Sigma^j\right\rVert_2 + \max_{k\in [K]} \left\lVert\bm\mu^k - \bar{\bm\mu}\right\rVert_2^2 \right) \right).
\end{align*}
Since this bound holds for any unit vector $\bm u$, the result follows from the definition of a sub-Gaussian random vector. \qed

\end{proof}

We now proceed to show the derivation of Theorem~\ref{thm:out_of_sample}.

\begin{proof}{Proof of Theorem \ref{thm:out_of_sample}}
We first establish a high probability bound on the 2-norm of the sample points $\{\bxi_{t,n}\}_{t\in[T],n\in[N]}$. Since $\bxi_{t,n}\in\RR^R$  sampled from a sub-Gaussian distribution with variance proxy $ \sigma^2 \le \max_{j \in [K]} \left\lVert\bm\Sigma^j\right\rVert_2 + \sum_{k=1}^K p^k \left\lVert\bm\mu^k - \bar{\bm\mu}\right\rVert_2^2$ , then by \cite[Theorem 1.19]{rigollet201518}, we have with probability at least $1-\delta$,
\begin{equation}
\label{eq:bound_norm_xi}
\|\bxi_{t,n}\|\leq \bar\xi\coloneqq 4\sigma\sqrt{R}+2\sigma\sqrt{2\log\tfrac{NT}{\delta} }+\gamma \quad\forall n\in[N]\;\;\forall t\in[T].
\end{equation}
Next, we ensure that the bound in \eqref{eq:fix_x_bound} holds uniformly for all $\x\in\X_t(\cdot)$.
To simplify the exposition,  we define 
\begin{equation*}
\begin{array}{rl}
\tau\coloneqq& \lbar\left[\;C_{Q+R}\left( \bar\xi^2  + \beta + \left(\gamma+\frac{\beta}{\alpha}(\bar\xi+\gamma)\right)^2\right)+C'_{Q+R}\sqrt{ \bar\xi^2  + \beta + \left(\gamma+\frac{\beta}{\alpha}(\bar\xi+\gamma)\right)^2} + C''_{Q+R}\right] \\
&\;\; + \frac{\lbar\fbar}{\fl} \left(C_R \left(\beta+\gamma^2\right)+C'_R\sqrt{\beta+\gamma^2}+C''_R+C_{Q}\bar\xi^2+C'_{Q}\bar\xi + C''_{Q}\right)
\end{array}
\end{equation*}
Note that by construction, $\tau$ upper bounds the error expression in \eqref{eq:fix_x_bound} for all sample points satisfying \eqref{eq:bound_norm_xi}.
Thus, by \eqref{eq:fix_x_bound}, we have 
\begin{align*}
\label{eq:fix_x_bound_simple}
&\nonumber\displaystyle\;\;\;\;\left|\EE\left[V_{t+1}(\x,\txi_{t+1})\big|\txi_t=\bxi_t\right]- \hat\EE\left[V_{t+1}(\x,\txi_{t+1})\big|\txi_t=\bxi_t\right]\right|
\leq  \tau+\frac{\lbar\fbar^2}{\fl^2}\sqrt{\frac{2}{N}\log\left(\frac{8}{\delta}\right)}
\end{align*}
with probability at least $1-\delta$. 
Using a covering number argument, we obtain the uniform bound
\begin{align*}
&\nonumber\displaystyle\;\;\;\;\left|\EE\left[V_{t+1}(\x,\txi_{t+1})\big|\txi_t=\bxi_t\right]- \hat\EE\left[V_{t+1}(\x,\txi_{t+1})\big|\txi_t=\bxi_t\right]\right|
\leq  \tau+\frac{\lbar\fbar^2}{\fl^2}\sqrt{\frac{2}{N}\log\left(\frac{\mathcal O(\overline D/\eta)^R}{\delta}\right)}+2L\eta 
\end{align*}
for all $\x_t\in\X_t(\x_{t-1},\bxi_t)$ with probability at least $1-\delta$. Following the recursive induction argument in the proof of \citep[Theorem 1]{park2022data}, we propagate the stagewise approximation errors backward from $t = T$ to $t = 1$, yielding:
\begin{align*}
    &\bm \ell(\widehat{\bm x}_1^N,\xii_1)+ {\mathcal V}_{2}(\widehat{\bm x}_1^N, \bm \xi_1)
    - \left(\min_{\bm x_1 \in \X_1(\x_0,\bm \xi_1) }
    \ell(\x_1,\xii_1) + {\mathcal V}_{2}({\bm x}_1, \bm \xi_1)\right)\\
    &\quad\quad      \leq 2(T-1)\tau + 4(T-1)L\eta + 2\sum_{t=2}^{T}
        \frac{\lbar\fbar^2}{\fl^2}
         \sqrt{\frac{2}{N} \log \left(\frac{O(1) N^{t-2} (\overline D / \eta)^{R(t-1)} }{ \delta }\right) }
\end{align*}
with probability at least $1-(T-1)\delta$.
Finally, applying the union bound with \eqref{eq:bound_norm_xi} and Assumption (D), we conclude that the result holds with probability at least $1 - (T+1)\delta$. Finally, setting $(T+1)\delta \rightarrow \delta$ completes the proof.
\qed
\end{proof}

\section{Details of Experiments}\label{sec:experiment_appendix}
\subsection{Inventory Management}\label{subsec:inventory_appendix}

For the contextual newsvendor problem, the distributionally robust optimization model derived from Proposition~\ref{prop:dro_reformulation} is formulated as the following SOCP:
\begin{equation*} 
\begin{aligned}
\inf\;\;&\displaystyle\lambda\varepsilon^2 + \frac{1}{M}\sum_{m\in[M]} \gamma_m \\
\st&\displaystyle\lambda\in\RR_+, \bm \gamma\in\RR^Q,\theta\in\RR_+, q\in\RR_+ \\
&\displaystyle \left\|\begin{bmatrix}2\lambda\xi_m+\theta-h \\ 
\displaystyle\lambda\xi_m^2-hq+\gamma_m-\lambda\end{bmatrix} \right\|_2 \leq \lambda\xi_m^2-hq+\gamma_m+\lambda & \displaystyle\forall m\in[M] \\
&\displaystyle \left\|\begin{bmatrix}2\lambda\xi_m+\theta+b \\
\lambda\xi_m^2+bq+\gamma_m-\lambda\end{bmatrix} \right\|_2 \leq \lambda\xi_m^2+bq+\gamma_m+\lambda & \forall m\in[M] \\
& \displaystyle\lambda\xi_m^2-hq+\gamma_m \geq 0 &  \displaystyle\forall m\in[M] \\
& \displaystyle\lambda\xi_m^2+bq+\gamma_m \geq 0 & \displaystyle\forall m\in[M]
\end{aligned}
\end{equation*}

The implementation details of the GMM-based methods are as follows. We first standardize the joint observations, consisting of the side information and realized demand, to improve numerical stability. We then implement the plain-vanilla GMM method. Specifically, we fit a GMM to the standardized joint samples $(\bm s,\xi)$. The number of mixture components is selected from the candidate set $\{1,2,3\}$ using the Akaike Information Criterion (AIC)~\citep{akaike2025akaike}, a common metric for balancing goodness-of-fit and model complexity~\citep{CHAKRABARTI2011583}. After obtaining the GMM parameters, we compute the conditional distribution of $\tilde\xi$ given the test covariate $\bm s$ using Proposition~\ref{prop:conditional_param} and generate $M=100$ Monte Carlo samples from this conditional distribution. These samples are then used in the Wasserstein-robust newsvendor problem to obtain the optimal decision.

We next implement the GMM-NF method. We use the GM distribution fitted by the plain-vanilla GMM method as the base distribution and train a normalizing flow $T_\btheta$ to map this base distribution to the standardized data distribution. We use a lightweight Autoregressive Rational Quadratic Spline (ARQS) architecture for $T_\btheta$, with 32 hidden nodes, a single hidden layer, one flow block, and 8 spline bins. Here, one flow block refers to a single invertible autoregressive transformation module in the normalizing-flow composition, and the spline bins specify the number of intervals used to parameterize each monotone rational-quadratic spline transformation~\citep{papamakarios2021normalizing}. To avoid overfitting, we split the training dataset into an 80\% subtraining set and a 20\% validation set and use early stopping. The training process is halted if there is no improvement in the validation loss for 50 consecutive epochs, with the maximum number of training epochs capped at 500. Demand samples are generated by conditioning on the observed side information and then mapping the conditional samples back to the original demand scale through the learned flow.

The Wasserstein radius $\epsilon$ is selected by cross-validation over the grid $\{a\times10^b: a\in\{1,5,9\},\, b\in\{-2,-1,0,1\}\}$. To reduce the computational burden in the largest-sample configuration, we set $\epsilon=0$ rather than cross-validating $\epsilon$ when $N=400$. All methods are evaluated using the same test covariate and the same out-of-sample Monte Carlo seed within each trial to ensure a fair comparison.

We also report the running time of different methods under varying contextual dimensions $Q$ and training sample sizes $N$ in Table~\ref{tab:runtime_vs_N}. For each configuration, we repeat the experiment over 100 independent trials and report the average runtime together with the 10th and 90th percentiles. To ensure a fair comparison of the computational cost of the decision-generation pipeline, we exclude all cross-validation procedures from the timing evaluation and use fixed hyperparameter settings for all methods. In particular, the number of GMM components and the Wasserstein radius are fixed across trials. 

\revise{
\begin{table}[htbp]
\centering
\scriptsize
\resizebox{\textwidth}{!}{%
\begin{tabular}{clrrrrrrrrr}
\toprule
\multirow{2}{*}{$N$} & \multirow{2}{*}{Method}
& \multicolumn{3}{c}{$Q=1$}
& \multicolumn{3}{c}{$Q=5$}
& \multicolumn{3}{c}{$Q=20$} \\
\cmidrule(lr){3-5} \cmidrule(lr){6-8} \cmidrule(lr){9-11}
& & mean & $P_{10\text{-th}}$ & $P_{90\text{-th}}$
  & mean & $P_{10\text{-th}}$ & $P_{90\text{-th}}$
  & mean & $P_{10\text{-th}}$ & $P_{90\text{-th}}$ \\
\midrule
\multirow{6}{*}{50}
& GMM-NF     & 0.757 & 0.679 & 0.930 & 0.722 & 0.685 & 0.804 & 0.701 & 0.686 & 0.731 \\
& GMM        & 0.522 & 0.383 & 0.562 & 0.533 & 0.391 & 0.571 & 0.534 & 0.389 & 0.573 \\
& ResDRO     & 0.150 & 0.106 & 0.265 & 0.170 & 0.106 & 0.267 & 0.172 & 0.106 & 0.271 \\
& RNW        & 0.122 & 0.076 & 0.223 & 0.099 & 0.076 & 0.221 & 0.101 & 0.076 & 0.226 \\
& LDR        & 0.005 & 0.004 & 0.005 & 0.005 & 0.005 & 0.005 & 0.005 & 0.005 & 0.006 \\
& NN-DR      & 0.365 & 0.331 & 0.395 & 0.358 & 0.336 & 0.405 & 0.348 & 0.298 & 0.405 \\
\midrule
\multirow{6}{*}{100}
& GMM-NF     & 0.995 & 0.672 & 1.519 & 0.890 & 0.719 & 1.136 & 0.874 & 0.797 & 1.004 \\
& GMM        & 0.547 & 0.541 & 0.565 & 0.558 & 0.553 & 0.573 & 0.576 & 0.561 & 0.591 \\
& ResDRO     & 0.279 & 0.213 & 0.400 & 0.277 & 0.215 & 0.403 & 0.293 & 0.217 & 0.410 \\
& RNW        & 0.285 & 0.153 & 0.307 & 0.286 & 0.152 & 0.308 & 0.246 & 0.149 & 0.310 \\
& LDR        & 0.005 & 0.005 & 0.005 & 0.005 & 0.005 & 0.006 & 0.007 & 0.007 & 0.007 \\
& NN-DR      & 0.542 & 0.451 & 0.666 & 0.529 & 0.382 & 0.649 & 0.482 & 0.266 & 0.652 \\
\midrule
\multirow{6}{*}{200}
& GMM-NF     & 1.824 & 0.907 & 2.803 & 1.152 & 0.931 & 1.384 & 1.069 & 0.907 & 1.210 \\
& GMM        & 0.576 & 0.564 & 0.606 & 0.581 & 0.574 & 0.610 & 0.584 & 0.579 & 0.623 \\
& ResDRO     & 0.718 & 0.630 & 0.852 & 0.734 & 0.642 & 0.873 & 0.725 & 0.645 & 0.875 \\
& RNW        & 0.411 & 0.301 & 0.468 & 0.411 & 0.300 & 0.471 & 0.416 & 0.301 & 0.478 \\
& LDR        & 0.006 & 0.006 & 0.006 & 0.007 & 0.007 & 0.007 & 0.010 & 0.009 & 0.010 \\
& NN-DR      & 0.638 & 0.461 & 0.806 & 0.599 & 0.398 & 0.820 & 0.619 & 0.362 & 0.933 \\
\midrule
\multirow{6}{*}{400}
& GMM-NF     & 3.118 & 1.823 & 4.631 & 1.830 & 1.447 & 2.239 & 1.634 & 1.534 & 1.738 \\
& GMM        & 0.454 & 0.386 & 0.640 & 0.601 & 0.402 & 0.666 & 0.641 & 0.645 & 0.666 \\
& ResDRO     & 1.487 & 1.418 & 1.519 & 1.522 & 1.436 & 1.786 & 1.470 & 1.440 & 1.507 \\
& RNW        & 0.929 & 0.803 & 0.982 & 0.850 & 0.800 & 0.990 & 0.825 & 0.801 & 0.860 \\
& LDR        & 0.007 & 0.007 & 0.008 & 0.010 & 0.009 & 0.010 & 0.015 & 0.015 & 0.016 \\
& NN-DR      & 0.953 & 0.568 & 1.394 & 0.822 & 0.468 & 1.399 & 0.960 & 0.479 & 1.723 \\
\bottomrule
\end{tabular}
}
\caption{Runtime comparison of different approaches across varying contextual dimensions ($Q$) and training sample sizes ($N$). Each entry reports the average runtime in seconds, together with the 10th and 90th percentiles.}
\label{tab:runtime_vs_N}
\end{table}
}

Table~\ref{tab:runtime_vs_N} shows that the GMM method is computationally efficient and stable as the sample size and contextual dimension vary. The GMM-NF method requires more time as the sample size increases, primarily because of the additional training cost of the normalizing flow. Nevertheless, its runtime remains manageable across all reported configurations. We also observe a slight decrease in the running time of GMM-NF as the covariate dimension grows. This occurs because the early stopping criterion is typically triggered earlier in high-dimensional settings, mitigating the risk of overfitting. 

\subsection{Portfolio Optimization}\label{subsec:portfolio_appendix}

We next describe the implementation details for the portfolio experiment. By applying the 2-Wasserstein robust formulation from Proposition~\ref{prop:dro_reformulation}, our model can be formulated as the following SOCP:
\begin{equation*}
\begin{aligned}
\inf \quad & 
\displaystyle\lambda \varepsilon^2 + \frac{1}{M} \sum_{m\in[M]} \bm \gamma_m \\
\st \quad & \displaystyle\lambda \in \RR_+,\quad \x\in\R^D_+, \quad \bm\gamma \in \RR^{M \times D},\quad \beta \in \RR \\
\displaystyle& \left\| \begin{bmatrix} 2\lambda - \left(\frac{1}{\tau} + \eta\right)\bm x \\
\displaystyle\lambda \bxi_m^\top \bxi_m + \bm \gamma_m + \frac{\beta}{\tau} - \beta - \lambda\end{bmatrix}\right\| \leq \lambda \bxi_m^\top \bxi_m + \bm \gamma_m + \frac{\beta}{\tau} - \beta + \lambda & \forall m \in [M] \\
\displaystyle& \left\| \begin{bmatrix} 2\lambda \bxi_m - \eta \bm x \\
\displaystyle\lambda \bxi_m^\top \bxi_m + \bm \gamma_m - \beta - \lambda\end{bmatrix}\right\|  \leq \lambda \bxi_m^\top \bxi_m + \bm \gamma_m - \beta + \lambda &\forall m \in [M] \\
\displaystyle& \lambda \bxi_m^\top \bxi_m + \bm \gamma_m \geq -\frac{\beta}{\tau}+\beta &\forall m \in [M] \\
\displaystyle& \lambda \bxi_m^\top \bxi_m + \bm \gamma_m \geq \beta &\forall m \in [M] \\
\displaystyle& \sum_{d\in[D]} x_d = 1
\end{aligned}
\end{equation*} 

We first obtain the raw side information and asset return data from publicly available datasets and standardize each block separately to ensure numerical stability. Next, we fit a GMM to the concatenated standardized joint samples. The number of mixture components is selected by AIC from the candidate set $\{1,2,3\}$. We then adopt the fitted GMM as the base distribution and train an ARQS normalizing flow~\citep{papamakarios2021normalizing} to learn the transformation from the data distribution to this base distribution. The flow architecture $T_\btheta$ is fixed across all runs, with 64 hidden nodes, 2 hidden layers, 2 spline blocks, and 8 bins. To prevent overfitting, we split the training dataset into an 80\% subtraining set and a 20\% validation set and use early stopping. The training process is halted if there is no improvement in the validation loss for 30 consecutive epochs, with the maximum number of training epochs capped at 500. Finally, 100 conditional samples of future returns are generated by conditioning on the observed side information within the latent space and then mapped back to the original space through the learned flow. 

The hyperparameter tuning procedure is carried out on a quarterly rolling basis. For each out-of-sample quarter, we use the immediately preceding calendar month as the validation period and the prior two years of data as the estimation window. We tune the Wasserstein radius over $\{0.01, 0.05, 0.1, 0.5, 1.0, 1.5, 2.0\}$. Once the hyperparameters are chosen, the corresponding model is re-estimated on the full two-year training window.

\subsection{Wind Energy Multi-Stage Optimization}\label{subsec:Wind_energy_appendix}

We adopt the multistage wind-energy storage framework of~\citet{park2022data}. At each decision stage, the decision maker chooses a 24-hour day-ahead commitment vector, and the next-day power realization determines the recourse decisions and the next storage state. Let $\bm u_t\in\mathbb R^{24}$ denote the 24-hour day-ahead commitment vector at stage \(t\). We consider three storage units indexed by \(\ell\in[3]\). Storage unit \(\ell\) has capacity \(\bar\zeta^\ell\), leakage factor \(\rho^\ell\), charging efficiency \(\rho_c^\ell\), and discharging efficiency \(\rho_d^\ell\). In our implementation, the storage capacity vector is \(\bar{\bm\zeta}=(5000,2000,1000)\). The leakage, charging, and discharging coefficients are \(\bm\rho=(0.98,0.99,0.995)\), \(\bm\rho_c=(0.8,0.9,1.0)\), and \(\bm\rho_d=(0.8,0.9,1.0)\), respectively, and the initial storage level is \(\bm \zeta_0=(500,500,500)\). For each child scenario \(n\) and hour \(h\), let \(e_{t,n}^{s,h}\) denote the generated energy used to satisfy the commitment, \(e_{t,n}^{u,h}\) the unmet commitment, and \(e_{t,n}^{d,h}\) the curtailed energy. In addition, \(e_{t,n}^{+,\ell,h}\) and \(e_{t,n}^{-,\ell,h}\) denote the charging and discharging amounts for storage unit \(\ell\), respectively.

The random state is defined as the concatenation of daily power and price profiles,
\(\bxi_t=(\bm w_t,\bm p_t)\in\mathbb R^{48}\), where \(\bm w_t,\bm p_t\in\mathbb R^{24}\). In the stage-\(t\) subproblem at current node \(j\), the price vector \(\bm p_{t}^j\) appears in the objective, while the next-day power profiles \(\{\bm w_{t+1,n}\}_{n=1}^N\) define the child-scenario recourse constraints. We solve the equivalent negative-profit minimization problem.

To construct the GMM likelihood-ratio transition weights, we fit stagewise joint Gaussian mixture models to consecutive power-price state pairs \((\bxi_t,\bxi_{t+1})\), where \(\bxi_t=(\bm w_t,\bm p_t)\). Given nominal transition weights \(\bm q_t^j=(q_{t,1}^j,\ldots,q_{t,N}^j)\), the stage problem has the form
\[
\begin{aligned}
\mathcal V_t(\bm\zeta_t,\bxi_t^j)=\min\quad
&-(\bm p_t^j)^\top \bm u_t
+2\sum_{n=1}^N q_{t,n}^j(\bm p_t^j)^\top \bm e_{t,n}^{u}
+\mathcal R_{\eps}(\bm q_t^j,\bm\alpha_t)\\
\text{s.t.}\quad
& \bm u_t\in\mathbb R_+^{24},\quad \bm\alpha_t\in\mathbb R^N,\\
& \bm e_{t,n}^{s},\bm e_{t,n}^{u},\bm e_{t,n}^{d}\in\mathbb R_+^{24},\quad
\bm e_{t,n}^{+,\ell},\bm e_{t,n}^{-,\ell}\in[0,\bar\zeta^\ell]^{24}
&& \forall \ell\in[3],\; n\in[N],\\
& w_{t+1,n}^{h}=e_{t,n}^{s,h}+\sum_{\ell\in[3]}e_{t,n}^{+,\ell,h}+e_{t,n}^{d,h}
&& \forall h\in[24],\; n\in[N],\\
& u_t^{h}=e_{t,n}^{s,h}+\sum_{\ell\in[3]}e_{t,n}^{-,\ell,h}+e_{t,n}^{u,h}
&& \forall h\in[24],\; n\in[N],\\
& \zeta_{t+1,n}^{\ell,1}=\rho^\ell\zeta_t^\ell+\rho_c^\ell e_{t,n}^{+,\ell,1}
-\frac{1}{\rho_d^\ell}e_{t,n}^{-,\ell,1}
&& \forall \ell\in[3],\; n\in[N],\\
& \zeta_{t+1,n}^{\ell,h}=\rho^\ell \zeta_{t+1,n}^{\ell,h-1}
+\rho_c^\ell e_{t,n}^{+,\ell,h}
-\frac{1}{\rho_d^\ell}e_{t,n}^{-,\ell,h}
&& \forall h=2,\ldots,24,\; \ell\in[3],\; n\in[N],\\
& 0\le \zeta_{t+1,n}^{\ell,h}\le \bar\zeta^\ell
&& \forall h\in[24],\; \ell\in[3],\; n\in[N],\\
& \alpha_{t,n}\ge a_{r,n}+\bm b_{r,n}^{\top}\bm\zeta_{t+1,n}
&& \forall r\in\mathcal K_{t+1,n},\; n\in[N].
\end{aligned}
\]

Here, \(\bm\zeta_{t+1,n}=(\zeta_{t+1,n}^{1,24},\zeta_{t+1,n}^{2,24},\zeta_{t+1,n}^{3,24})\) denotes the end-of-day storage vector under child scenario \(n\), and \(\mathcal K_{t+1,n}\) is the set of SDDP cuts generated so far for stage \(t+1\) at child node \(n\). The variable \(\alpha_{t,n}\) is the epigraph variable for the SDDP approximation of the next-stage cost-to-go function at child node \(n\). The pair \((a_{r,n}, \bm b_{r,n})\) represents the intercept and slope of the \(r\)-th affine cut used to approximate the downstream value function at child node \(n\). 
The term $\mathcal R_{\eps}(\bm q_t^j,\bm\alpha_t)$ is the distributionally robust future-value functional. Instead of using the nominal expectation $\sum_{n=1}^N q_{t,n}^j\alpha_{t,n}$, we use the worst-case expectation over a modified-$\chi^2$ ambiguity set around the nominal transition vector:
\[
\mathcal R_{\eps}(\bm q_t^j,\bm\alpha_t)
=
\sup_{\bm\pi\in\mathcal U_{\eps}(\bm q_t^j)}
\sum_{n=1}^N \pi_n\alpha_{t,n}.
\]

We implement this robust future-value term through its dual reformulation in the SDDP subproblem. For each stage \(t\) and current node \(j\), the modified-\(\chi^2\) ambiguity set around the nominal transition vector \(\bm q_t^j\) is defined as
\[\mathcal U_{\eps}(\bm q_t^j) = \left\{ \sum_{n \in [N]} \pi_n \cdot \delta_{\bxi_{t+1,n}}:  \bm\pi\in\Delta^N, \left(\left(\frac{\pi_n-q_{t,n}^j}{\sqrt{q_{t,n}^j}} \right)_{n\in[N]}, \sqrt{\eps} \right) \in \mathcal C^{N+1} \right\}, \]
where $\delta_{\bxi_{t+1,n}}$ represents the Dirac distribution that assigns a unit mass at the sample point \(\bxi_{t+1,n}\), and \(\mathcal C^{N+1}=\{(\bm y,\tau)\in\mathbb R^N\times\mathbb R:\|\bm y\|_2\le \tau\}\) denotes the second-order cone. 

We select the hyperparameter by cross-validation. For each quarter, we use the first 75\% of the training set as the subtraining set and the remaining 25\% as the validation set. We search for the ambiguity-set radius $\eps$ in the candidate set $\{0.001, 0.005, 0.01, 0.05\}$. For each candidate, we estimate the GMM transition model and solve the SDDP policy using the subtraining set, and then evaluate the resulting policy on the validation set based on average profit. We select the candidate with the best validation performance. Finally, using the selected parameter, we refit the GMM transition model and solve the final SDDP policy using the entire training set.

We also report the average out-of-sample evaluation time for each method in Table~\ref{tab:runtime_sddp}. The N-SDDP method achieves the fastest evaluation times, reflecting the benefit of its neural network warm-start approach. The proposed GMM-SDDP method maintains computational speeds comparable to those of IND-SDDP and NW-SDDP. This indicates that the GMM-based transition model can be incorporated into the SDDP framework without materially compromising online evaluation speed.

\begin{table}[h]
\centering
\begin{tabular}{ l l c }
\toprule
 Data set & Model & Online evaluation time (s) $\downarrow$ \\
\midrule
 \multirow{4}{*}{Ohio} 
 & IND-SDDP & 81.027 \\
 & N-SDDP & \textbf{42.579} \\
 & NW-SDDP & 94.492 \\
 & GMM-SDDP & 100.818 \\
\midrule
 \multirow{4}{*}{North Carolina} 
 & IND-SDDP & 108.896 \\
 & N-SDDP & \textbf{55.447} \\
 & NW-SDDP & 128.661 \\
 & GMM-SDDP & 128.086  \\
\bottomrule
\end{tabular}
\caption{Out-of-sample evaluation time of SDDP-based methods.}
\label{tab:runtime_sddp}
\end{table}

\newpage
\end{APPENDICES}

\end{document}